%% file: ms.tex
\documentclass[12pt]{article}
\input{packages}
\input{theorem}

\input{title}

\input{paddings}

\input{commands}

\newcommand{\Addresses}{{
  \bigskip
  \footnotesize
  \noindent \textsc{Siarhei Finski, Institut Fourier - Université Grenoble Alpes, France.}\par\nopagebreak
  \noindent  \textit{E-mail }: \texttt{finski.siarhei@gmail.com}.
}}

\newenvironment{sciabstract}{}

\begin{document} 

\maketitle

\begin{sciabstract}
  \textbf{Abstract.}
 	We establish an asymptotic relation between the spectrum of the discrete Laplacian associated to discretizations of a half-translation surface with a flat unitary vector bundle and the spectrum of the Friedrichs extension of the Laplacian with von Neumann boundary conditions.
   \par As an interesting byproduct of our study, we obtain Harnack-type estimates on “almost harmonic" discrete functions, defined on the graphs, which approximate a given surface.
   \par The results of this paper will be later used to relate the asymptotic expansion of the number of spanning trees, spanning forests and weighted cycle-rooted spanning forests on the discretizations to the corresponding zeta-regularized determinants.
\end{sciabstract}

\tableofcontents

\section{Introduction}\label{sect_intro}
	Many natural combinatorial invariants of a finite graph $G$ can be expressed through the spectrum $\spec(\laplcomp_{G})$ of the combinatorial Laplacian $\laplcomp_{G}$, defined as the difference between degree and adjacency operators. For example, the number of connected components corresponds to the multiplicity of 0 in $\spec(\laplcomp_{G})$. 
	The largest eigenvalue of $\laplcomp_{G}$ in a regular graph $G$ is twice the degree of the vertices if and only the graph is bipartite.
	The smallest positive eigenvalue of $\spec(\laplcomp_{G})$ is related to the Cheeger isoperimetric constant, which detects “bottlenecks", see Cheeger \cite{CheegConst}.
	\par 
	According to a famous Matrix-tree theorem of Kirchoff, the product of non-zero eigenvalues of $\laplcomp_{G}$ corresponds to the number of marked spanning trees on $G$. This theorem has been generalized by Forman in \cite[Theorem 1]{Forman} to the setting of a unitary line bundle on a graph and by Kenyon in \cite[Theorems 8,9]{KenyonFlatVecBun} to vector bundles of rank $2$ endowed with ${\rm{SL}}_2(\comp)$-connections (we give a precise meaning for those notions in Section \ref{sect_sq_t_s_disc}). 
	They showed that a sum of cycle-rooted spanning forests (cf. \cite[\S 4]{KenyonFlatVecBun} for related definitions), weighted by some function depending on the monodromy of the connection, calculated over the cycle, can be expressed through the determinant of the Laplacian associated with the vector bundle (see (\ref{eq_comb_lapl})).
	\par 
	Now, if instead of considering a single graph $G$ with a vector bundle over it, we consider a family $\Psi_n$, $n \in \nat^*$ of graphs endowed with vector bundles $F_n$, constructed as approximations of a certain surface $\Psi$ with a flat unitary vector bundle $F$, we could ask ourselves how the spectral invariants of $(\Psi_n, F_n)$ would behave, as $n \to \infty$.
	\par 
	Such families of graphs arise naturally in many contexts in mathematics, when one studies some continuous quantity or process by using scaling limits of grid-based approximations.
	\par
\begin{figure}[h]
		\includegraphics[width=0.4\textwidth]{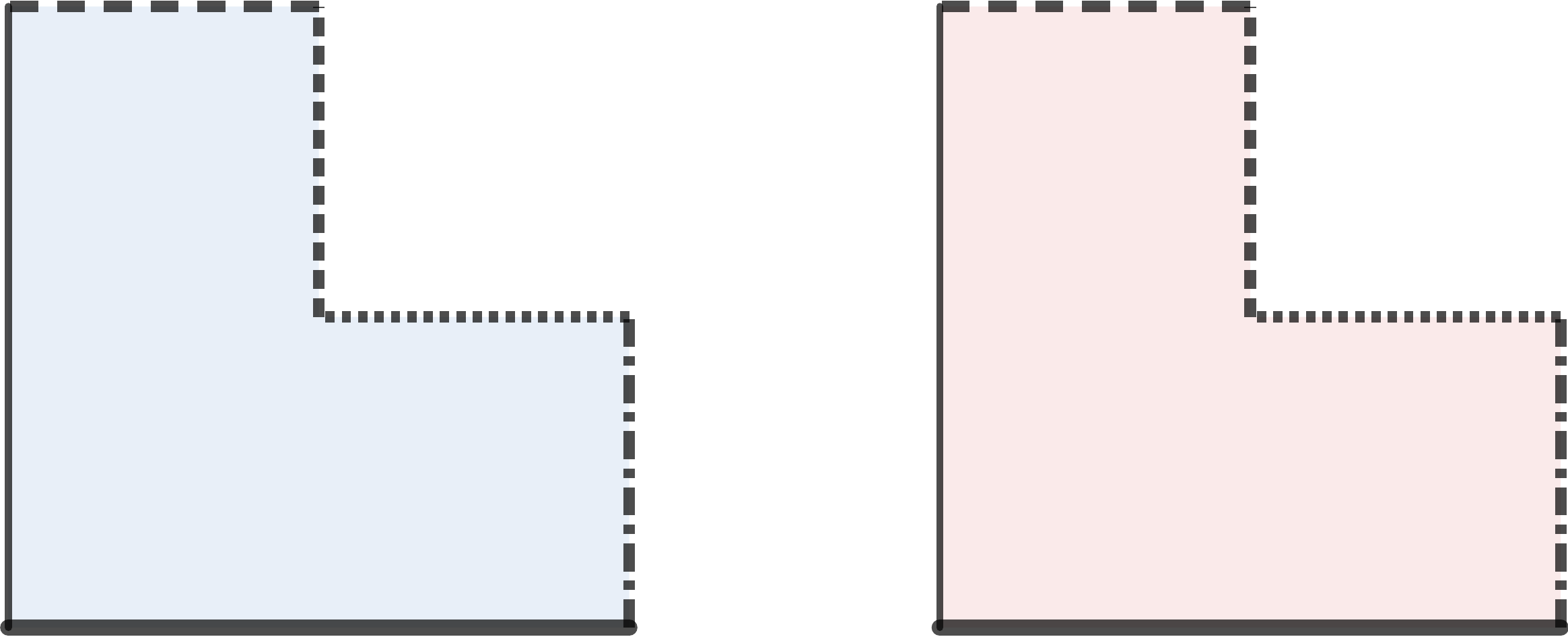}	
		\centering
		\caption{An example of a flat surface. The two domains are glued together along the respective arcs. The resulting surface has 5 conical points of angles $\pi$ and one conical point of angle $3 \pi$. Clearly, from the topological point of view, the resulting surface is a sphere.}
		\label{fig_ex_flat_surf}
	\end{figure}

 The surfaces we consider in this article are called \textit{half-translation} surfaces $(\Psi, g^{T \Psi})$. This means that the metric $g^{T \Psi}$ over $\Psi$ is flat and has conical singularities ${\rm{Con}}(\Psi)$ of angles $k \pi$, $k \in \nat^*$, $k \in \nat^* \setminus \{2\}$, see Figure \ref{fig_ex_flat_surf} for an example. The name half-translation comes from the fact that it has a ramified double cover, which is a translation surface, i.e. it can be obtained by identifying parallel sides in some polygonal domain in $\comp$ (see Section \ref{sect_sq_t_s_disc}).
	\par 
	We endow a \textit{half-translation} surface $(\Psi, g^{T \Psi})$ with piecewise geodesic boundary $\partial \Psi$ with a flat unitary vector bundle $(F, h^F, \nabla^F)$.
	We suppose that $(\Psi, g^{T \Psi})$ can be tiled by euclidean squares of equal sizes.
	We show that for a suitably chosen discretization $\Psi_n$ of $\Psi$, the combinatorial spectral theory of the graph Laplacian $\laplcomp_{\Psi_n}^{F_n}$ on $\Psi_n$, associated with the discretization $(F_n, h^{F_n}, \nabla^{F_n})$ of $(F, h^F, \nabla^F)$, is an approximation, up to a renormalization, of the spectral theory of the Friedrichs extension of the Laplacian $\laplcomp_{\Psi}^{F}$ with von Neumann boundary conditions on $\partial \Psi$.
	\par 
	More precisely, the \textit{Laplacian} $\laplcomp_G^{V}$ associated to a graph $G$ and a vector bundle $V$ over $G$ with a connection $\nabla^V$ is the linear operator, acting on the sections $f$ of $V$ as follows
	\begin{equation}\label{eq_comb_lapl}
		\laplcomp_G^{V} f(v) = \sum_{(v, v') \in E(G)} \big( f(v) - \phi_{v' v} f(v') \big), \qquad v \in V(G),
	\end{equation}
	where $\phi_{v' v}$ are the parallel transports along the edges. 
	In the simplest case when the vector bundle is trivial, we get the standard graph Laplacian.
	Remark that the Laplacian $\laplcomp_G^{V}$ doesn't depend on any Riemannian data, it only depends on the connection $\nabla^F$ of $F$.
	\par We fix a half-translation surface $(\Psi, g^{T \Psi})$ with piecewise geodesic boundary. 
	We also fix a flat unitary vector bundle $(F, h^F, \nabla^F)$ on the compactification 
	\begin{equation}\label{eq_compact_defn}
		\overline{\Psi} := \Psi \cup {\rm{Con}}(\Psi),
	\end{equation}	 
	of $\Psi$ (i.e. we require $(\nabla^F)^2 = 0$ and we suppose that the connection $\nabla^F$ preserves the metric $h^F$). 
	In particular, we suppose that the monodromies around conical points are trivial.
	\par 
	We suppose that $\Psi$ can be tiled completely and without overlaps over subsets of positive Lebesgue measure by euclidean squares of area $1$. In particular, the boundary $\partial \Psi$ is tiled by the boundaries of the tiles, and the angles of the boundary corners  are of the form $\frac{k \pi}{2}$, $k \in \nat^* \setminus \{2\}$.
	Such surfaces are also called \textit{pillowcase covers}, and in case if there is no boundary, they can be characterized as certain ramified coverings of $\mathbb{CP}^1$, see Section \ref{sect_sq_t_s_disc} for  more details.
	\par 
	For example, if $\Psi$ is a torus, then it can be tiled by euclidean squares of the same size if and only if the ratio of its periods is rational.
	If $\Psi$ is a rectangular domain in $\comp$, we are basically requiring that the ratios between the lengths of the sides of $\Psi$ are rational. 
	\begin{figure}[h]
		\includegraphics[width=0.3\textwidth]{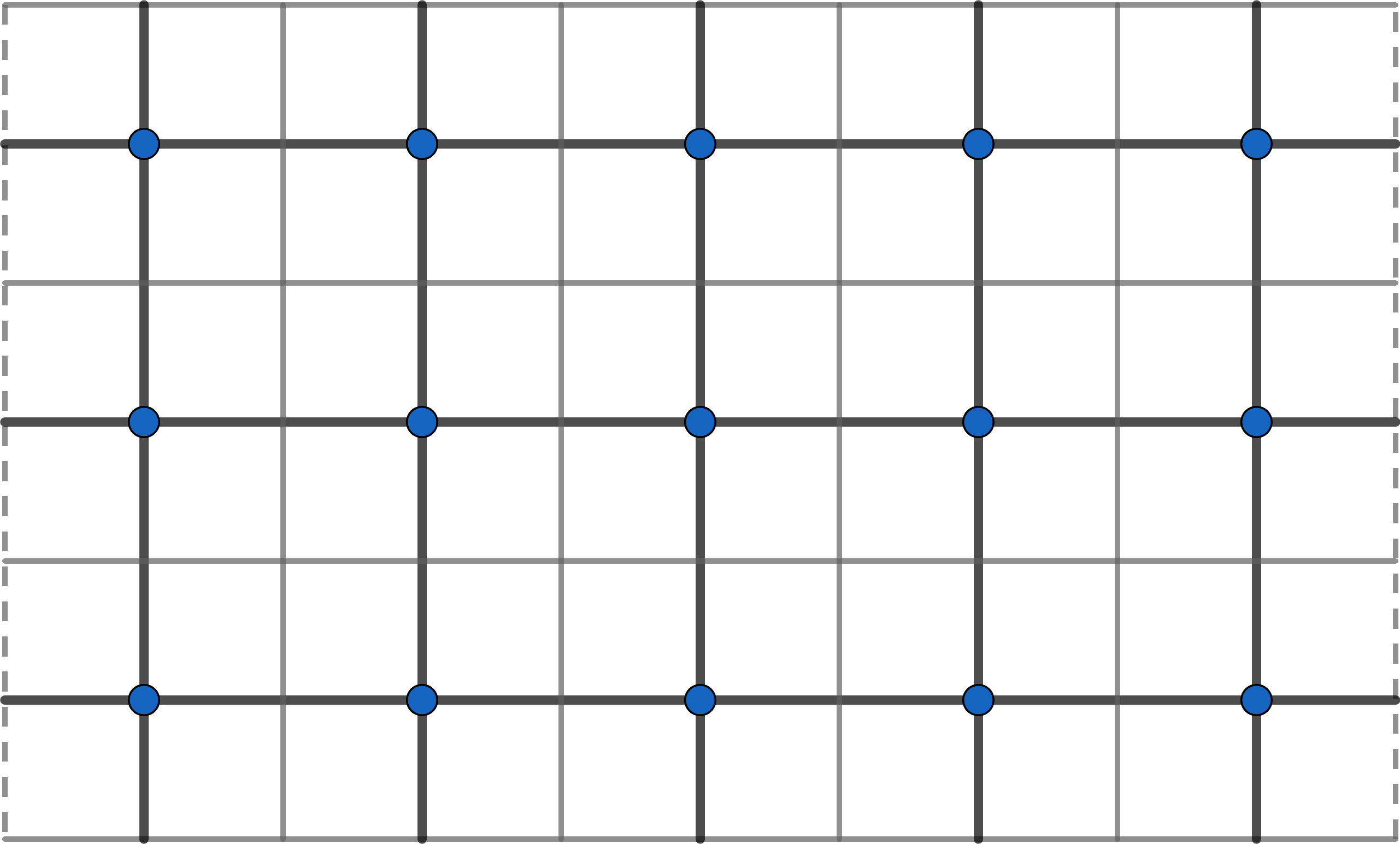}	
		\centering
		\caption{A discretization of torus. The border lines of the same pattern are glued together and the edges passing through points, which are identified, represent a single edge.}
		\label{fig_tore}
	\end{figure}
	\par We fix a tiling of $\Psi$. 
	We construct a graph $\Psi_1 = (V(\Psi_1), E(\Psi_1))$ by taking vertices $V(\Psi_1)$ as the centers of tiles and edges $E(\Psi_1)$ in such a way that the resulting graph $\Psi_1$ is the \textit{nearest-neighbor graph} with respect to the flat metric on $\Psi$.
	This means that an edge connects two vertices if and only if they are the closest neighbors with respect to the metric $g^{T \Psi}$.
	\par 
	The vector bundle $F_1$ over $\Psi_1$ and the Hermitian metric $h^{F_1}$ on $F_1$ are constructed by the restriction from $F$ and $h^F$. 
	The connection $\nabla^{F_1}$ is constructed using the parallel transport of $\nabla^F$ with respect to the straight path between the vertices.
	It is a matter of a simple verification to see that since $(F, h^F, \nabla^F)$ is unitary, the vector bundle $(F_1, h^{F_1}, \nabla^{F_1})$ is unitary as well.
	By considering regular subdivisions of tiles into $n^2$ squares, $n \in \nat^*$, and repeating the same procedure, we construct a family of graphs $\Psi_n = (V(\Psi_n), E(\Psi_n))$ with unitary vector bundles $(F_n, h^{F_n}, \nabla^{F_n})$ over $\Psi_n$, for $n \in \nat^*$. 
	Note that we have a natural injection
	\begin{equation}\label{eq_inj_vert}
		V(\Psi_n) \hookrightarrow \Psi.
	\end{equation}
	\par
\begin{figure}[!htbp]%
    \centering
    \subfloat{{\includegraphics[width=0.2\textwidth]{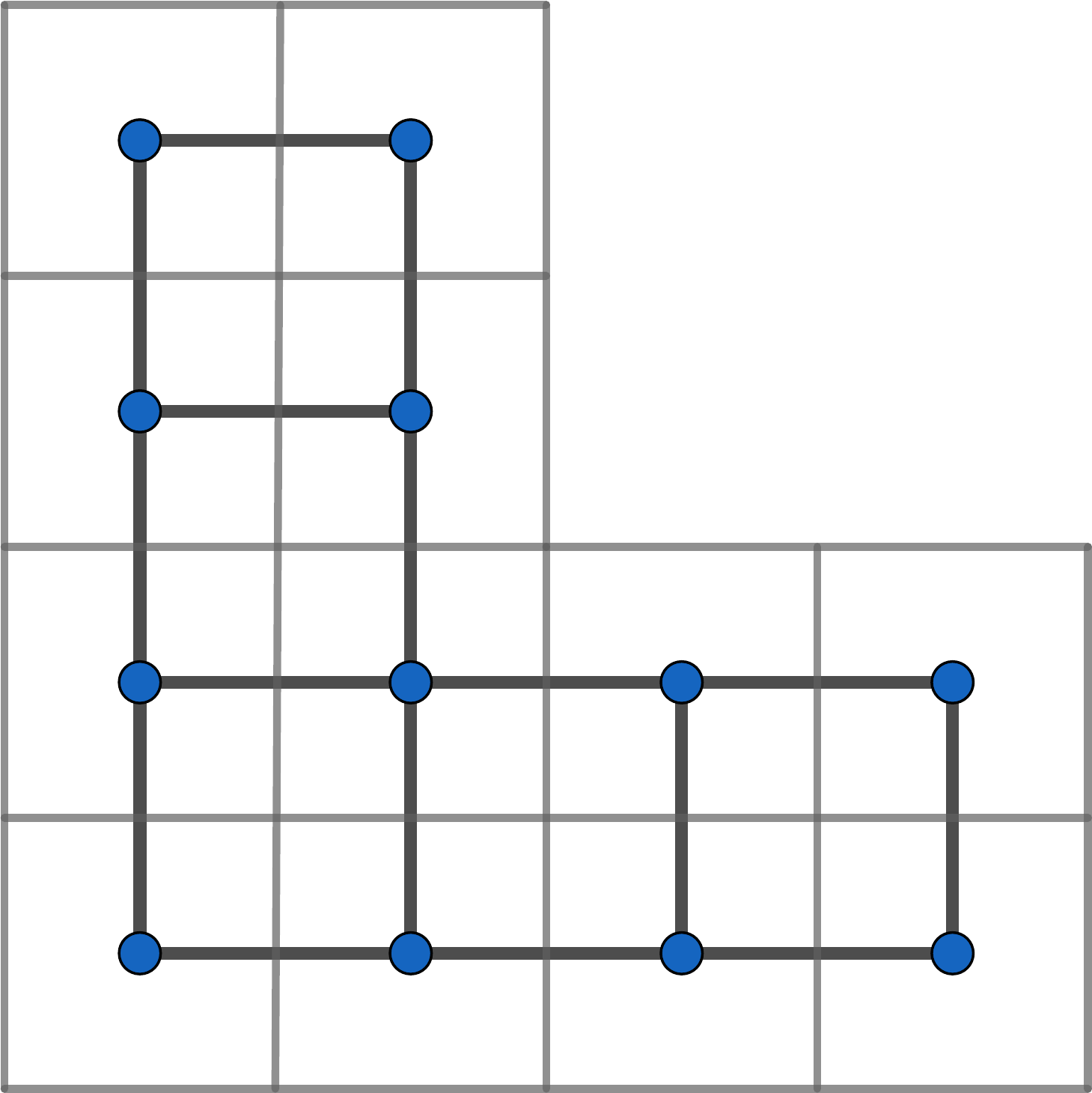}}}%
    \qquad \qquad \qquad \qquad
    \subfloat{{\includegraphics[width=0.2\textwidth]{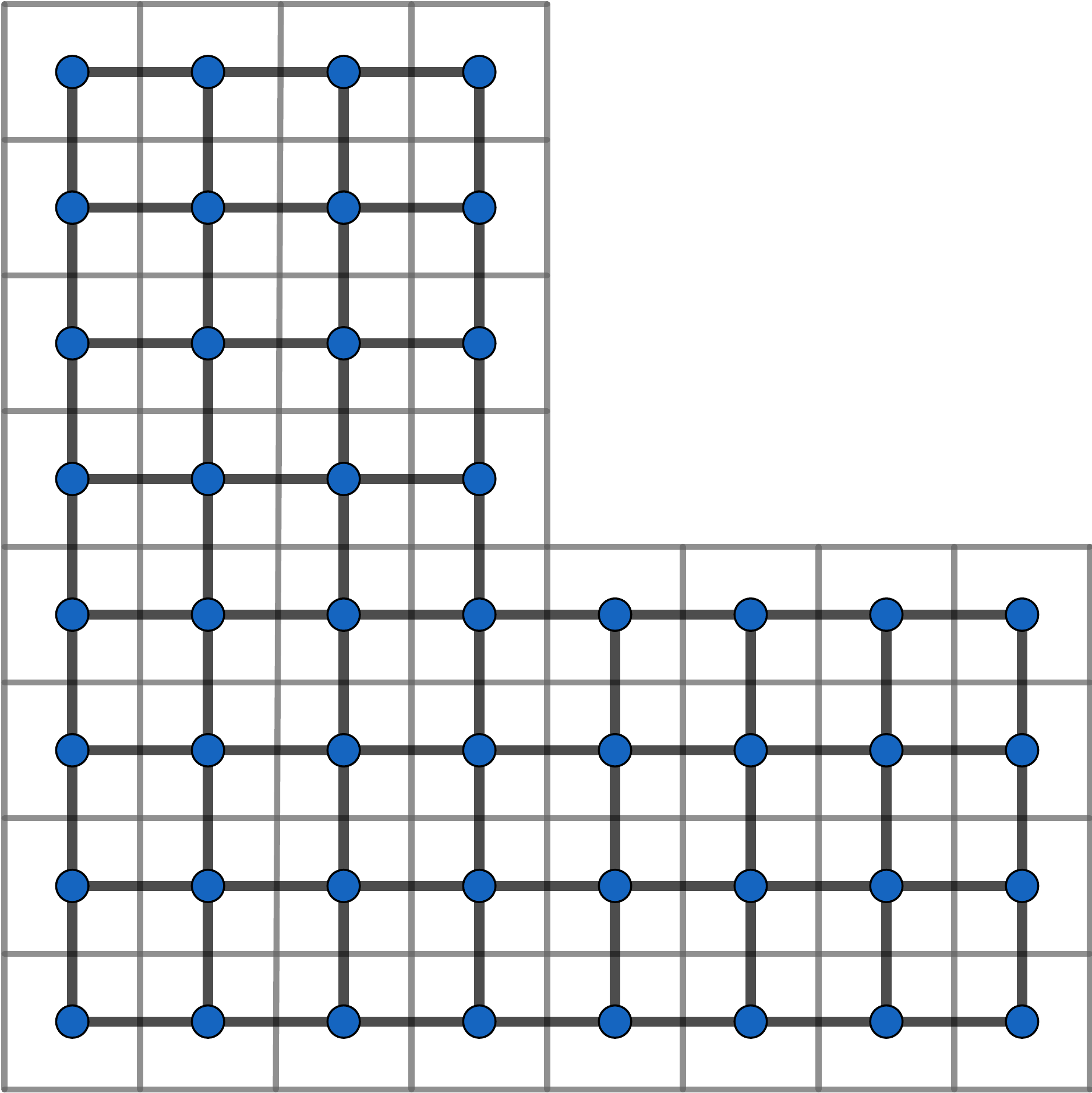}}}%
    \caption{An $L-$shape and its discretizations. It is a rectangular domain in $\comp$ with a single corner with angle $\frac{3 \pi}{2}$ and 5 corners  with angles $\frac{\pi}{2}$.}%
    \label{fig_l_shape}%
	\end{figure}

 For example, in case if $\Psi$ is a rectangular domain in $\comp$ with integer vertices, the family of graphs $\Psi_n$ coincides with subgraphs of $\frac{1 + \imun}{2n} + \frac{1}{n} \integ^2$, which stay inside of $\Psi$. See Figure \ref{fig_l_shape}.
	\par 
	We denote by $(\nabla^F)^*$ the formal adjoint of $\nabla^F$ with respect to the $L^2$-metric induced by $g^{T \Psi}$ and $h^F$.
	We denote by $\laplcomp_{\Psi}^{F}$ the scalar Laplacian on $\Psi$ associated with $(F, h^{F}, \nabla^{F})$. 
	It is a differential operator acting on the smooth sections of $F$ by
	\begin{equation}\label{eq_lapl}
		\laplcomp_{\Psi}^{F} := (\nabla^F)^* \nabla^F.
	\end{equation}
	If $(F, h^{F}, \nabla^{F})$ is trivial, $\laplcomp_{\Psi}^{F}$ coincides with the usual Laplacian, given by the formula $- \frac{\partial^2}{\partial x^2} - \frac{\partial^2}{\partial y^2}$.
	\par 
	In this paper we always consider $\laplcomp_{\Psi}^{F}$ with \textit{von Neumann boundary conditions} on $\partial \Psi$. In other words, the sections $f$ from the domain of our Laplacian satisfy
	\begin{equation}\label{eq_vn_bound_cond}
		\nabla^{F}_{n} f = 0 \quad \text{over } \partial \Psi,
	\end{equation}
	where $n$ is the normal to the boundary.
	\par It is well-known that unlike for smooth manifolds, the Laplacian $\laplcomp_{\Psi}^{F}$ is not necessarily \textit{essentially self-adjoint}. 	
	Thus, to define the spectrum of $\laplcomp_{\Psi}^{F}$, we will be obliged to specify the self-adjoint extension of $\laplcomp_{\Psi}^{F}$ we are working with.
	We choose the Friedrichs extension and, by abuse of notation, we denote it by the same symbol $\laplcomp_{\Psi}^{F}$. See Section \ref{sect_zeta_flat_fun} for some properties and definitions on Friedrichs extension and Section \ref{sect_funct_an_flat} for an explicit description of its domain.
	\par There are many ways to motivate this choice of a self-adjoint extension. To name one of them, this extension is positive (cf. \cite[Theorem X.23]{ReedSimonII}), which is rather handy since the discrete Laplacians we consider are positive as well (see the end of Section \ref{sect_sq_t_s_disc}). The main statements of this article, Theorems \ref{thm_eigval_convergence}, \ref{thm_eigvec_convergence}, could also be seen as another justification for such a choice.
	\par 
	Similarly to the classical case of smooth domains, the spectrum of $\laplcomp_{\Psi}^{F}$ is discrete (cf. Proposition \ref{prop_spec_discr}), in other words (by our convention, the sets take into account the multiplicity)
	\begin{equation}\label{eq_spec_list}
		\spec(\laplcomp_{\Psi}^{F}) = \{ \lambda_1, \lambda_2, \cdots \},
	\end{equation}
	where $\lambda_i$, $i \in \nat^*$ form a non-decreasing sequence.
	\par 
	The main goal of this article is to study the relationship between $\spec(\laplcomp_{\Psi_n}^{F_n})$ and $\spec(\laplcomp_{\Psi}^{F})$.
	\par 
	From (\ref{eq_comb_lapl}) and the fact that the degree of every vertex is bounded by $4$, we see that 
	\begin{equation}\label{eq_unif_bound_comb_lapl}
		\spec(\laplcomp_{\Psi_n}^{F_n}) \subset [0, 8 \cdot \rk{F}].
	\end{equation}
	Also, it is clear that the set $\spec(\laplcomp_{\Psi}^{F})$ is unbounded in $\real$.
	Thus, it makes more sense to look at the rescaled spectrum $\spec(n^2 \cdot \laplcomp_{\Psi_n}^{F_n})$, $n \in \nat^*$, which we denote by
	\begin{equation}\label{eq_spec_list_discr}
		\spec(n^2 \cdot \laplcomp_{\Psi_n}^{F_n}) = \{\lambda_{1}^{n}, \lambda_{2}^{n}, \cdots \},
	\end{equation}
	where $\lambda_{i}^{n}$, $i \in \nat^*$ form a non-decreasing sequence for any $n \in \nat^*$.
	\par 
	The first main result of this article is the following
	\begin{thm}\label{thm_eigval_convergence}
		Let $(\Psi, g^{T \Psi})$ be a half-translation surface endowed with a flat unitary vector bundle $(F, h^F, \nabla^F)$.
		Suppose $\Psi$ can be tiled by euclidean squares of area 1.
		Construct the family of graphs $\Psi_n$, $n \in \nat^*$, as above.
		Denote by $(F_n, h^{F_n}, \nabla^{F_n})$ the induced unitary vector bundles on $\Psi_n$.
		\par 
		In the notations of (\ref{eq_spec_list}) and (\ref{eq_spec_list_discr}), for any $i \in \nat$, as $n \to \infty$, the following limit holds
		\begin{equation}\label{eq_eigval_convergence}
			\lambda_{i}^{n} \to \lambda_i.
		\end{equation}
	\end{thm}
	\begin{rem}
	In particular, from Theorem \ref{thm_eigval_convergence}, we see that any spectral invariant depending on a finite number of eigenvalues with fixed indices is related in the limit with the spectrum of $\laplcomp_{\Psi}^{F}$.
	\end{rem}
	\par 
	The second result is a similar statement in realms of eigenvectors. Morally, it says that the eigenvectors of $n^2 \cdot \Delta_{\Psi_n}^{F_n}$ converge in $L^2(\Psi, F)$ to the eigenvectors of $\Delta_{\Psi}^{F}$.
	However, as the eigenvalues might have some multiplicity, and the finitely dimensional space ${\rm{Map}}(V(\Psi_n), F_n)$ doesn't inject in $L^2(\Psi, F)$ canonically, we need to introduce some further notation for a rigorous statement.
	\par 
	Assume that the eigenvalue $\lambda_i$, $i \in \nat^*$ of $\laplcomp_{\Psi}^{F}$ has multiplicity $m_i$.
	Let $f_{i, j} \in L^2(\Psi, F)$, $j = 1, \ldots, m_i$ be the orthonormal basis of eigenvectors of $\laplcomp_{\Psi}^{F}$ corresponding to the eigenvalue $\lambda_i$.
	By Theorem \ref{thm_eigval_convergence}, we conclude that there is a series of eigenvalues $\lambda_{i, j}^{n}$, $j = 1, \ldots, m_i$ of $n^2 \cdot \laplcomp_{\Psi_n}^{F_n}$, converging to $\lambda_i$, as $n \to \infty$. 
	Moreover, no other eigenvalue of $n^2 \cdot \laplcomp_{\Psi_n}^{F_n}$, $n \in \nat^*$ comes close to $\lambda_i$ asymptotically.
	In the beginning of Section \ref{sect_harn_in}, we define a “linearization" functional 
	\begin{equation}
		L_n : {\rm{Map}}(V(\Psi_n), F_n) \to L^2(\Psi, F).
	\end{equation}
	One should think of it as a sort of linear interpolation, (\ref{eq_inj_vert}), which “blurs" the function near the set  ${\rm{Con}}(\Psi)$ of conical points of $\Psi$ and the set ${\rm{Ang}}(\Psi)$ of non-smooth points of the boundary $\partial \Psi$.
	\par 
	We denote by $\norm{\cdot}_{L^2(\Psi_n, F_n)}$ the $L^2$-norm on ${\rm{Map}}(V(\Psi_n), F_n)$, given by 
	\begin{equation}
		\norm{f}_{L^2(\Psi_n, F_n)}^{2} = \sum_{v \in V(\Psi_n)} h^{F_n} \big(f(v), f(v) \big).
	\end{equation}
	\begin{thm}\label{thm_eigvec_convergence}
		We use the same notation as in Theorem \ref{thm_eigval_convergence}.
		For any $i \in \nat$, there is $N$ such that for any $n \geq N$, there are $f_{i, j}^{n} \in {\rm{Map}}(V(\Psi_n), F_n)$,$1 \leq j \leq m_i$, which are pairwise orthogonal, satisfy $\| f_{i, j}^{n} \|_{L^2(\Psi_n, F_n)}^{2} = n^2$, and which are in the span of the eigenvectors of $n^2 \cdot \laplcomp_{\Psi_n}^{F_n}$, corresponding to the eigenvalues $\lambda_{i, j}^{n}$, $j = 1, \ldots, m_i$,   such that, as $n \to \infty$, in $L^2(\Psi, F)$, the following limit holds 
		\begin{equation}\label{eq_eigvec_convergence}
			L_n(f_{i, j}^{n}) \to f_{i, j}.
		\end{equation}
	\end{thm}
	We remark that in the proofs of Theorems \ref{thm_eigval_convergence}, \ref{thm_eigvec_convergence}, we were inspired by the approximation theory of Dodziuk, \cite{Dodziuk}, and Dodziuk-Patodi, \cite{DodziukPatodi}.
	Note, however, that there are two big differences between their theory and ours. 
	First, they work with simplicial complexes, and their Laplacian on the discrete model is defined as $d_{\Psi_n}^{*} d_{\Psi_n}$, where $d_{\Psi_n}^{*}$ is the adjoint of $d_{\Psi_n}$ with respect to the pull-back of the $L^2$-metric on $\Psi$ through so-called \textit{Whitney map} (cf. Dodziuk \cite[\S 1]{Dodziuk}).  In particular, the Laplacian from \cite{Dodziuk}, \cite{DodziukPatodi} has nothing to do with the combinatorial Laplacian (\ref{eq_comb_lapl}) we are considering in this article (which depends only on the combinatorics of the approximation graph and doesn't depend on the metric).
	Second, the boundary of their manifold is smooth. Non-smoothness of the boundary in our case raises some technical problems which did not appear in the articles \cite{Dodziuk}, \cite{DodziukPatodi}.
	On continuous side, to overcome the lack of elliptic estimates near conical singularities and corners, which entail the non-differentiability of the eigenvectors there, we use weak elliptic regularity results on polygons due to Grisvard, \cite{Grisvard}. 
	On discrete side, we develop \textit{Harnack-type estimates} in Theorem \ref{thm_harn_type} for the eigenvectors corresponding to small eigenvalues of $\laplcomp_{\Psi_n}^{F_n}$ to prove that the discrete eigenvectors are “asymptotically continuous".
	To establish the Harnack-type estimates, we rely on the potential theory on lattices, introduced by Duffin \cite{Duffin} in dimension $3$ and developped by Kenyon \cite{Kenyon2002Invent} in dimension $2$.
	\par 
	This article is organized as follows. In Section \ref{sec_disc_cont} we introduce the main notions related to flat surfaces and their discretizations. We define the Friedrichs extension of the Laplacian and give an explicit description of its domain. In Section \ref{sect_fd_meth}, we prove Theorems \ref{thm_eigval_convergence}, \ref{thm_eigvec_convergence}, modulo some Harnack-type estimates, which are proved in Section \ref{sect_harn_pf}.
	\par  This paper will be used in \cite{FinAsDet} to relate the asymptotic expansion of the number of spanning trees, spanning forests and weighted cycle-rooted spanning forests on the discretizations of flat surfaces to the corresponding zeta-determinants.
	The results of those papers are announced in \cite{FinDetCRAS}.
	\par 
	\textbf{Notation.}
	For a graph $G$, we denote by $V(G)$, $E(G)$ the sets of vertices and edges of $G$ respectively.
	For a Hermitian vector bundle $(V, h^V)$ on a finite graph $G$, we denote by $\scal{\cdot}{\cdot}_{L^2(G, V)}$ the $L^2$-scalar product on the set ${\rm{Map}}(V(G), V)$, defined by the following formula
	\begin{equation}
		\scal{f}{g}_{L^2(G, V)} = \sum_{v \in V(G)} h^V(f(v), g(v)), \qquad f, g \in {\rm{Map}}(V(G), V)
	\end{equation}
	Recall that the divergence operator $d_{G}^{*} : {\rm{Map}}(E(G), \comp) \to {\rm{Map}}(V(G), \comp)$ is defined by
	\begin{equation}\label{eq_div_defn}
		(d_{G}^{*} g)(P) := \sum_{\substack{e \in E(G), \\ t(e) = P}} g(h(e)), \qquad g \in {\rm{Map}}(E(G), \comp).
	\end{equation}
	Then the following identities can be verified directly
	\begin{equation}\label{eq_scal_laplg_dg}
		\laplcomp_{G} = d_{G}^{*} d_{G},
		\qquad \qquad
		\scal{\laplcomp_{G} f}{g}_{L^2(G)} = \scal{d_{G} f}{d_{G} g}_{L^2(G)}.
	\end{equation}
	For $k \in \nat$, we denote by 
	\begin{equation}\label{eq_defn_ck_b}
		\ccal^{k}(\overline{\Psi}, F) := \Big \{ 
			f \in \ccal^{k}(\Psi, F) : \nabla^l f \in L^{\infty}(\Psi, F), \text{ for }l \leq k
		\Big \},
	\end{equation}
	\par 
	\textbf{Acknowledgements.} I would like to thank Dmitry Chelkak, Yves Colin de Verdière for related discussions and their interest in this article, and especially Xiaonan Ma for important comments and remarks. I also would like to thank the colleagues and staff from Institute Fourier, Université Grenoble Alpes, where this article has been written, for their hospitality.

\section{Analysis on flat surfaces}\label{sec_disc_cont}
	In this section we recall some results about functional analysis on flat surfaces and introduce the main objects of this article. More precisely, in Section \ref{sect_sq_t_s_disc}, we recall the basics of flat surfaces, pillowcase covers, and we  briefly introduce the main notions on vector bundles over graphs.
	In Section \ref{sect_zeta_flat_fun}, we define the Friedrichs extension of the Laplacian and study some of its spectral properties.
	Finally, in Section \ref{sect_funct_an_flat}, we give an alternative description of the domain of the Friedrichs extension by describing the singularities of the functions from it.
	
	\subsection{Pillowcase covers and their discretizations}\label{sect_sq_t_s_disc}
	\par Here we recall the definition of flat surfaces, pillowcase covers, explain some properties of discretizations of pillowcase covers and give a short introduction to vector bundles over graphs.
	\par 
	By Gauss-Bonnet theorem, the only closed Riemann surface admitting a flat metric has the topology of the torus. 
	However, any Riemann surface can be endowed with a flat metric having a finite number of cone-type singularities. 
	Let's explain this point more precise.
	\par 
	A cone-type singularity is a Riemannian metric
	\begin{equation}\label{eq_conical_metric}
		ds^2 = dr^2 + r^2 dt^2,
	\end{equation}
	on the manifold
	\begin{equation}\label{eq_v_theta}
		C_{\theta} := \{(r, t) : r > 0; t \in \real / \theta \integ \},
	\end{equation}
	where $\theta > 0$.
	In what follows, when we speak of cones, we assume $\theta \neq 2\pi$.
	By a \textit{flat metric with a finite number of cone-type singularities} we mean a metric defined away from a finite set of points such that there is an atlas for which the metric looks either like the standard metric on $\real^2$, or like the conical metric (\ref{eq_conical_metric}) on an open subset of (\ref{eq_v_theta}).
	\par 
	Let $M$ be a compact surface endowed with a flat metric with cone-type singularities of angles $\theta_1, \ldots, \theta_m$ at points  $P_1, \ldots, P_m$. 
	By Gauss-Bonnet theorem (cf. Troyanov \cite[Proposition 3]{Troyanov}):
	\begin{equation}\label{eq_gauss_bonnet}
		\sum_{i = 1}^{m} (2\pi - \theta_i) = 2 \pi \chi(M),
	\end{equation}
	where $\chi(M)$ is the Euler characteristic of $M$. 
	Although we will not need this in what follows, by a theorem of Troyanov \cite[Théorème, \S 5]{Troyanov}, this is the only obstruction for the existence of a metric with cone-type singularities of angles $\theta_1, \ldots, \theta_m$ at points  $P_1, \ldots, P_m$ in a conformal class of $M$.
	\par 
 	In this article we are primarily interested in compact surfaces endowed with a metric with a finite number of cone-type singularities of angles $\pi k$, $k \in \nat^*$. 
 	Those Riemann surfaces can also be described by possession of an atlas, defined away from the singularities of the metric, such that the transition maps between charts are given by $z \to \pm z + {\rm{const}}$ (cf. \cite[\S 3.3]{ZorichFlatSurf}).
	In literature, such surfaces are called \textit{half-translation surfaces}.
 	In case if all the angles are of the form $2 \pi k$, $k \in \nat^*$, it can be proven that the atlas can be chosen in such a way so that the transition maps between charts are given by $z \to z + {\rm{const}}$, and the surface in this case is called a \textit{translation surface}.
 	Clearly, from any closed half-translation surface, we can construct a translation surface as a double cover, ramified over the conical points with angles $(2 k + 1) \pi$, $k \in \nat^*$.
 	\par 
 	We denote by ${\rm{Con}}(\Psi)$ the conical points of the surface $\Psi$, and by ${\rm{Ang}}(\Psi)$ the points where two different smooth components of the boundary meet (corners). 
	We denote by $\angle : {\rm{Con}}(\Psi) \to \real$ the function which associates to a conical point its angle and by $\angle : {\rm{Ang}}(\Psi) \to \real$ the function which associates the interior angle between the smooth components of the boundary.
 	\par 
 	It is easy to see that for a half-translation surface, a notion of a \textit{straight line} makes sense, and for translation surfaces, even a notion of a \textit{ray} (or a \textit{direction}) is well-defined.
 	\par There is an alternative description of closed \textit{half-translation surfaces} with prescribed line (resp. \textit{translation surfaces} with prescribed direction) in terms of Riemann surfaces endowed with a meromorphic quadratic differential with at most simple poles (resp. a holomorphic differential). The zeros\footnote{Here we interpret a pole of a quadratic differential as a zero of order $-1$.} of order $k$ in this description correspond to the conical points with angles $(k + 2) \pi$ (resp. $2 \pi (k + 1)$). In this description, one could easily identify the moduli spaces of closed half-translation surfaces with line (resp. translation surfaces with direction) with corresponding Hodge bundles on the moduli of curves (cf. \cite[\S 8.1]{ZorichFlatSurf}). The orders of zeros of a meromorphic quadratic differential (resp. of a holomorphic differential) induce stratifications of those moduli spaces.
 	\par 
	For the most part of this paper, we will be interested in considering special type of half-translation surfaces which are called \textit{pillowcase covers}. Those are surfaces with a fixed tiling by squares of equal sizes.
	In case if the surface has no boundary, those surfaces can be thought as rational points in so-called period coordinates (cf. \cite[\S 7.1]{ZorichFlatSurf}) inside the respective stratas of the moduli space, and thus, the set of pillowcase covers form a dense subset of the corresponding moduli spaces (similarly to the set of tori with rational periods inside the moduli space of tori).
	\par
\begin{figure}[h]
		\includegraphics[width=0.2\textwidth]{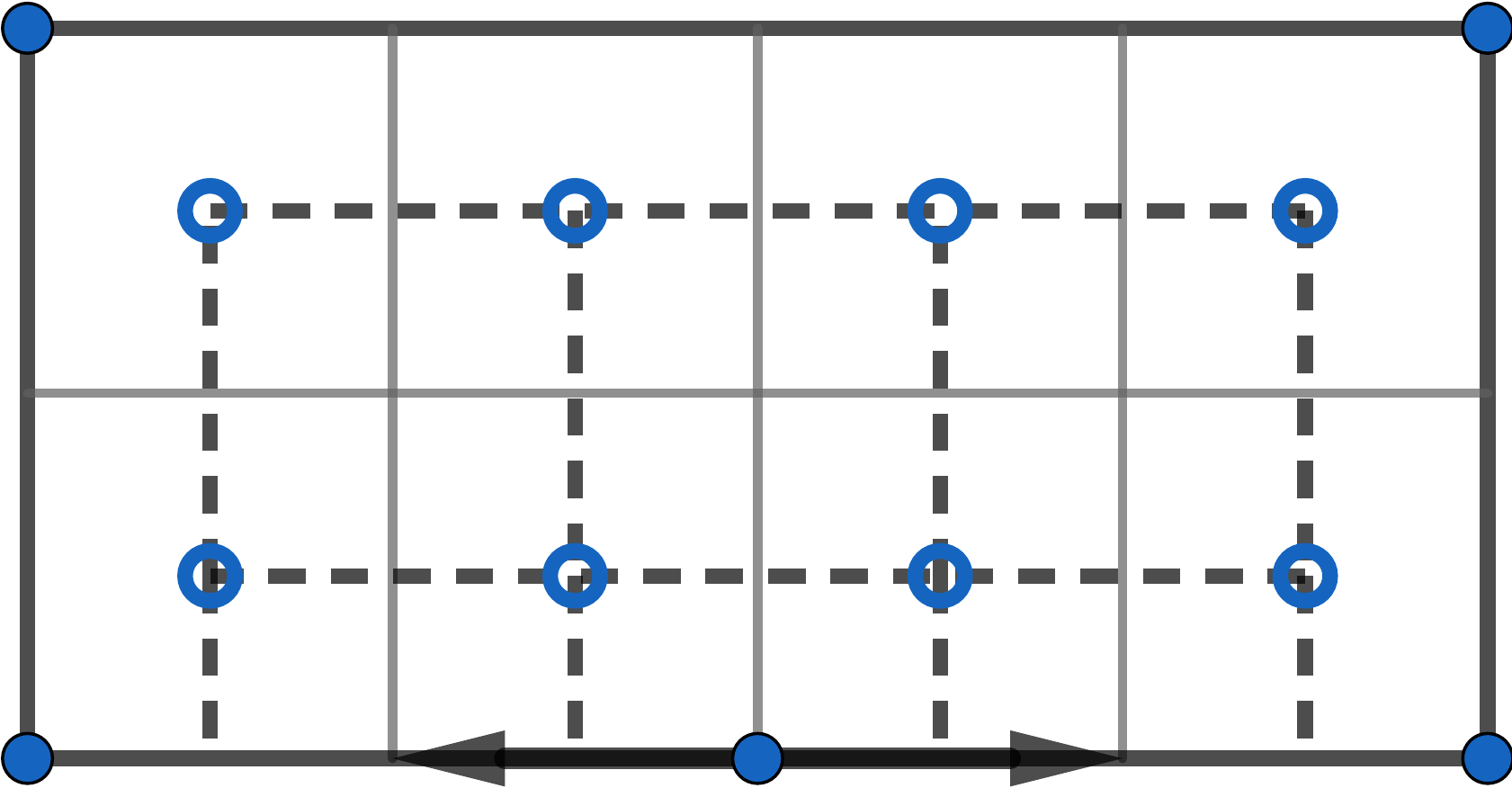}	
		\centering
		\caption{An example of multiple edges in a discretization. The surface is obtained by identifying the edges along the directed lines. The edges of the graph are dotted. }
		\label{fig_cpi_graph}
	\end{figure}

 Let's now set up some notations associated with discretization $\Psi_n$ of a pillowcase cover $\Psi$, constructed as in Introduction.
For $P \in {\rm{Con}(\Psi)} \cup {\rm{Ang}(\Psi)}$, we define the set $V_n(P)$ as follows
	\begin{equation}\label{eq_defn_vp}
		V_n(P) := \Big\{ Q \in V(\Psi_n) : \dist_{\Psi}(P, Q) = \frac{1}{2n} \Big\},
	\end{equation}
	in other words, $V_n(P)$ is the set of the nearest neighbors of $P$ from the vertex set $V(\Psi_n)$ with respect to the flat metric on $\Psi$. 
	It's easy to verify that for $n \geq 2$, we have 
	\begin{equation}
		\# V_n(P) = \frac{2 \angle(P)}{\pi}.
	\end{equation}
	Let's define for $P \in {\rm{Con}(\Psi)} \cup {\rm{Ang}(\Psi)}$ the open subset $U_n(P) \subset \Psi$ as follows
	\begin{equation}\label{eq_defn_up}
		U_n(P) := \Big\{ z \in \Psi : {\rm{dist}}_{\Psi}(z, P) < \frac{1}{2n} \Big\}.
	\end{equation}
	Remark that the points $V_n(P)$ lie on the boundary of $U_n(P)$.
	\begin{rem}
		For $n \geq 2$, the edges of $\Psi_n$ have at most double multiplicity, see Figure \ref{fig_cpi_graph}.
		Moreover, the number of edges with double multiplicity is equal to $\# \{ P \in {\rm{Con}}(\Psi) : \angle(P) = \pi\}$. We are, thus, working with multigraphs, but by abuse of notation, we call them graphs.
	\end{rem}
	\par Now, by a \textit{vector bundle} $V$ on a graph $G = (V(G), E(G))$, we mean the choice of a vector space $V_v$ for any $v \in V(G)$ so that for any $v, v' \in V(G)$, the vector spaces $V_v$ and $V_{v'}$ are isomorphic.
	The set of \textit{sections} ${\rm{Map}}(V(G), V)$ of $V$ is defined by
	\begin{equation}
		{\rm{Map}}(V(G), V) = \oplus_{v \in V(G)} V_v.
	\end{equation}
	A \textit{connection} $\nabla^{V}$ on a vector bundle $V$ is the choice for each edge $e =(v, v') \in E(G)$ of an isomorphism $\phi_{v v'}$ between the corresponding vector spaces $\phi_{v v'} : V_v \to V_{v'}$, with the property that $\phi_{v v'} = \phi_{v' v}^{-1}$. This isomorphism is called the \textit{parallel transport} of vectors in $V_v$ to vectors in $V_{v'}$.
	\par 
	A \textit{Hermitian metric} $h^{V}$ on the vector bundle $V$ is a choice of a positive-definite Hermitian metric $h_{v}$ on $V_v$ for each $v \in V(G)$. We say that a connection $\nabla^{V}$ is unitary with respect to $h^V$ if its parallel transports preserve $h^{V}$.
	\par The \textit{Laplacian} $\laplcomp_G^{V}$ on $(V, \nabla^V)$ is the linear operator $\laplcomp_G^{V} : {\rm{Map}}(V(G), V) \to {\rm{Map}}(V(G), V)$, defined for $f \in {\rm{Map}}(V(G), V)$ by (\ref{eq_comb_lapl}).
	Remark that unlike Laplace–Beltrami operator on a smooth manifold, we don't use the metric to define the Laplacian (\ref{eq_comb_lapl}).
	\par Consequently, in general, the operator $\laplcomp_G^{V}$ is not self-adjoint, see for example \cite[\S 3.2 and equation (1)]{KenyonFlatVecBun}.
	However, if one assumes that the connection is unitary with respect to $h^{V}$, then it becomes self-adjoint (cf. Kenyon \cite[\S 3.3]{KenyonFlatVecBun}). 
	\par 
	More precisely, one can extend the definition of a vector bundle to the edges of $G$. 
	A vector bundle $V'$ over $V(G) \oplus E(G)$ is a choice of a vector spaces $V_e$ for each edge $e \in E(G)$ as well as $V_v$ for each vertex $v \in V(G)$. 
	A connection $\nabla^{V'}$ on $V'$ is a choice of a connection $\nabla^V$ on $V$ as well as connection isomorphisms $\phi_{ve} : V_v \to V_e$, $\phi_{ev} : V_e \to V_v$ for each $v \in V(G)$ and $e \in E(G)$, and satisfying $\phi_{ve} = \phi_{ev}^{-1}$ and
	$\phi_{v v'} = \phi_{ev'} \circ \phi_{ve}$.
	Similarly to ${\rm{Map}}(V(G), V)$, we defined the set of sections of ${\rm{Map}}(E(G), V)  := \oplus_{e \in E(G)} V_e$.
	\par Quite easily, for any vector bundle $V$ and a connection $\nabla^V$ on $V(G)$, we may extend it to a vector bundle $V'$ and a connection $\nabla^{V'}$ on $E(G) \oplus V(G)$.
	Note, however, that such a choice would not be unique.
	If the initial vector bundle $V$ is endowed with a Hermitian metric $h^V$, for which the connection $\nabla^V$ is unitary, then one might endow the vector bundles $V_e$, $e \in E(G)$ with metrics and choose the connections $\phi_{ev}$, $e \in E(G)$, $v \in V(G)$ so that $\nabla^{V'}$ is unitary as well.
	\par There is a natural map $\nabla^V_{G} : {\rm{Map}}(V(G), V) \to {\rm{Map}}(E(G), V)$, defined as follows 
	\begin{equation}
		(\nabla^V_{G} f)(e)
		=
		\phi_{t(e)e} f(t(e)) - \phi_{h(e)e} f(h(e)), 
		\qquad
		f \in {\rm{Map}}(V(G), V), 
	\end{equation}	 
	where $t(e)$ and $h(e)$ are tail and head respectively of an oriented edge $e$.
	We also define the operator $(\nabla^V_{G})^* : {\rm{Map}}(E(G), V) \to {\rm{Map}}(V(G), V)$ by the formula
	\begin{equation}
		((\nabla^V_{G})^* f)(v) = \sum_{\substack{e \in E(G) \\ t(e) = v}} \phi_{e v} f(e).
	\end{equation}
	It is an easy verification (cf. Kenyon \cite[\S 3.3]{KenyonFlatVecBun}) that for the Laplacian, defined by (\ref{eq_comb_lapl}), we have
	\begin{equation}\label{eq_lapl_self_adjoint}
		\laplcomp_G^{V} = (\nabla^V_{G})^* \nabla^V_{G}.
	\end{equation}
	Note, however, that in general $(\nabla^V_{G})^*$ is not the adjoint of $\nabla^V_{G}$ with respect to the appropriate $L^2$-metrics. 
	But if the connection $\nabla^V$ is unitary, it is indeed the case, cf. \cite[\S 3.3]{KenyonFlatVecBun}.
	\par 
	In this article, all our connections are unitary, and thus, by (\ref{eq_lapl_self_adjoint}), the associated discrete Laplacians are self-adjoint and positive.

	\subsection{Properties of Friedrichs extension of the Laplacian}\label{sect_zeta_flat_fun}
		In this section we study the Laplacian $\laplcomp_{\Psi}^{F}$ associated to a flat surface $\Psi$ with conical singularities and piecewise geodesic boundary endowed with a flat unitary vector bundle $(F, h^F, \nabla^F)$ over $\overline{\Psi}$, (\ref{eq_compact_defn}).
		The content of this section is certainly not new, but we weren't able to find a complete reference for all the results contained here.
		\par 
		We consider $\laplcomp_{\Psi}^{F}$ as an operator acting on the functional space $\ccal^{\infty}_{0, vN}(\Psi, F)$, where
		\begin{equation}
			\ccal^{\infty}_{0, vN}(\Psi, F) := \Big\{ f \in \ccal^{\infty}_{0} \big(\Psi \setminus  {\rm{Ang}(\Psi)} \big) : \nabla_n f = 0 \text{ over } \partial \Psi \Big\},
		\end{equation}	
		where $n$ is the normal to $\partial \Psi$.
		Unlike in the case of a manifold with smooth boundary, the operator $\laplcomp_{\Psi}^{F}$ is in general not essentially self-adjoint.
		\par 
		Let ${\rm{Dom}}_{\max}(\laplcomp_{\Psi}^{F})$ denote the maximal closure of $\laplcomp_{\Psi}^{F}$. In other words, for $u \in L^2(\Psi, F)$, we have $u \in {\rm{Dom}}_{\max}(\laplcomp_{\Psi}^{F})$ if and only if $\laplcomp_{\Psi}^{F} u \in L^2(\Psi, F)$, where $\laplcomp_{\Psi}^{F} u$ is viewed as a current.
		\par 
	 	Let's denote by $W^k_{p}(\Psi, F)$ the Sobolev space on $\Psi$, defined as
	 	\begin{equation}
	 		W^k_{p}(\Psi, F) = \Big\{ u \in L^p(\Psi, F) : \nabla^l u \in L^p(\Psi, F) \quad \text{ for } \quad 0 \leq l \leq k  \Big\}.
	 	\end{equation}
	 	We denote by $\norm{\cdot}_{W^k_{p}(\Psi, F)}$ the norm on $W^k_{p}(\Psi, F)$, given for $u \in W^k_{p}(\Psi, F)$ by
	 	\begin{equation}
	 		\norm{u}_{W^k_{p}(\Psi, F)} = \sum_{0 \leq l \leq k} \big\| \nabla^l u \big\|_{L^p(\Psi, F)}.
	 	\end{equation}
	 	We use the following shorthand notations $H^{k}(\Psi, F) := W^k_{2}(\Psi, F)$ and $\norm{\cdot}_{H^{k}(\Psi, F)} := \norm{\cdot}_{W^k_{2}(\Psi, F)}$.
		 \begin{thm}[Rellich-Kondrachov]\label{thm_rell_kondr}
	 		The inclusion $H^1(\Psi, F) \hookrightarrow L^2(\Psi, F)$ is compact.
	 	\end{thm}
	 	\begin{proof}
	 		The mentioned inclusion is trivially continuous. Let's prove that it is compact.
	 		\par 
	 		Geometrically it is trivial that one can choose a finite number of functions $f_i: \Psi \to \{0, 1\}$, $i \in I$ such that for any $x \in \Psi$ away from a negligible set, there exists exactly one $i \in I$ such that $f_i(x) = 1$, and, for any $i \in I$,  $U_i := {\rm{supp}}(f_i)$ is isomorphic to a subdomain of $\comp$ with piecewise smooth boundary. 
	 		In particular, in the notations of Adams \cite[p. 66]{Adams}, $U_i$, $i \in I$, satisfy the cone property.
	 		By \cite[Theorem 6.2, Page 144]{Adams}, Rellich-Kondrachov theorem holds for any subsets of $\comp$, satisfying the cone property. Thus, the inclusions $H^1(U_i) \hookrightarrow L^2(U_i)$ are compact for any $i \in I$.
	 		Note that we assumed that $(F, h^F, \nabla^F)$ is flat unitary, and thus we can trivialize it over any star-like domain. In other words, we may interpret the section of $F$ over $U_i$ as functions.
	 		Thus, the inclusions $H^1(U_i, F) \hookrightarrow L^2(U_i, F)$ are compact as well.
	 		\par Now, let $u_j$, $j \in J$ be some bounded sequence in $H^1(\Psi, F)$. 
	 		For $i \in I$, consider a sequence $f_i u_j \in H^1(U_i, F)$, $j \in J$. Since the set $I$ is finite, by the mentioned version of Rellich-Kondrachov theorem, one could choose a subsequence of $u_k$, $k \in K \subset J$ such that for any $i \in I$, the sequence $f_i u_j$, $j \in K$ converges in $L^2(U_i, F)$.
	 		However, by the choice of $f_i$, the following identity
	 		\begin{equation}
	 			u_k = \sum_{i \in I} f_i u_k
	 		\end{equation}
	 		holds in $L^2(\Psi, F)$. Thus, the sequence $u_k$, $k \in K$ converges in $L^2(\Psi, F)$, and as a consequence, the inclusion $H^1(\Psi, F) \hookrightarrow L^2(\Psi, F)$ is compact.
	 	\end{proof}
	 	For any positive symmetric operator, one can construct in a canonical way a self-adjoint extension, called Friedrichs extensions, through the completion of the associated quadratic form, cf. \cite[Theorem X.23]{ReedSimonII}.
	 	Once the definition is unraveled, the domain ${\rm{Dom}}_{Fr}(\laplcomp_{\Psi}^{F})$ of the Friedrichs extensions of the Laplacian $\laplcomp_{\Psi}^{F}$ on $\Psi$ with von Neumann boundary conditions on $\partial \Psi$ is given by
	 	\begin{equation}\label{eq_dom_friedr}
	 		{\rm{Dom}}_{Fr}(\laplcomp_{\Psi}^{F}) = {\rm{Dom}}_{\max}(\laplcomp_{\Psi}^{F})  \cap H^{1}_{0, vN}(\Psi, F),
	 	\end{equation}
	 	where $H^{1}_{0, vN}(\Psi, F)$ is the closure of $\ccal^{\infty}_{0, vN}(\Psi, F)$ in $H^1(\Psi, F)$.
	 	The value of the Friedrichs extensions of the Laplacian $\laplcomp_{\Psi}^{F}$ on $f \in {\rm{Dom}}_{Fr}(\laplcomp_{\Psi}^{F})$ is defined in the distributional sense. By the definition of ${\rm{Dom}}_{\max}(\laplcomp_{\Psi}^{F})$, it lies in $L^2(\Psi, F)$.
	 	\begin{prop}\label{prop_spec_discr}
	 		The spectrum of $\laplcomp_{\Psi}^{F}$ is discrete.
	 	\end{prop}
	 	 \begin{proof}
	 		Since the Friedrichs extension is non-negative by construction (cf. \cite[Theorem X.23]{ReedSimonII}), the kernel of $1 + \laplcomp_{\Psi}^{F}$ is empty.
	 		Thus, the inverse $(1 + \laplcomp_{\Psi}^{F})^{-1}$ is well-defined. By (\ref{eq_dom_friedr}), the image of  $(1 + \laplcomp_{\Psi}^{F})^{-1}$ lies in $H^1(\Psi, F)$. 
	 		By this and Theorem \ref{thm_rell_kondr}, we deduce that $(1 + \laplcomp_{\Psi}^{F})^{-1}$ is a compact operator.
	 		In particular, it has a discrete spectrum with only one possible accumulation point at $0$, which of course implies that the spectrum of $\laplcomp_{\Psi}^{F}$ is discrete.
	 	\end{proof}
		Now, recall that for a smooth codimension $1$ submanifold $\Gamma \subset \Psi$, transversal to the boundary $\partial \Psi$, and $k \geq \frac{1}{p}$, the \textit{trace} (or restriction) operator
		\begin{equation}\label{eq_trace_op}
			W^k_{p}(\Psi, F) \to W^{k - \frac{1}{p}}_{p}(\Gamma, F),  \qquad \qquad f \mapsto f|_{\Gamma},
		\end{equation}
		is well-defined, cf. \cite[Theorem 1.5.1.1]{Grisvard}. 
		In the following proposition and after, we apply (\ref{eq_trace_op}) implicitly when we mention an integration over $\Gamma$ of a function from $W^{k}_{p}(\Psi, F)$, $k \geq \frac{1}{p}$.
		\begin{prop}[Green's identity]\label{prop_green_identity}
			For any open subset $U \subset \Psi$ with piecewise smooth boundary $\partial U$ not passing through ${\rm{Con}}(\Psi)$ and ${\rm{Ang}}(\Psi)$, and any $u, v \in {\rm{Dom}}_{Fr}(\laplcomp_{\Psi}^{F})$, we have
			\begin{equation}\label{eq_green_identity}
				\scal{\laplcomp_{\Psi}^{F} u}{v}_{L^2(U, F)} = \scal{\nabla^F u}{\nabla^F v}_{L^2(U, F)} - \int_{\partial U} \nabla_n^{F} u  \cdot v dv_{\partial U},
			\end{equation}
			where $n$ is the outward normal to the boundary $\partial U$, and to simplify the notations, we omit the pointwise scalar product induced by $h^F$ in the last integral.
		\end{prop}
		\begin{rem}
			For $u, v \in {\rm{Dom}}_{\max}(\laplcomp_{\Psi}^{F})$, the statement (\ref{eq_green_identity}) is false. For example, take a trivial vector bundle, and for a radial function $r$, centered in some conical singularity, take $u = 1$ and $v = \log(r)$, multiplied by some bump function away from the singularity.
		\end{rem}
		\begin{proof}
			First of all, let's explain why all the terms on the right-hand side of (\ref{eq_green_identity}) are well-defined.
			We fix a smooth neighborhood $V$ of $\partial U$ with smooth boundary $\partial V$ such that $V \cap ({\rm{Con}}(\Psi) \cup {\rm{Ang}}(\Psi)) = \emptyset$.
			Since $V$ is a smooth domain, from elliptic estimates and the fact that $u, \laplcomp_{\Psi}^{F} u \in L^2(\Psi, F)$, we deduce that the restriction $u|_{V}$ of $u$ to $V$, viewed as a distribution, lies in $H^2(W, F)$.
			Thus by the trace theorem (\ref{eq_trace_op}), for the normal $n$ of $\partial U$, we have $\nabla_n^{F} u \in H^{\frac{1}{2}}(\partial U, F)$ and $v \in H^{\frac{3}{2}}(\partial U, F)$. 
			In particular, the integral $\int_{\partial U} \nabla_n^{F} u \cdot v dv_{\partial U}$ is well-defined.
			By (\ref{eq_dom_friedr}), we see that the product $\scal{\nabla^F u}{\nabla^F v}_{L^2(U, F)}$ is well defined as well.
			\par Now, let's prove (\ref{eq_green_identity}) for $u \in {\rm{Dom}}_{Fr}(\laplcomp_{\Psi}^{F})$ and $v \in \ccal^{\infty}_{0, vN}(\Psi, F)$. 
			Let $W \subset \Psi$, $W \cap ({\rm{Con}}(\Psi) \cup {\rm{Ang}}(\Psi)) = \emptyset$ be an open set with smooth boundary $\partial W$, such that 
			\begin{equation}\label{eq_supp_v}
				{\rm{supp}}(v) \subset W.
			\end{equation}
			As before, we deduce that the restriction $u|_W$ of $u$ to $W$ lies inside of $H^2(W, F)$. 
			By Grisvard \cite[Theorem 1.5.3.3]{Grisvard}, we see that
			\begin{equation}\label{eq_suppl_green_0}
				\scal{\laplcomp_{\Psi}^{F} u}{v}_{L^2(W, F)} = \scal{\nabla^F u}{\nabla^F v}_{L^2(W, F)} - \int_{\partial W} \nabla_n^{F} u \cdot v dv_{\partial W}.
			\end{equation}
		 	Now, by (\ref{eq_supp_v}) and the fact that $\laplcomp_{\Psi}^{F} u \in L^2(\Psi, F)$, we deduce
		 	\begin{equation}\label{eq_suppl_green_1}
		 		\scal{\laplcomp_{\Psi}^{F} u}{v}_{L^2(W, F)} = \scal{\laplcomp_{\Psi}^{F} u}{v}_{L^2(\Psi, F)}.
		 	\end{equation}
		 	By (\ref{eq_supp_v}) and the fact that $u \in H^1(\Psi, F)$, we deduce
		 	\begin{equation}\label{eq_suppl_green_2}
		 		\scal{\nabla^F u}{\nabla^F v}_{L^2(W, F)} = \scal{\nabla^F u}{\nabla^F v}_{L^2(\Psi, F)}.
		 	\end{equation}
		 	By (\ref{eq_supp_v}), the fact that $u \in H^2(W, F)$, thus, $\nabla_n^{F} u \in H^{\frac{1}{2}}(\partial W, F) \subset L^2(\partial W, F)$,  and the identity $\nabla_n u = 0$ over $\partial \Psi$, we deduce that 
		 	\begin{equation}\label{eq_suppl_green_3}
		 		\int_{\partial W} \nabla_n^{F} u \cdot v dv_{\partial W}
		 		=
		 		\int_{\partial U} \nabla_n^{F} u \cdot v dv_{\partial U}.
		 	\end{equation}
		 	From (\ref{eq_suppl_green_0})-(\ref{eq_suppl_green_3}), we conclude that (\ref{eq_green_identity}) holds for $u \in {\rm{Dom}}_{Fr}(\laplcomp_{\Psi}^{F})$ and $v \in \ccal^{\infty}_{0, vN}(\Psi, F)$. 
		 	\par Now, let's argue why (\ref{eq_green_identity}) also holds for $u \in {\rm{Dom}}_{Fr}(\laplcomp_{\Psi}^{F})$ and $v \in H^1_{0, vN}(\Psi, F)$. Indeed, by definition of the space $H^1_{0, vN}(\Psi, F)$, there is a sequence $v_m \in \ccal^{\infty}_{0, vN}(\Psi, F)$ such that $v_m \to v$ in $H^1(\Psi, F)$. Then, since $ \nabla_n^{F} u \in L^2(\partial U, F)$, we see by the continuity of the trace operator that 
		 	\begin{equation}
		 		\int_{\partial U} \nabla_n^{F} u \cdot v_m dv_{\partial U} \to  \int_{\partial U} \nabla_n^{F} u \cdot v dv_{\partial U}.
		 	\end{equation}
		 	Similarly, the first and the second terms of (\ref{eq_green_identity}) associated to $v_m$ converge to the respective terms associated to $v$. In the limit we obtain that (\ref{eq_green_identity}) holds for $u \in {\rm{Dom}}_{Fr}(\laplcomp_{\Psi}^{F})$ and $v \in H^1_{0, vN}(\Psi, F)$. 
		 	By (\ref{eq_dom_friedr}), we see that it implies that (\ref{eq_green_identity}) holds for $u, v \in {\rm{Dom}}_{Fr}(\laplcomp_{\Psi}^{F})$.
		\end{proof}
	 	\begin{cor}\label{cor_kernel_const}
	 		The kernel of  $\laplcomp_{\Psi}^{F}$ consists of flat sections of $F$. In other words,
	 		\begin{equation}
	 			\ker \laplcomp_{\Psi}^{F} \simeq H^0(\Psi, F).
	 		\end{equation}
	 	\end{cor}
	 	\begin{proof}
	 		 First of all, it's easy to see that the flat sections are from ${\rm{Dom}}_{Fr}(\laplcomp_{\Psi}^{F})$, and they are trivially from the kernel of $\laplcomp_{\Psi}^{F}$. Now, let $u \in \ker(\laplcomp_{\Psi}^{F})$. Then from Proposition \ref{prop_green_identity}, we deduce that 
	 		 \begin{equation}
	 		 	\scal{\laplcomp_{\Psi}^{F} u}{u}_{L^2(\Psi, F)} = \scal{\nabla^F u}{\nabla^F u}_{L^2(\Psi, F)}.
	 		 \end{equation}
	 		 Thus, we see that $\nabla^F u = 0$, which means that $u$ is a flat section. 
	 	\end{proof}
	 	
	\subsection{Domain of Friedrichs Laplacian through singularities}\label{sect_funct_an_flat}
	In this section we give a more explicit description of the domain of the Friedrichs extension of the Laplacian. This description will later play an important role in the proof of Theorem \ref{thm_eigval_convergence}. We don't claim originality on this section but again we weren't able to find a complete reference for it.
	\par
	We describe${\rm{Dom}}_{Fr}(\laplcomp_{\Psi}^{F})$ by prescribing some asymptotical behavior near ${\rm{Con}(\Psi)} \cup {\rm{Ang}(\Psi)}$ to the functions from it. All our considerations are local, so by choosing local flat frames, we may and we will assume in all the proofs that our flat vector bundle is a trivial line bundle.
	\par 
	To begin, we suppose that $\partial \Psi = \emptyset$ and that $(F, h^F, \nabla^F)$ is a trivial line bundle. 
	Then it is sufficient to consider a surface $\Psi$ with only one conical point $P$ of the conical angle $\theta = \angle(P)$. 
	
	\par For $k \in \nat$, we introduce the functions $C^k_{\pm, \theta}$ on the model cone $C_{\theta}$, (\ref{eq_v_theta}), by 
	\begin{equation}
		\begin{aligned}
			& C^k_{\pm, \theta}(r, \phi) = r^{\pm \frac{2 \pi k}{\theta}} \exp \Big(\imun \frac{2 \pi k \phi}{\theta} \Big), \qquad \text{for } k > 0,
			\\
			& C^0_{+, \theta}(r, \phi) = 1, \qquad \qquad  \qquad C^0_{-, \theta}(r, \phi) = \log(r).
		\end{aligned}
	\end{equation}
	This is easily verifiable that those functions are formal solutions to the homogeneous problem $\laplcomp_{C_{\theta}} u = 0$ on the cone $C_{\theta}$, (\ref{eq_v_theta}).
	Notice that  the functions $C^k_{\pm, \theta}$, $k \in \nat$ belong to the $L^2$-space $L^2(C_{\theta})$ of the cone $C_{\theta}$, with respect to the conical metric (\ref{eq_conical_metric}).
	\par 
	Let $\chi_P$ be a smooth function on $\Psi$ which is equal to $1$ near $P$ and such that in a vicinity of the support of $\chi_P$, the manifold $\Psi$ is isometric to $C_{\theta}$.
	We use this isometry to view the function $\chi_P C^k_{\pm, \theta}$ as a function on $\Psi$.
	\begin{sloppypar}		 
		Let ${\rm{Dom}}_{\min}(\laplcomp_{\Psi})$ denote the minimal closure of $\laplcomp_{\Psi}$, viewed as an operator on $\ccal^{\infty}_{0, vN}(\Psi)$. In other words, for $u \in L^2(\Psi)$, we have $u \in {\rm{Dom}}_{\min}(\laplcomp_{\Psi})$ if and only if there exists a sequence of functions $u_m \in \ccal^{\infty}_{0, vN}(\Psi)$ such that $u_m \to u$ in $L^2(\Psi)$ and $\laplcomp_{\Psi} u_m \to w$ for some $w \in L^2(\Psi)$. Clearly, in this case, we have $\laplcomp_{\Psi} u = w$ on the level of currents.
	\end{sloppypar}
	Recall that the maximal domain ${\rm{Dom}}_{\max}(\laplcomp_{\Psi})$ was defined in the beginning of Section \ref{sect_zeta_flat_fun}.
	\begin{prop}[{Mooers \cite[Proposition 2.3]{Mooers}} ]\label{prop_max_domain_mooers}
		The following identity holds
		\begin{equation}
			{\rm{Dom}}_{\max}(\laplcomp_{\Psi}) = {\rm{Dom}}_{\min}(\laplcomp_{\Psi}) 
			+
			 \sum_{0 \leq k < \frac{\theta}{2 \pi}} \langle \chi_P C^k_{+, \theta} \rangle 
			+
			 \sum_{0 \leq k < \frac{\theta}{2 \pi}} \langle \chi_P C^k_{-, \theta} \rangle,
		\end{equation}
		where by $\langle v \rangle$ we mean a vector space spanned by the vector $v$.
	\end{prop}
	\par 
	Now, the description of the set of all self-adjoint extensions of $\laplcomp_{\Psi}$ looks as follows, cf. Mooers \cite[Theorem 2.1]{Mooers}. 
	Denote by $\mathcal{M}$ the linear subspace of $L^2(\Psi, F)$ spanned by the functions $\chi_P C^k_{\pm, \theta}$ with $0 \leq k < \frac{\theta}{2 \pi}$.
	The dimension, $2d$, of $\mathcal{M}$ is even. 
	To get a self-adjoint extension of $\laplcomp_{\Psi}$ one chooses a subspace $\mathcal{N}$ of $\mathcal{M}$ of dimension $d$ such that for any $u, v \in \mathcal{N}$, we have
	\begin{equation}
		\scal{\laplcomp_{\Psi} u}{v}_{L^2(\Psi)} - \scal{u}{\laplcomp_{\Psi} v}_{L^2(\Psi)} 
		=
		0.
	\end{equation}
	To any such $\mathcal{N}$ there corresponds a self-adjoint extension of $\laplcomp_{\Psi}$ with domain ${\rm{Dom}}_{\min}(\laplcomp_{\Psi}) + \mathcal{N}$.
 	\par The domain of the \textit{Friedrichs extension} ${\rm{Dom}}_{Fr}(\laplcomp_{\Psi})$ corresponds to the choice of functions $\chi_P C^k_{+, \theta}$ with $0 \leq k < \frac{\theta}{2 \pi}$.
 	To verify this, by Proposition \ref{prop_max_domain_mooers} and the description above, since the number of those functions is exactly equal to $\frac{1}{2} \dim ({\rm{Dom}}_{\max}(\laplcomp_{\Psi}) / {\rm{Dom}}_{\min}(\laplcomp_{\Psi}))$, it is enough to prove that those functions are inside of ${\rm{Dom}}_{Fr}(\laplcomp_{\Psi})$. 
 	But this follows by an explicit construction of a sequence of functions $f_{\epsilon}^{k} \in \ccal^{\infty}_{0}(\Psi \setminus \{ P\})$, $\epsilon > 0$ satisfying $f_{\epsilon}^{k} \to \chi_P C^k_{+, \theta}$ and $\nabla f_{\epsilon}^{k} \to \nabla(\chi_P C^k_{+, \theta})$ in $L^2(\Psi)$, as $\epsilon \to 0$.
 	For example, take 
 	\begin{equation}\label{eq_f_eps_conv_fam}
		f_{\epsilon}^{k} =  \big( \chi_P(r, \phi) - \chi_P(r/\epsilon, \phi) \big)  C^k_{+, \theta}.
 	\end{equation}
 	\par 
 	To make further notation easier, for $P \in {\rm{Con}}(\Psi)$, we denote
 	\begin{equation}\label{eq_defn_ckp}
 		C^k_{P} := \chi_P C^k_{+, \angle(P)}.
 	\end{equation}
 	\par Now let's come back to the original surface $\Psi$ with piecewise geodesic boundary $\partial \Psi$ and a flat unitary vector bundle $(F, h^F, \nabla^F)$ over $\overline{\Psi}$. 
 	Consider a double manifold $\widetilde{\Psi} := \Psi \sqcup \Psi^*$, where $\Psi^*$ is isomorphic to $\Psi$ but with the opposite orientation, and $\Psi$ and $\Psi^*$ are glued along the boundary $\partial \Psi$ in an obvious way.
 	Then its easy to see that $\widetilde{\Psi}$ has a structure of a flat surface with conical angles $\alpha$, where $\alpha \in \angle({\rm{Con}}(\Psi))$ is of double multiplicity coming from the conical angles of $\Psi$ and $\Psi^*$, and $2 \beta$, where $\beta \in \angle({\rm{Ang}}(\Psi))$, coming from corners of the boundary. 
 	\par 
 	The manifold $\widetilde{\Psi}$ has a natural involution $i$, interchanging $\Psi$ and $\Psi^*$. The fixed point set of $i$ is equal to $\partial \Psi \hookrightarrow \widetilde{\Psi}$. We denote by $\pi : \widetilde{\Psi} \to \Psi$ the obvious projection, for which we have $\pi \circ i = \pi$.
 	The pullback $\pi^*(F, h^F, \nabla^F)$ gives a flat unitary vector bundle $(\widetilde{F}, h^{\widetilde{F}}, \nabla^{\widetilde{F}})$ over $\widetilde{\Psi}$.
 	\par 
 	Define the functions $A^k_{Q}$, $Q \in {\rm{Ang}}(\Psi)$, $k \in \nat$, by
 	\begin{equation}\label{eq_defn_akq}
 		A^k_{Q}(r, \phi) := \chi_Q (C^{k}_{+, 2 \beta}(r, \phi) + C^{k}_{+, 2 \beta}(r, -\phi)),
 	\end{equation}
 	where $\chi_Q$ is a smooth function on $\Psi$ which is equal to $1$ near $Q$ and such that in a vicinity of the support of $\chi_Q$, the manifold $\Psi$ is isometric to a neighborhood of a standard angle
 	\begin{equation}\label{eq_a_theta}
		A_{\theta} := \{(r, t) : r \geq 0; 0 \leq t \leq \theta \}/_{(0, t) \sim (0, t')},
	\end{equation}
	endowed with the metric (\ref{eq_conical_metric}).
 	\par 
 	From (\ref{eq_dom_friedr}), we see that the domain ${\rm{Dom}}_{Fr}(\laplcomp_{\Psi}^{F})$ of the Friedrichs extensions of the Laplacian $\laplcomp_{\Psi}^{F}$ on $\Psi$ with von Neumann boundary conditions on $\partial \Psi$ is related to the domain ${\rm{Dom}}_{Fr}(\laplcomp_{\widetilde{\Psi}}^{\widetilde{F}})$ of the Friedrichs extensions of the Laplacian $\laplcomp_{\widetilde{\Psi}}^{\widetilde{F}}$ as follows: $u \in {\rm{Dom}}_{Fr}(\laplcomp_{\Psi}^{F})$ if and only if $\pi^* u \in {\rm{Dom}}_{Fr}(\laplcomp_{\widetilde{\Psi}}^{\widetilde{F}})$. Moreover, the following identity holds
 	\begin{equation}\label{eq_pi_lapl_psi}
 		\pi^* (\laplcomp_{\Psi}^{F} u) = \laplcomp_{\widetilde{\Psi}}^{\widetilde{F}} (\pi^* u).
 	\end{equation}
 	From the above description of the case with empty boundary and (\ref{eq_pi_lapl_psi}), we deduce
 	\begin{prop}
 		The domain ${\rm{Dom}}_{Fr}(\laplcomp_{\Psi}^{F})$ of Friedrichs extension can be described as
 		 \begin{equation}\label{eq_dom_fr_descr}
	 		{\rm{Dom}}_{Fr}(\laplcomp_{\Psi}^{F}) = {\rm{Dom}}_{\min}(\laplcomp_{\Psi}^{F}) 
	 		+ 
	 		\sum_{P \in {\rm{Con}}(\Psi)} \sum_{0 \leq k < \frac{\angle(P)}{2 \pi} } \langle C^k_{P} \rangle
	 		+
	 		\sum_{Q \in {\rm{Ang}}(\Psi)} \sum_{0 \leq k < \frac{\angle(Q)}{\pi} } \langle A^k_{Q} \rangle,
	 	\end{equation}
	 	where by $\langle v \rangle$ we mean a vector space spanned by the vectors $v \cdot e_i$, where $e_i$, $i = 1, \ldots, \rk{F}$ is a flat local frame in the neighborhood of a fixed point over which the summation is done. 
	\end{prop}

 	\begin{prop}\label{prop_conic_angle_sobolev}
 		For $P \in {\rm{Con}}(\Psi)$ and $Q \in {\rm{Ang}}(\Psi)$, we have
 		\begin{equation}\label{eq_ckp_sobolll}
 			\begin{aligned}
 				& C^{k}_{P} \in W^l_{p}(\Psi) \quad \text{ if and only if } \quad k = 0 \quad \text{ or } \quad \frac{2 \pi k}{\angle(P)} > l - \frac{2}{p},
 				\\
 				& A^{k}_{Q} \in W^l_{p}(\Psi) \quad \text{ if and only if } \quad k = 0 \quad \text{ or } \quad \frac{\pi k}{\angle(Q)} > l - \frac{2}{p}.
 			\end{aligned}
 		\end{equation}
 		Also, for any $k \in \nat$, $P \in {\rm{Con}}(\Psi)$ and $Q \in {\rm{Ang}}(\Psi)$, we have
 		\begin{equation}\label{eq_incl_c_a_fr}
 			C^{k}_{P}, A^{k}_{Q} \in {\rm{Dom}}_{Fr}(\laplcomp_{\Psi}).
		\end{equation}	 		 
 	\end{prop}
 	\begin{proof}
 		The statement (\ref{eq_ckp_sobolll}) follows by an explicit evaluation of derivatives, cf. Grisvard \cite[Preface and Lemma 4.4.3.5]{Grisvard}.
 		The statement (\ref{eq_incl_c_a_fr}) follows by an explicit construction of the family of functions in $\ccal^{\infty}_{0, vN}(\Psi)$ converging to $C^{k}_{P}$ and $A^{k}_{Q}$ in $H^1(\Psi)$ as it was done in (\ref{eq_f_eps_conv_fam}).
 	\end{proof}
 	\begin{prop}\label{prop_dom_sobol}
 		For any $(\Psi, g^{T \Psi})$ and $(F, h^F, \nabla^F)$ as before, the following inclusion holds
 		\begin{equation}\label{prop_dom_sobol_2}
 		 	{\rm{Dom}}_{Fr}(\laplcomp_{\Psi}^{F}) \subset L^{\infty}(\Psi, F).
 		\end{equation}
 		Also, there exists $p_0 > 1$, which depends only on the sets $\angle({\rm{Con}}(\Psi))$ and $\angle({\rm{Ang}}(\Psi))$, such that
 		\begin{equation}\label{prop_dom_sobol_1}
 			{\rm{Dom}}_{Fr}(\laplcomp_{\Psi}^{F}) \subset W^2_{p_0}(\Psi, F).
 		\end{equation}
 	\end{prop}
 	\begin{rem}
 		In fact, the Friedrichs extension is the only extension satisfying (\ref{prop_dom_sobol_2}). This can be easily seen from the general description of self-adjoint extensions, given after Proposition \ref{prop_max_domain_mooers}.
 	\end{rem}
 	\begin{proof}
 		In \cite[Theorem 4.3.1.4]{Grisvard}, Grisvard proved the following a priori estimate: there exists $C > 0$ such that for any $u \in H^2(\Psi, F)$, the following elliptic estimate holds
 		\begin{equation}\label{eq_grisv_apriori}
 			\norm{u}_{H^2(\Psi, F)} \leq C 
 			\big( 
 				\norm{u}_{L^2(\Psi, F)} + \| \laplcomp_{\Psi}^{F} u \|_{L^2(\Psi, F)}
 			\big).
 		\end{equation}
 		We note that Grisvard proved this estimate only for polygons with no vector bundles, but the proof remains obviously valid for general flat surfaces with piecewise geodesic boundary and flat unitary vector bundle.
 		From (\ref{eq_grisv_apriori}), we deduce that
 		\begin{equation}\label{eq_d_min_h2}
 			{\rm{Dom}}_{\min}(\laplcomp_{\Psi}^{F}) \subset H^2(\Psi, F).
 		\end{equation}
 		From Proposition \ref{prop_conic_angle_sobolev}, we see that one can choose $p > 1$ very close to $1$ so that $C^{k}_{P} \in W^2_{p}(\Psi)$ for $P \in {\rm{Con}}(\Psi)$ and $0 \leq k < \frac{\angle(P)}{2 \pi}$
 		and $A^{k}_{Q} \in W^2_{p}(\Psi)$ for $Q \in {\rm{Ang}}(\Psi)$ and $0 \leq k < \frac{\angle(Q)}{\pi} $. By this, (\ref{eq_dom_fr_descr}) and (\ref{eq_d_min_h2}), we deduce (\ref{prop_dom_sobol_1}).
 		\par 
 		Now, by Adams \cite[Theorem 5.4 Part II, Case C'']{Adams}, the following inclusion holds
 		\begin{equation}\label{eq_dom_sobol_aux_1}
 			H^2(\Psi, F) \subset L^{\infty}(\Psi, F).
 		\end{equation}
 		We note that Adams proves (\ref{eq_dom_sobol_aux_1}) only for planar domains satisfying cone property (cf. Adams \cite[p. 66]{Adams} for a definition of the cone property), but since a flat surface can be decomposed into sectors, which satisfy the cone property, the inclusion (\ref{eq_dom_sobol_aux_1}) continues to hold for all $\Psi$.
 		Now, (\ref{prop_dom_sobol_2}) follows from (\ref{eq_dom_sobol_aux_1}) and the trivial fact that $C^{k}_{P}, A^{k}_{Q} \in L^{\infty}(\Psi)$.
 	\end{proof}
 	
	\section{Finite difference method, proofs of Theorems \ref{thm_eigval_convergence}, \ref{thm_eigvec_convergence}}\label{sect_fd_meth}

	In this section we investigate the extension of finite difference method from lattices $\frac{1}{n} \integ^2$, $n \in \nat^*$ in $\real^2$ to general graphs $\Psi_n$ “approximating" pillowcase covers $\Psi$ (see Section \ref{sect_sq_t_s_disc}).
	In particular, we prove Theorems \ref{thm_eigval_convergence}, \ref{thm_eigvec_convergence}.
	We conserve the notation from Theorems \ref{thm_eigval_convergence}.
	\par 
	This section is organized as follows.
	The main result of Section \ref{sect_eigvec_reg} is
	\begin{thm}\label{thm_bnd_fin_cont_above}
		For any $i \in \nat$, the following bound holds
		\begin{equation}\label{eq_bnd_fin_cont_above}
			\lambda_i \geq  \limsup_{n \to \infty} \lambda_i^{n}.
		\end{equation}
	\end{thm}
	In Section \ref{sect_harn_in}, modulo some Harnack-type inequality, Theorem \ref{thm_harn_type}, which we prove in Section \ref{sect_harn_pf}, we prove the “inverse" statement
	\begin{thm}\label{thm_bnd_fin_cont_below}
		For any $i \in \nat$, the following bound holds
		\begin{equation}\label{eq_bnd_fin_cont_below}
			\liminf_{n \to \infty} \lambda_i^{n} \geq  \lambda_i.
		\end{equation}
	\end{thm}		
	Clearly, Theorems \ref{thm_bnd_fin_cont_above}, \ref{thm_bnd_fin_cont_below} and Theorem \ref{thm_eigval_convergence} are equivalent. However, we prefer to state them separately, since the techniques we use in their proofs are rather different.
	\par 
	Then in Section \ref{sect_pf_eigenvect}, we show that Theorem \ref{thm_eigvec_convergence} follows almost formally from the arguments, developed in the course of the proof of Theorem \ref{thm_eigval_convergence}.

	\subsection{Regularity of eigenvectors on a flat surface, a proof of Theorem \ref{thm_bnd_fin_cont_above}}\label{sect_eigvec_reg}
		The main goal of this section is to prove Theorem \ref{thm_bnd_fin_cont_above}.
		The main idea is to start from an eigenvector $f_i$ corresponding to an eigenvalue $\lambda_i$ of the Friedrichs extensions of the Laplacian $\laplcomp_{\Psi}^{F}$, i.e.
		\begin{equation}
			\laplcomp_{\Psi}^{F} f_i = \lambda_i f_i,
		\end{equation}		 
		and to construct a section $R_n f_i \in {\rm{Map}}(V(\Psi_n), F_n)$ by restriction, (\ref{eq_inj_vert}). We prove that the Rayleigh quotients  associated to $f_i$ and $R_n f_i$ are close enough. Then Theorem \ref{thm_bnd_fin_cont_above} would follow from a simple application of the min-max theorem. 
		The main difficulty is that due to the lack of elliptic regularity near singularities, $f_i$ does not extend smoothly up to ${\rm{Con}}(\Psi) \cup {\rm{Ang}}(\Psi)$.
		\begin{sloppypar}
			More precisely, we define the restriction operator as follows
			\begin{equation}\label{eq_defn_restr_fun}
				R_n: \ccal^{0}(\Psi \setminus  \partial \Psi, F) \to {\rm{Map}}(V(\Psi_n), F_n), 
				\qquad  
				(R_n f)(v) := f(v),
			\end{equation}
			where we implicitly used the injection (\ref{eq_inj_vert}).
		\end{sloppypar}		 
		\par
		Recall that in the whole article we suppose that $\Psi$ is tiled by euclidean squares of area 1.
		The following proposition shows that in some sense, the rescaled discrete Laplacians $n^2 \cdot \laplcomp_{\Psi_n}^{F_n}$ approximate weakly the smooth Laplacian $\laplcomp_{\Psi}^{F}$ with von Neumann boundary conditions. This fact relies heavily on our choice of discretizations $\Psi_n$ and the injection (\ref{eq_inj_vert}).
		\begin{prop}\label{prop_find_diff_vn_cond}
			There is a constant $C > 0$ such that for any $f \in \ccal^{3}(\Psi, F)$, satisfying (\ref{eq_vn_bound_cond}), any $n \in \nat^*$ and any $P \in V(\Psi_n)$, the following bound holds
			\begin{equation}\label{eq_fin_diff_ccal3}
				\Big| n^2 \cdot \laplcomp_{\Psi_n}^{F_n} (R_n f)(P) -  (\laplcomp_{\Psi}^{F} f) (P) \Big| 
				\leq
				\frac{C}{n} \norm{f}_{\ccal^{3}(B_{\Psi}( \frac{2}{n}, P ))}.
			\end{equation}
		\end{prop}
		\begin{proof}
			Our considerations are purely local, so we may assume that $(F, h^F, \nabla^F)$ is trivial.
			\par 
			The are essentially three different cases to consider. Those are the cases when the vertex $P \in V(\Psi_n)$ has degree $2$, $3$ or $4$.
			We will concentrate on the most cumbersome case, which is of degree $3$. This means that the vertex $P$ lies on the “boundary" of $\Psi_n$.
			\par 
			Denote by $Q, R, S \in V(\Psi_n)$ the neighbors of $P$. Normalize the coordinates $x, y$ so that $P$ corresponds to $(0, 0)$, $Q$ to $(-\frac{1}{n}, 0)$, $R$ to $(\frac{1}{n}, 0)$ and $S$ to $(0, -\frac{1}{n})$ as in Figure \ref{fig_pqrst_bound}.
			\begin{figure}[h]
				\includegraphics[width=0.4\textwidth]{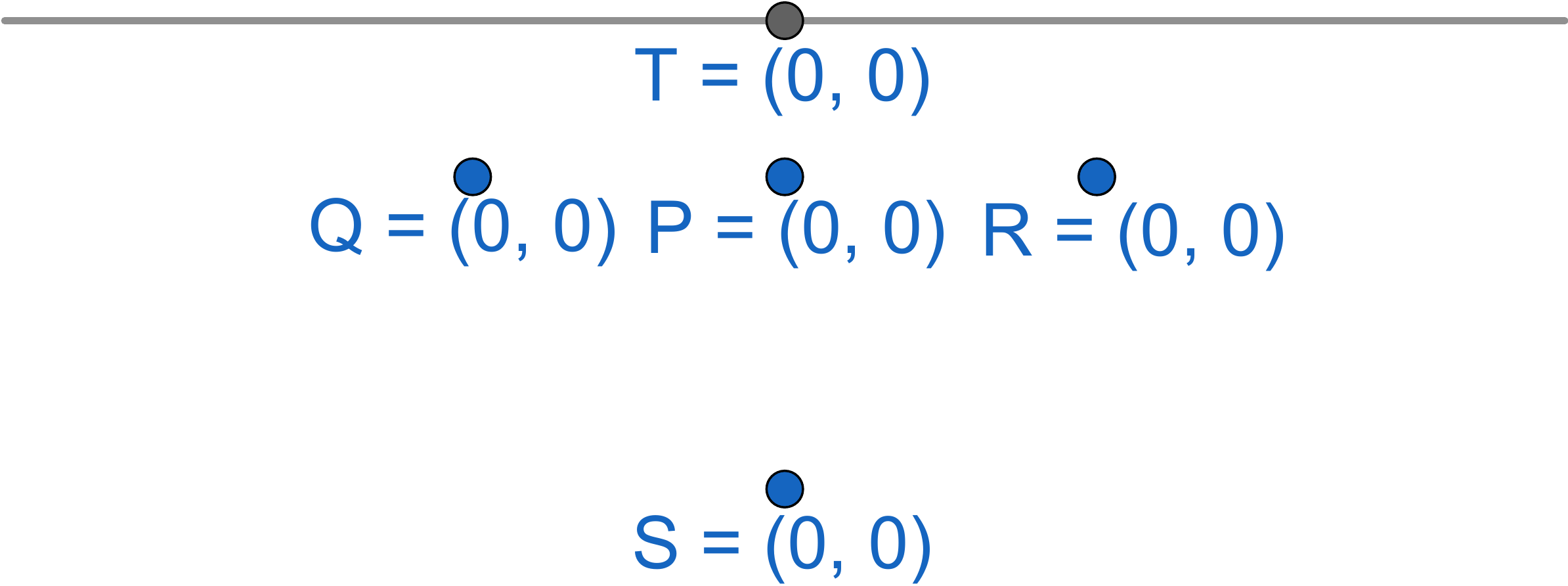}	
				\centering
				\caption{An arrangement of vertices $P, Q, R, S$. The boundary $\partial \Psi$ is in grey.}
				\label{fig_pqrst_bound}
			\end{figure}
			Denote by $T$ the point with coordinates $(0, \frac{1}{2n})$.
			Now, by Taylor expansion, we have
			\begin{align}
				& \label{eq_lapl_hor_fd} 2 f(P) - f(Q) - f(R) = - \frac{1}{n^2} \frac{\partial^2}{\partial x^2} f(P)
				+
				O\Big(
					\frac{1}{n^3} \norm{f}_{\ccal^{3}(B_{\Psi}( \frac{2}{n}, P ))}
				\Big),
				\\
				& \label{eq_taylor_bound}
				\frac{\partial}{\partial y} f(T) =  \frac{\partial}{\partial y} f(P) + \frac{1}{2n} \frac{\partial^2}{\partial y^2} f(P) + O\Big( \frac{1}{n^2} \norm{f}_{\ccal^{3}(B_{\Psi}( \frac{2}{n}, P ))} \Big).
			\end{align}
			From (\ref{eq_vn_bound_cond}) and (\ref{eq_taylor_bound}), we obtain
			\begin{equation}\label{eq_lapl_ver_fd}
				f(P) - f(S) = - \frac{1}{n^2} \frac{\partial^2}{\partial y^2} f(P) 
				+ 				
				O\Big(
					\frac{1}{n^3} \norm{f}_{\ccal^{3}(B_{\Psi}( \frac{2}{n}, P ))}
				\Big).
			\end{equation}
			From (\ref{eq_lapl_hor_fd}) and (\ref{eq_lapl_ver_fd}), we conclude.
		\end{proof}
		\par 
		Now, we consider the eigenspace $V_k$, $k \in \nat$, associated with the first $k$ eigenvalues of $\laplcomp_{\Psi}^{F}$.
		In the end of this section, by studying explicitly the structure of the singularities of $f \in V_k$ near ${\rm{Con}}(\Psi)$ and ${\rm{Ang}}(\Psi)$, we show the following
		\begin{thm}\label{thm_reg_main_limit_2}
			For any $f, g \in V_k$, as $n \to \infty$, we have
			\begin{equation}\label{eq_reg_main_limit_2}
				\scal{\laplcomp_{\Psi_n}^{F_n} (R_n f)}{R_n g}_{L^2(\Psi_n, F_n)} \to \scal{\laplcomp_{\Psi}^{F} f}{g}_{L^2(\Psi, F)}.
			\end{equation}
		\end{thm}
		Let's see how Theorem \ref{thm_bnd_fin_cont_above} follows from Theorem \ref{thm_reg_main_limit_2}.
		\begin{proof}[Proof of Theorem \ref{thm_bnd_fin_cont_above}.]
			By elliptic regularity, we know that the eigenvectors $f_i$, $i \in \nat$ of $\laplcomp_{\Psi}^{F}$ are smooth in the interior of $\Psi$.
			By this and (\ref{prop_dom_sobol_2}), we conclude
			\begin{equation}\label{eq_fi_linf}
				f_i \in \ccal^{\infty}(\Psi, F) \cap L^{\infty}(\Psi, F)
			\end{equation}
			From (\ref{eq_fi_linf}) and the fact that the tiles have area $1$, we see that for any $f, g \in V_k$, as $n \to \infty$:
			\begin{equation}\label{eq_lim_cont_linf}
				\frac{1}{n^2} \scal{R_n f}{R_n g}_{L^2(\Psi_n, F_n)} \to \scal{f}{g}_{L^2(\Psi, F)}.
			\end{equation}
			Construct a vector space $V_k^{n}$ as follows
			\begin{equation}
			V_k^{n} := R_n(V_k).
			\end{equation}
			From Theorem \ref{thm_reg_main_limit_2} and (\ref{eq_lim_cont_linf}), we deduce that
			\begin{equation}\label{eq_reg_aux_1}
				\lim_{n \to \infty} 
				\sup_{f \in V_k^{n}} 
				\bigg\{ 
					\frac{\scal{n^2 \cdot (\laplcomp_{\Psi_n}^{F_n}f)}{f}_{L^2(\Psi_n, F_n)}}{\scal{f}{f}_{L^2(\Psi_n, F_n)}}
				\bigg\}
				=
				\lambda_k.
			\end{equation}
			Now, by the characterization of the eigenvalues through Rayleigh quotient, we have
			\begin{equation}\label{eq_rayleigh_2}
				\lambda_{k}^{n} = \inf_{\substack{V \subset {\rm{Map}}(V(\Psi_n), F_n) } } \sup_{f \in V} 
				\bigg\{ 
					\frac{\scal{n^2 \cdot (\laplcomp_{\Psi_n}^{F_n} f)}{f}_{L^2(\Psi_n, F_n)}}{\scal{f}{f}_{L^2(\Psi_n, F_n)}}
					 : \dim V = k		
				\bigg\}.
			\end{equation}
			Clearly, for $n$ big enough, we have $\dim V_k^{n} = k$.
			By this, from (\ref{eq_rayleigh_2}), we see that
			\begin{equation}\label{eq_reg_aux_22}
				\lambda_{k}^{n} \leq \sup_{f \in V_k^{n}} 
				\bigg\{ 
					\frac{\scal{n^2 \cdot \laplcomp_{\Psi_n}^{F_n}(f)}{f}_{L^2(\Psi_n, F_n)}}{\scal{f}{f}_{L^2(\Psi_n, F_n)}}
				\bigg\}
			\end{equation}
			By (\ref{eq_reg_aux_1}) and (\ref{eq_reg_aux_22}), we deduce Theorem \ref{thm_bnd_fin_cont_above}.
		\end{proof}
		To finish the proof of Theorem \ref{thm_bnd_fin_cont_above}, we only need to prove Theorem \ref{thm_reg_main_limit_2}.
		The main difficulty lies in the fact that in general, the finite differences $n^2 \cdot \laplcomp_{\Psi_n}^{F_n} R_n(f_i)$ might “explode" near ${\rm{Con}}(\Psi) \cup {\rm{Ang}}(\Psi)$. We use the theory of elliptic regularity for polygons developed by Grisvard, \cite{Grisvard}, to prove that this divergence poses no problem once we are concerned with the $L^2$-product in (\ref{eq_reg_main_limit_2}).
		\par
		Recall that the functions $C_P^{k}$ and $A_Q^{k}$ for $P \in {\rm{Con}}(\Psi)$, $Q \in {\rm{Ang}}(\Psi)$ and $k \in \nat$ were defined in (\ref{eq_defn_ckp}) and (\ref{eq_defn_akq}).
		The main result of this section is the following 
		\begin{thm}\label{thm_reg_eig}
			There are coefficients $\gamma_P^{k, j} \in \real$ for $P \in {\rm{Con}}(\Psi)$, $k \in \nat$, $0 < \frac{2 \pi k}{\angle(P)}  \leq 4 - \frac{2}{p}$, $1 \leq j \leq \rk{F}$ and $\alpha_Q^{k, j} \in \real$ for $Q \in {\rm{Ang}}(\Psi)$, $k \in \nat$, $0 < \frac{\pi k}{\angle(Q)}  \leq 4 - \frac{2}{p}$, $1 \leq j \leq \rk{F}$ such that the section $g_i$, defined as
			\begin{equation}\label{thm_reg_eig_111}
				g_i := f_i - \sum_{P \in {\rm{Con}}(\Psi)} \sum_{0 <  \frac{2 \pi k}{\angle(P)}  \leq 4 - \frac{2}{p}} \gamma_P^{k, j} \cdot e_j \cdot C_P^{k}  - \sum_{Q \in {\rm{Ang}}(\Psi)} \sum_{0 <  \frac{\pi k}{\angle(Q)}  \leq 4 - \frac{2}{p}} \alpha_Q^{k, j} \cdot e_j \cdot A_Q^{k},
			\end{equation}
			satisfies the following regularity properties (see (\ref{eq_defn_ck_b}))
			\begin{equation}
				g_i \in \ccal^{2}(\overline{\Psi}, F).
			\end{equation}
			In (\ref{thm_reg_eig_111}) we used Einstein summation convention for local flat frames $e_j$ based at respective points.
		\end{thm}
		Let's now see how Theorem \ref{thm_reg_main_limit_2} follows from Theorem \ref{thm_reg_eig}. After that, we prove Theorem \ref{thm_reg_eig}.
		\begin{proof}[Proof of Theorem \ref{thm_reg_main_limit_2}.]
			By elliptic regularity, the functions $f$, $g$ are smooth away from ${\rm{Con}}(\Psi) \cup {\rm{Ang}}(\Psi)$. So, by Proposition \ref{prop_find_diff_vn_cond}, as $n \to \infty$, the following convergence holds uniformly away from a neighborhood of ${\rm{Con}}(\Psi) \cup {\rm{Ang}}(\Psi)$:
			\begin{equation}
			\begin{aligned}
				& n^2 \cdot \laplcomp_{\Psi_n}^{F_n} (R_n f)(P_n) \to \laplcomp_{\Psi}^{F} f(P), 
				\\
				& (R_n g)(P_n) \to g(P),
			\end{aligned}
			\end{equation}
			here $P_n \in V(\Psi_n)$ converges to $P \in \Psi$, and we implicitly used parallel transport associated with $\nabla^F$.
			From this, it's clear that to prove (\ref{eq_reg_main_limit_2}), it is enough to prove that for any $\epsilon > 0$, there is $c > 0$ such that for any $P \in {\rm{Con}}(\Psi) \cup {\rm{Ang}}(\Psi)$, we have
			\begin{align}
				& \scal{\laplcomp_{\Psi}^{F} f_i}{f_j}_{L^2(B_{\Psi}(c, P))} \leq \epsilon.\label{eq_reg_aux_3}
				\\ 
				& \scal{\laplcomp_{\Psi_n}^{F_n} (R_n f_i)}{R_n f_j}_{L^2(B_{\Psi_n}(cn, P))} \leq \epsilon  \label{eq_reg_aux_2}
			\end{align}
			The bound (\ref{eq_reg_aux_3}) clearly follows from Theorem \ref{thm_reg_eig} and the fact that near ${\rm{Con}}(\Psi) \cup {\rm{Ang}}(\Psi)$, we have $\laplcomp_{\Psi} C_P^{k} = \laplcomp_{\Psi} A_Q^{k} = 0$.
			Let's now prove (\ref{eq_reg_aux_2}). 
			In fact, by (\ref{eq_fin_diff_ccal3}) and the identity $\laplcomp_{\Psi} C_P^{k} = 0$, which holds near $P$, we see by (\ref{eq_fin_diff_ccal3}) that there is $C > 0$ such that for any $n \in \nat^*$, we have
			\begin{equation}\label{eq_reg_aux_4}
				\Big| \laplcomp_{\Psi_n} (R_n C_P^{k})(Q) \Big| \leq \frac{C}{n^3} \Big( \frac{\dist_{\Psi_n}(P, Q)}{n} \Big)^{\frac{2 \pi k}{\angle(P)} - 3}.
			\end{equation}
			Thus, we conclude that for any $k, l \in \nat$, and $\alpha := \frac{2 \pi k}{\angle(P)} + \frac{2 \pi l}{\angle(P)}$, we have
			\begin{equation}\label{eq_reg_aux_5}
				\scal{\laplcomp_{\Psi_n} (R_n C_P^{k})}{R_n C_{P}^{l}}_{L^2(B_{\Psi_n}(cn, P))}  
				\leq 
				\frac{1}{n^{\alpha}} \sum_{Q \in B_{\Psi_n}(cn, P)} \dist_{\Psi_n}(P, Q)^{\alpha - 3}.
			\end{equation}
			However, trivially, for any $\epsilon > 0$, there is $c > 0$ such that
			\begin{equation}\label{eq_reg_aux_6}
				 \sum_{Q \in B_{\Psi_n}(cn, P)} \dist_{\Psi_n}(P, Q)^{\alpha - 3}
				\leq 
				4 \frac{\angle(P)}{2 \pi} \sum_{i = 0}^{cn} i^{\alpha - 2} \leq \epsilon n^{\alpha}.
			\end{equation}
			Now, by the fact that $g_i \in \ccal^{2}(\overline{\Psi})$, we conclude that there is $C > 0$ such that
			\begin{equation}\label{eq_reg_aux_7}
				\big\| 
					n^2 \cdot \laplcomp_{\Psi_n}^{F_n} (R_n g_i)
				\big\|_{L^{\infty}(\Psi, F)}
				\leq 
				C.
			\end{equation}
			By (\ref{eq_reg_aux_5}),  (\ref{eq_reg_aux_6}) and analogous estimates for $A_Q^{k}$, $Q \in {\rm{Ang}}(\Psi)$, $k \in \nat$, (\ref{eq_reg_aux_7}), we deduce that (\ref{eq_reg_aux_2}) holds. This finishes the proof.
		\end{proof}
		\begin{proof}[Proof of Theorem \ref{thm_reg_eig}.]
			By Proposition \ref{prop_dom_sobol}, there is $p > 1$ such that we have
			\begin{equation}
				f_i \in {\rm{Dom}}_{Fr}(\laplcomp_{\Psi}^{F}) \subset W^2_{p}(\Psi, F).
			\end{equation}
			Then, by Grisvard \cite[Theorem 5.1.3.1]{Grisvard}, we know that there exists a formal solution $u \in L^2(\Psi, F)$ of the equation 
			\begin{equation}
				\laplcomp_{\Psi}^{F} u = \lambda_i f_i,
			\end{equation}
			satisfying von Neumann boundary conditions on $\partial \Psi$ and for which there are coefficients $\gamma_P^{k, j} \in \real$ for $P \in {\rm{Con}}(\Psi)$, $k \in \nat$, $0 < \frac{2 \pi k}{\angle(P)}  \leq 4 - \frac{2}{p}$, $1 \leq j \leq \rk{F}$ and $\alpha_Q^{k, j} \in \real$ for $Q \in {\rm{Ang}}(\Psi)$, $k \in \nat$, $0 < \frac{\pi k}{\angle(Q)}  \leq 4 - \frac{2}{p}$, $1 \leq j \leq \rk{F}$ such that
			\begin{equation}\label{eq_u_reg_coeff}
				u - \sum_{P \in {\rm{Con}}(\Psi)} \sum_{0 <  \frac{2 \pi k}{\angle(P)}  \leq 4 - \frac{2}{p}} \gamma_P^{k, j} \cdot e_j \cdot C_P^{k}  - \sum_{Q \in {\rm{Ang}}(\Psi)} \sum_{0 <  \frac{\pi k}{\angle(Q)}  \leq 4 - \frac{2}{p}} \alpha_Q^{k, j} \cdot e_j \cdot A_Q^{k} \in W^4_{p}(\Psi, F),
			\end{equation}
			where we used Einstein summation convention for local flat frames $e_j$ based at respective points.
			\par We remark that Grisvard's result was proved for polygons, and not for flat surfaces.
			Let's explain why it continues to hold for flat surfaces as well.
			First, since Grisvard never actually uses the restriction $< 2 \pi$ on the angles of the polygon (see \cite[Remark 4.3.2.7 and \S 4.2]{Grisvard}), his result applies to any flat surface with piecewise geodesic boundary. 
			Now, one could obtain a cone by gluing an angle with Dirichlet boundary conditions and another angle with von Neumann boundary conditions (cf. for example \cite[\S 4]{MazzRow}). 
			Since the result of Grisvard was proved by using local techniques (see in particular \cite[p. 252, 253]{Grisvard}) and it holds with \textit{arbitrary boundary conditions}, by decomposing $f_i$ into von Neumann and Dirichlet parts, we see that his proof still holds for general flat surfaces with conical singularities and piecewise geodesic boundary.
			\par 
			Recall that by Sobolev embedding theorem (cf. Adams \cite[Part II, Theorem 6.2]{Adams}), we have
			\begin{equation}\label{eq_inj_w4_h1}
				W^4_{p}(\Psi, F) \subset H^1(\Psi, F).
			\end{equation}
			Remark that (\ref{eq_inj_w4_h1}) was proved in \cite{Adams} only for domains satisfying cone property and without vector bundles. However, similarly to the proof of Theorem \ref{thm_rell_kondr}, by decomposing $\Psi$ into sectors, which obviously satisfy cone property, we see that (\ref{eq_inj_w4_h1}) holds.
			\par 
			By Proposition \ref{prop_conic_angle_sobolev}, (\ref{eq_u_reg_coeff}) and (\ref{eq_inj_w4_h1}), we see that $u \in H^1(\Psi, F)$. However, by Grisvard \cite[Lemma 4.4.3.1]{Grisvard}, the solutions of $\laplcomp_{\Psi}^{F} u = \lambda_i f_i$, $u \in H^1(\Psi, F)$, satisfying von Neumann boundary conditions are unique up an element from $\ker \laplcomp_{\Psi}^{F}$. 
			By (\ref{eq_dom_friedr}), however, we see that $f_i \in H^1(\Psi, F)$, thus, by the argument above, there is an element $k \in \ker \laplcomp_{\Psi}^{F}$ such that 
			\begin{equation}\label{eq_u_fi_rel}
				u = f_i + k.
			\end{equation}						
			Now, by Sobolev embedding theorem, cf. \cite[p.97]{Adams}, for $p > 1$, the following inclusion holds
			\begin{equation}\label{eq_sob_w4_emb}
				W^4_{p}(\Psi, F) \hookrightarrow \ccal^{2}(\overline{\Psi}, F).
			\end{equation}
			Again, Adams in \cite[p.97]{Adams} gives a proof of the inclusion (\ref{eq_sob_w4_emb}) for domains verifying cone property, but by decomposing $\Psi$ into sectors which clearly satisfy this property, the inclusion (\ref{eq_sob_w4_emb}) extends to the general case as well.
			We conclude by (\ref{eq_u_reg_coeff}), (\ref{eq_u_fi_rel}) and (\ref{eq_sob_w4_emb}).
		\end{proof}
		
	\subsection{Linearization functional, a proof of Theorem \ref{thm_bnd_fin_cont_below}}\label{sect_harn_in}
		The main goal of this section is to prove Theorem \ref{thm_bnd_fin_cont_below}.
		The main idea is similar to the one from Section \ref{sect_eigvec_reg}, but the methods are crucially different.
		We start from an eigenvector $f_i^{n}$ corresponding to an eigenvalue $\lambda_i^{n}$ of the discrete model $n^2 \cdot \laplcomp_{\Psi_n}^{F_n}$ and construct by “linearization" $L_n ( f_i^{n} ) \in L^2(\Psi, F)$.
		We prove that the associated Rayleigh quotients are close enough. Then Theorem \ref{thm_bnd_fin_cont_below} would follow from a simple application of the min-max theorem.
		\par We note that all the statements of this section are local in nature, so we will constantly write them in local flat frames of $F$ without saying it explicitly. 
		All the final objects can be written using parallel transport and they do not depend on the choice of the frame.
		\begin{sloppypar}
		Let's define the “linearization" functional $L_n : {\rm{Map}}(V(\Psi_n), F_n) \to \ccal^0(\Psi, F)$.
		Recall that the sets $V_n(P)$ and $U_n(P)$ were defined in (\ref{eq_defn_vp}), (\ref{eq_defn_up}).
		First, for $f \in {\rm{Map}}(V(\Psi_n), F_n)$, let's define $f^{avg} \in {\rm{Map}}(V(\Psi_n), F_n)$ by averaging the function $f$ on $V_n(P)$, i.e.
		\begin{equation}\label{eq_favg_defn}
			f^{avg}(v) :=
			 \begin{cases} 
      			\hfill \frac{1}{\# V_n(P)} \sum_{Q \in V_n(P)} f(Q),  & \text{ if there is $P \in {\rm{Con}}(\Psi) \cup {\rm{Ang}}(\Psi), v \in V_n(P)$, }  \\
      			\hfill f(v),   & \text{ otherwise. } \\
 			\end{cases}
		\end{equation}
		By definition, the functional $L_n$ satisfies the following property
		\begin{equation}
			L_n(f) = L_n(f^{avg}).
		\end{equation}
		Now, we define $L_n$ by describing explicitly the value of $(L_n f)(z)$ for any $z \in \Psi$ and any $f \in {\rm{Map}}(V(\Psi_n), F_n)$, satisfying $f = f^{avg}$.
		\par 
		Suppose $z \in U_n(P)$ for some $P \in {\rm{Con}}(\Psi) \cup {\rm{Ang}}(\Psi)$. We define 
		\begin{equation}
			(L_n f)(z) := f(Q), \quad \text{where} \quad Q \in V_n(P).
		\end{equation}
		By the assumption, we have $f = f^{avg}$, so $(L_n f)(z)$ doesn't depend on the choice of $Q$.
		\par 
		Next, suppose that ${\rm{dist}}(z, \partial \Psi) < \frac{1}{2n}$ and that for any $P \in {\rm{Con}}(\Psi) \cup {\rm{Ang}}(\Psi)$, we have $z \notin U_n(P)$. 
		Then let $P, Q \in V(\Psi_n)$ be the two closest points to $z$ in $V(\Psi_n)$. 
		In case if there are several choices, take any.
		The geometrical place of points $z$, satisfying ${\rm{dist}}(z, \partial \Psi) < \frac{1}{2n}$, and having $P$ and $Q$ as their closest points in $V(\Psi_n)$, is a rectangle. We denote this rectangle by $R$.
		The points $P$ and $Q$ have either the same $x$ or $y$ coordinates, where $x$ and $y$ are linear coordinates having axes parallel to the boundaries of tiles of $\Psi$. 
		Suppose that they share the same $y$ coordinate. 
		We renormalize $x$ coordinate so that it satisfies $x(P) = 1$ and $x(Q) = 0$. Then we define 
		\begin{equation}\label{eq_lin_rect_defn}
			(L_n f)(z) := f(P)x(z) + f(Q)(1-x(z)).
		\end{equation} 
		\par
\begin{figure}[h]
		\includegraphics[width=0.2\textwidth]{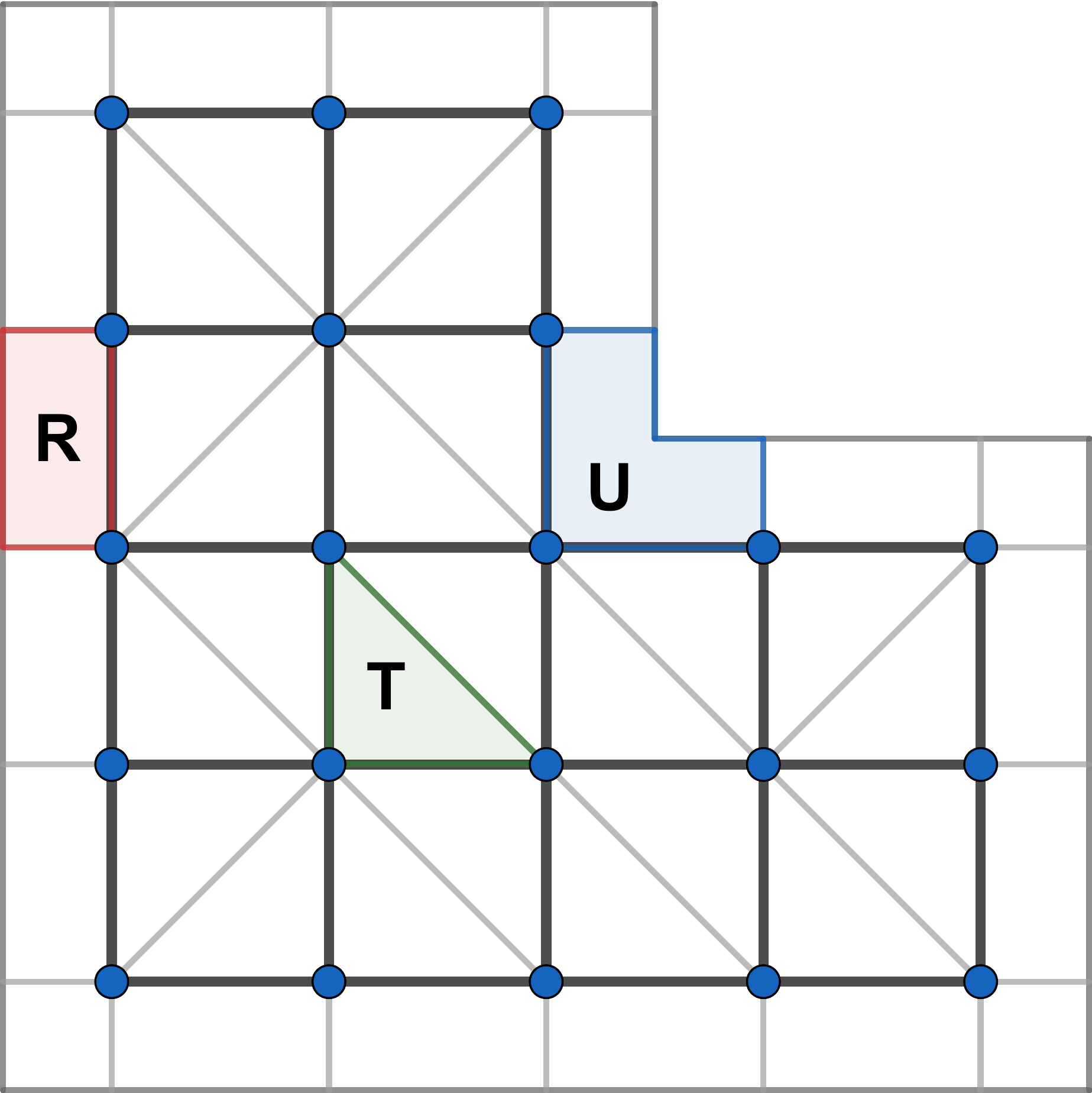}	
		\centering
		\caption{The functions from the image of $L_n$ are linear inside each subdomain. Over $U$, they are constant and over $R$ they are constant in the direction perpendicular to the boundary.}
		\label{fig_ln_func}
		\end{figure}

		Finally, suppose that $z$ is in none of the cases considered above. 
		Consider some triangulation of $\Psi \setminus (B_{\Psi}(\frac{1}{2n}, \partial \Psi) \cup (\cup_{P \in {\rm{Con}}(\Psi) \cup {\rm{Ang}}(\Psi)} U_n(P)))$ with vertices at $V(\Psi_n)$.
		Let $P, Q, R \in V(\Psi_n)$ be the vertices of the triangle $T$ containing $z$. 
		If there are several choices, take any. 
		We define $(L_n f)(z)$ as the value of the unique linear function $L$ at $z$, satisfying $L(P) = f(P)$, $L(Q) = f(Q)$ and $L(R) = f(Q)$ (with respect to the coordinates $x, y$ as in the previous step).
		This procedure describes $L_n$ completely.
		For a schematic description of the functional $L_n$, see Figure \ref{fig_ln_func}.
		\par Let's explain the construction of $L_n$. The choice of the value of $L_nf(z)$ for $z$ in the interior of $\Psi_n$ is very natural, as it is just a linear interpolation.
		\par
\begin{figure}[h]
			\includegraphics[width=0.3\textwidth]{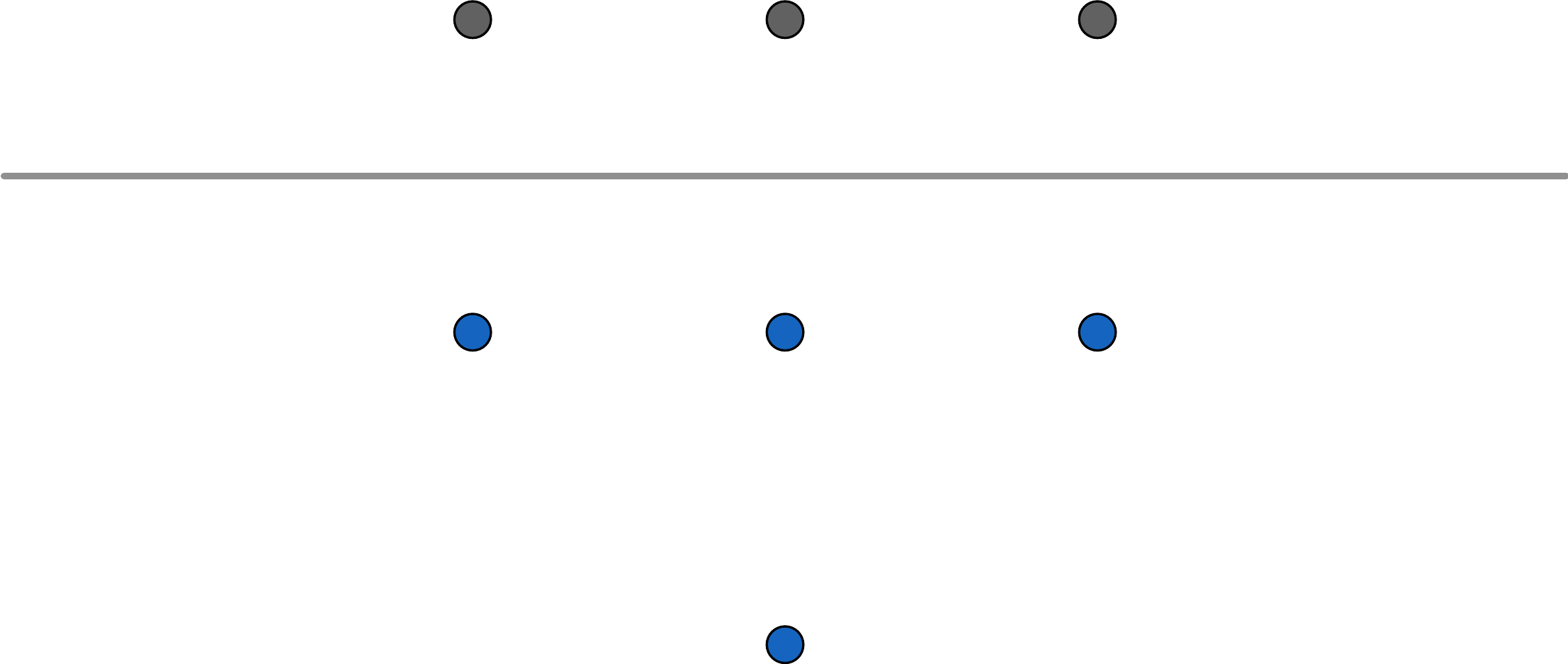}	
			\centering
			\caption{The extension of the graph $\Psi_n$. The line represents $\partial \Psi$. The empty vertices are from $\tilde{\Psi}_n$.}
			\label{fig_ext_graph}
		\end{figure}

 
		Now, the construction of $L_nf(z)$ for $z$ near the boundary of $\Psi$ can be interpreted as a linear interpolation as well, but done on a bigger graph. 
		More precisely, Proposition \ref{prop_find_diff_vn_cond} shows that in a certain sense the graph Laplacian is an approximation of the smooth Laplacian with von Neumann boundary conditions. 
		Now, let's suppose that the graph $\Psi_n$ is a part of some graph $\tilde{\Psi}_n$, which extends $\Psi$ beyond $\partial \Psi$, see Figure \ref{fig_ext_graph}.
		Near the boundary $\partial \Psi$ there is a locally defined involution with fixed point set $\partial \Psi$. 
		Extend $f \in {\rm{Map}}(V(\Psi_n), F_n)$ to $\tilde{f} \in {\rm{Map}}(V(\tilde{\Psi}_n), F_n)$  in such a way so that $\tilde{f}$ is invariant under this involution.
		It is a trivial but basic verification to see that for $P \in V(\Psi_n)$, we have $\laplcomp_{\Psi_n} f(P) = \laplcomp_{\tilde{\Psi}_n} \tilde{f}(P)$.
		But then if we treat the point $P$ as an interior point in $\tilde{\Psi}_n$ and constructs $L_n(f)$ by linear interpolation of $\tilde{f}$, as we did for regions of type $T$ in Figure \ref{fig_ln_func}, we get precisely  (\ref{eq_lin_rect_defn}).
		\par
\begin{figure}[h]
			\includegraphics[width=0.2\textwidth]{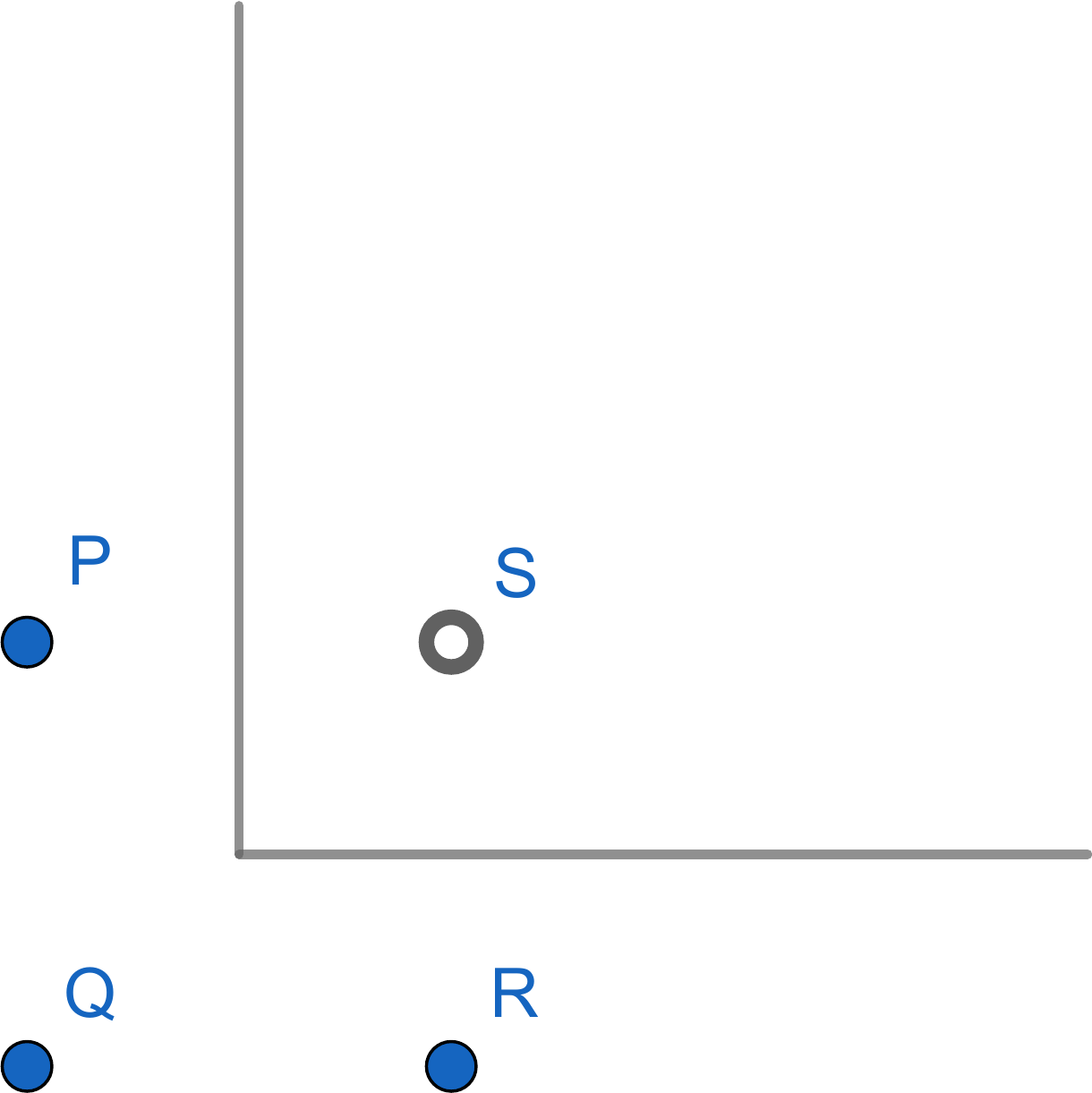}	
			\centering
			\caption{How to “extend" the function $f$ in $S$? Should we take the same value as $P$, $Q$ or $R$?}
			\label{fig_angle_problem}
		\end{figure}

 
		Now, if there is $Q \in {\rm{Con}}(\Psi) \cup {\rm{Ang}}(\Psi)$ such that $P \in V_n(Q)$, we encounter a problem with our previous explanation, depicted in Figure \ref{fig_angle_problem}. Basically, it is not clear how to extend the function $\tilde{f}$ over a point exterior to the domain.
		This is why $L_n$ averages $f$ near ${\rm{Con}}(\Psi) \cup {\rm{Ang}}(\Psi)$. 
		\par 
		The “linearization" functional satisfies a number of important properties. 
		To describe one of them, denote by ${\rm{Dom}}_{\min}^{vN}(\nabla^F)$ the minimal extension of the exterior derivative twisted by $F$, i.e.  
		\begin{equation}
			\nabla^F : \ccal^{\infty}_{0, vN}(\Psi, F) \to \ccal^{\infty}_{0, vN}(\Psi, T^* \Psi \otimes F).		
		\end{equation}
		In other words, a function $g \in L^2(\Psi, F)$ lies in ${\rm{Dom}}_{\min}^{vN}(\nabla^F)$ if and only if there is a sequence $f_n \in \ccal^{\infty}_{0, vN}(\Psi, F)$ such that $f_n \to g$ and $\nabla^F(f_n) \to h$ in $L^2(\Psi, F)$ for some $h \in L^2(\Psi, T^* \Psi \otimes F)$. If $g$ and $h$ are as described, in the distributional sense, we have $\nabla^F g = h$. Then
		\begin{prop}\label{prop_ln_dom_min}
			For any $f \in  {\rm{Map}}(V(\Psi_n), F_n)$, we have
			\begin{equation}\label{eq_dom_ln}
				L_n(f) \in {\rm{Dom}}_{\min}^{vN}(\nabla^F).
			\end{equation}
		\end{prop}
		The proof of Proposition \ref{prop_ln_dom_min} is rather direct, and it relies only the density results of smooth functions inside the Sobolev space $W^{1}_{\infty}$. It will be given in the end of this section.
		\end{sloppypar}
		\begin{prop}\label{prop_dmin_equal}
			Suppose that $f, g \in {\rm{Map}}(V(\Psi_n), F_n)$ satisfy $f = f^{avg}, g = g^{avg}$. Then the following identity holds
			\begin{equation}
				\scal{ \nabla^F \big(L_n(f) \big) }{ \nabla^F \big( L_n(g) \big) }_{L^2(\Psi, F)}  = \scal{\nabla^{F_n}_{\Psi_n} f}{\nabla^{F_n}_{\Psi_n} g}_{L^2(\Psi_n, F_n)}.
			\end{equation}
		\end{prop}
		The proof of Proposition \ref{prop_dmin_equal} is a direct calculation, and it is given in the end of this section.
		\par 
		Now, for an eigenvector $f_i^{n} \in {\rm{Map}}(V(\Psi_n), F_n)$, $i \in \nat$, $n \in \nat^*$, corresponding to the eigenvalue $\lambda_{i}^{n}$ of $n^2 \cdot \laplcomp_{\Psi_n}^{F_n}$, the functional $L_n$ satisfies the following proposition, the proof of which uses Theorem \ref{thm_bnd_fin_cont_above} and will be given in the end of this section.
		\begin{prop}\label{prop_ln_appr_easy}
			For any $\phi \in \ccal^{1}(\overline{\Psi})$, $i, j \in \nat$ fixed, as $n \to \infty$, the following estimation holds
			\begin{equation}\label{eq_ln_eig_appr_easy}
				\scal{\phi L_n(f_i^{n}) }{L_n(f_j^{n}) }_{L^2(\Psi, F)} = \frac{1}{n^2} \scal{\phi \cdot f_i^{n} }{f_j^{n}}_{L^2(\Psi_n, F_n)} + o(1).
			\end{equation}
		\end{prop}
		\begin{rem}
			In this article we will only apply Proposition \ref{prop_ln_appr_easy} for $\phi = 1$, but for further references we state it more generally.
		\end{rem}
		\par
		Now we can state the most important result of this section.
		\begin{thm}\label{thm_ln_eig_appr}
			For any $i, j \in \nat$ fixed, as $n \to \infty$, the following estimation holds
			\begin{equation}\label{eq_ln_eig_appr}
				\scal{ \nabla^F \big(L_n(f_i^{n}) \big) }{ \nabla^F \big( L_n(f_j^{n}) \big) }_{L^2(\Psi, F)} 
				= 
				\scal{ \nabla^{F_n}_{\Psi_n} f_i^{n} }{\nabla^{F_n}_{\Psi_n} f_j^{n} }_{L^2(\Psi_n, F_n)} + o(1).
			\end{equation}
		\end{thm}
		\begin{sloppypar}
		Our proof of Theorem \ref{thm_ln_eig_appr} relies on the following technical statement, to the proof of which we devote a separate Section \ref{sect_harn_pf}.
		\begin{thm}[Harnack-type inequality]\label{thm_harn_type}
			We fix $\lambda > 0$. Suppose that a sequence $f_n \in {\rm{Map}}(V(\Psi_n), F_n)$, $n \in \nat^*$, $\norm{f_n}_{L^2(\Psi_n, F_n)}^{2} = n^2$, satisfies
			\begin{equation}\label{eq_harn_bound_est}
				\big| \laplcomp_{\Psi_n}^{F_n} f_n \big| \leq \frac{\lambda}{n^2} f_n,
			\end{equation}
			Then, as $n \to \infty$, the following limit holds
			\begin{equation}\label{eq_harn_type_main_res}
				\max_{(P, Q) \in E(\Psi_n)} \big| f_n(P) - f_n(Q) \big| \to 0.
			\end{equation}
		\end{thm}
		\begin{rem}
			Essentially this theorem says that if a sequence of discrete functions is “nearly harmonic", then asymptotically this sequence is “continuous".
		\end{rem}
		\end{sloppypar}
		Let's see how Theorem \ref{thm_harn_type} and Proposition \ref{prop_dmin_equal} imply Theorem \ref{thm_ln_eig_appr}.
		\begin{proof}[Proof of Theorem \ref{thm_ln_eig_appr}]
			Recall that for  $f \in {\rm{Map}}(V(\Psi_n), F_n)$, we defined $f^{avg}$ in (\ref{eq_favg_defn}).
			By Proposition \ref{prop_dmin_equal}, we conclude
			\begin{equation}\label{eq_aux_thm_harn_1}
				\scal{ \nabla^F  \big(L_n(f_{i}^{n, avg}) \big) }{ \nabla^F  \big( L_n(f_{j}^{n, avg}) \big) }_{L^2(\Psi, F)}  = \scal{\nabla^{F_n}_{\Psi_n} (f_{i}^{n, avg} )}{\nabla^{F_n}_{\Psi_n} ( f_{j}^{n, avg} )}_{L^2(\Psi_n, F_n)}.
			\end{equation}
			However, by the construction of $f_{i}^{n, avg}$, there is $C > 0$, which depends purely on the sets $\angle({\rm{Con}}(\Psi))$ and $\angle({\rm{Ang}}(\Psi))$, such that we have
			\begin{multline}\label{eq_aux_thm_harn_2}
				\Big|	\scal{\laplcomp_{\Psi_n}^{F_n} (f_{i}^{n, avg})}{f_{j}^{n, avg}}_{L^2(\Psi_n)}
				- 
				\scal{\laplcomp_{\Psi_n}^{F_n} (f^n_{i})}{f^n_{j}}_{L^2(\Psi_n)}
				\Big| 
				\\
				\leq  
				C  \sum_{P \in {\rm{Con}}(\Psi) \cup {\rm{Ang}}(\Psi)}
				 \max_{\substack{R, Q \in V(\Psi_n), \\ {\rm{dist}}_{\Psi_n}(\{R, Q\}, V_n(P)) \leq 2}}   \big| f^n_{i}(R) - f^n_{i}(Q) \big| \cdot  \big| f^n_{j}(R) - f^n_{j}(Q) \big|,
			\end{multline}
			\vspace*{-0.3cm}
			\begin{multline}\label{eq_aux_thm_harn_3}
				\Big|
				\scal{ \nabla^F \big(L_n(f_{i}^{n, avg}) \big) }{  \nabla^F \big( L_n(f_{j}^{n, avg}) \big) }_{L^2(\Psi, F)} 
				- 
				\scal{ \nabla^F \big(L_n(f^n_{i}) \big) }{ \nabla^F \big( L_n(f^n_{j}) \big) }_{L^2(\Psi, F)} 
				\Big|
				\\
				\leq  
				C  \sum_{P \in {\rm{Con}}(\Psi) \cup {\rm{Ang}}(\Psi)}
				 \max_{\substack{R, Q \in V(\Psi_n), \\ {\rm{dist}}_{\Psi_n}(\{R, Q\}, V_n(P)) \leq 2}}   \big| f^n_{i}(R) - f^n_{i}(Q) \big| \cdot  \big| f^n_{j}(R) - f^n_{j}(Q) \big|.
			\end{multline}
			By Theorem \ref{thm_harn_type} and (\ref{eq_aux_thm_harn_1}), (\ref{eq_aux_thm_harn_2}), (\ref{eq_aux_thm_harn_3}), we conclude.
		\end{proof}
		
		Finally, after all the preparations, we are ready to prove
		\begin{proof}[Proof of Theorem \ref{thm_bnd_fin_cont_below}.]
		By Theorem \ref{thm_ln_eig_appr} and Proposition \ref{prop_ln_appr_easy}, applied for $\phi = 1$, we see that for any $k \in \nat$ fixed, as $n \to \infty$, we have
		\begin{multline}\label{eq_from_disc_to_ln}
			\lambda_k^{n}			
			=			
			\frac{\scal{n^2 \cdot \laplcomp_{\Psi_n}^{F_n}(f_k^{n})}{f_k^{n}}_{L^2(\Psi_n, F_n)}}{\scal{f_k^{n}}{f_k^{n}}_{L^2(\Psi_n, F_n)}}
			=
			\frac{n^2 \cdot \scal{\nabla^{F_n}_{\Psi_n}  (f_k^{n})  }{\nabla^{F_n}_{\Psi_n}  (f_k^{n}) }_{L^2(\Psi_n, F_n)}}{\scal{f_k^{n} }{f_k^{n}}_{L^2(\Psi_n, F_n)}}
			\\
			=
			\frac{\scal{ \nabla^F \big(L_n(f_k^{n}) \big) }{ \nabla^F \big( L_n(f_k^{n}) \big) }_{L^2(\Psi, F)} }{\scal{L_n(f_k^{n}) }{L_n(f_k^{n}) }_{L^2(\Psi, F)}}
			+ o(1).
		\end{multline}
		By Proposition \ref{prop_ln_dom_min}, we conclude that for any $i \in \nat$, $n \in \nat^*$, there is a sequence of functions $g^n_{i, k} \in \ccal^{\infty}_{0, vN}(\Psi, F)$, $k \in \nat$ such that in $L^2(\Psi, F)$, as $k \to \infty$, we have
		\begin{equation}
			g^n_{i, k} \to L_n(f_i^{n}), \qquad \qquad \nabla^F(g^n_{i, k}) \to \nabla^F \big( L_n(f_i^{n}) \big).
		\end{equation}
		Moreover, since the functions $L_n(f_i^{n})$ are constant in the neighborhood of ${\rm{Con}}(\Psi) \cup {\rm{Ang}}(\Psi)$, we see that the functions $g^n_{i, k}$ can be chosen to be constant as well.
		In particular, we have
		\begin{equation}\label{eq_gij_fr_dom_sq}
			g^n_{i, k} \in {\rm{Dom}}_{Fr}(\laplcomp_{\Psi}^{F}).
		\end{equation}
		Then for any $n \in \nat^*$, we can choose $k_{n, i}$ such that, as $n \to \infty$, the function $g^n_{i} := g_n^{i, k_{n, i}}$ satisfies
		\begin{equation}\label{eq_dmin_ln_dmin_g_n}
		\begin{aligned}
			&
			\scal{ \nabla^F \big(L_n(f_i^{n}) \big) }{ \nabla^F \big( L_n(f_j^{n}) \big) }_{L^2(\Psi, F)} 
			 = 
			\scal{ \nabla^F g^n_{i} }{ \nabla^F  g^n_{j} }_{L^2(\Psi, F)}  + o(1),
			\\
			&
			\scal{L_n(f_i^{n}) }{L_n(f_j^{n}) }_{L^2(\Psi, F)} = \scal{g^n_{i} }{g^n_{j} }_{L^2(\Psi, F)} + o(1).
		\end{aligned}
		\end{equation}
		Note, however, that by Proposition \ref{prop_green_identity} and (\ref{eq_gij_fr_dom_sq}), we have
		\begin{equation}\label{eq_dmin_dg}
			\scal{ \nabla^F  g^n_{i} }{ \nabla^F g^n_{j} }_{L^2(\Psi, F)} 
			=
			\scal{ \laplcomp_{\Psi}^{F} g^n_{i}  }{ g^n_{j} }_{L^2(\Psi, F)}.
		\end{equation}
		\par 
		Now, let's consider a vector space $V_k^{n} \subset \ccal^{\infty}_{0, vN}(\Psi, F)$, spanned by $g^{n}_{1}, \ldots, g^{n}_{k}$.
		By Proposition \ref{prop_ln_appr_easy}, applied for $\phi = 1$, we see that for $n$ big enough, we might choose $k_{n, i}$ big enough so that we have $\dim V_k^{n} = k$.
		We use the characterization of the eigenvalues of $\laplcomp_{\Psi}^{F}$ through Rayleigh quotient
		\begin{equation}\label{eq_rayleigh_cont}
			\lambda_{k} = \inf_{\substack{V \subset {\rm{Dom}}_{Fr}(\laplcomp_{\Psi}^{F}) } }  \sup_{f \in V} 
			\bigg\{ 
				\frac{\scal{\laplcomp_{\Psi}^{F} f}{f}_{L^2(\Psi, F)}}{\scal{f}{f}_{L^2(\Psi, F)}}
				 : \dim V = k			
			\bigg\}.
		\end{equation}
		In particular, by (\ref{eq_gij_fr_dom_sq}), we conclude
		\begin{equation}\label{eq_final_lin_pf_1}
			\lambda_{k} \leq \sup_{f \in V_k^{n}} 
				\frac{\scal{\laplcomp_{\Psi}^{F} f}{f}_{L^2(\Psi, F)}}{\scal{f}{f}_{L^2(\Psi, F)}}.
		\end{equation}
		However, by (\ref{eq_from_disc_to_ln}), (\ref{eq_dmin_ln_dmin_g_n}), (\ref{eq_dmin_dg}), we deduce that for any $k \in \nat$, we have
		\begin{equation}\label{eq_final_lin_pf_2}
			\liminf_{n \to \infty} \sup_{f \in V_k^{n}} 
				\frac{\scal{\laplcomp_{\Psi}^{F} f}{f}_{L^2(\Psi, F)}}{\scal{f}{f}_{L^2(\Psi, F)}}
			=
			\liminf_{n \to \infty} \lambda_{k}^{n}.
		\end{equation}
		From (\ref{eq_final_lin_pf_1}) and (\ref{eq_final_lin_pf_2}), we deduce Theorem  \ref{thm_bnd_fin_cont_below}.
		\end{proof}
		Now, to complete the proof of Theorem \ref{thm_bnd_fin_cont_below}, let's establish Propositions \ref{prop_ln_dom_min},  \ref{prop_dmin_equal} and \ref{prop_ln_appr_easy}.
		\begin{proof}[Proof of Proposition \ref{prop_ln_dom_min}.]
			First, take a cut-off function $\rho: \Psi \to [0, 1]$, which is equal to $1$ in the neighborhood of $\partial \Psi$ and which satisfies 
			\begin{equation}\label{eq_supp_rho}
				{\rm{supp}}(\rho) \subset B_{\Psi} \Big( \frac{1}{2n}, \partial \Psi \Big).
			\end{equation}
			For $f \in  {\rm{Map}}(V(\Psi_n), F_n)$, we decompose 
			\begin{equation}\label{eq_ln_dom_min_aux_1}
				L_n(f) = \rho L_n(f) + (1 - \rho) L_n(f).
			\end{equation}
			We will prove that $\rho L_n(f) \in {\rm{Dom}}_{\min}^{vN}(\nabla^F)$ and $(1 - \rho) L_n(f) \in {\rm{Dom}}_{\min}^{vN}(\nabla^F)$.
			\par 
			We start with $(1-\rho) L_n(f)$. 
			Take another cut-off function $\rho_1$, which is equal to $1$ near a small neighborhood ${\rm{Con}}(\Psi)$ and $0$ away from a neighborhood of it.
			Clearly, if one takes $\rho_1$ with very little support, then since $L_n(f)$ is constant over $U_n(P)$ for any point $P \in {\rm{Con}}(\Psi)$, there is a constant $C_P$ such that $\rho_1 (1-\rho) L_n(f)$ near $P$ is equal to $\rho_1 C_P$. 
			As a consequence, we have 
			\begin{equation}\label{eq_ln_dom_min_aux_11}
				\rho_1 (1-\rho) L_n(f) \in {\rm{Dom}}_{\min}^{vN}(\nabla^F).
			\end{equation}
			Now, let's take a smooth domain $W$ such that ${\rm{supp}} (1 - \rho_1) (1-\rho) \subset W$. 
			Then since $(1 - \rho_1) (1-\rho) L_n(f)$ is a Lipshitz function (it is a product of a smooth function and a piecewice linear function), by the density results (cf. \cite[Theorem 1.4.2.1]{Grisvard}) in Sobolev space $W^1_{\infty}$, we conclude that $(1 - \rho_1) (1 - \rho) L_n(f) \in {\rm{Dom}}_{\min}^{vN}(\nabla^F)$. 
			This, along with (\ref{eq_ln_dom_min_aux_11}) imply
			\begin{equation}\label{eq_ln_dom_min_aux_2}
				(1 - \rho) L_n(f) \in {\rm{Dom}}_{\min}^{vN}(\nabla^F).
			\end{equation}
			\par 
			Now let's study $\rho L_n(f)$. 
			Take a cut-off function $\rho_2$, which is equal to $1$ near a small neighborhood ${\rm{Ang}}(\Psi)$ and $0$ away from a neighborhood of it.
			Clearly, if one takes $\rho_2$ with very little support, then by the similar argument as for $\rho_1 (1-\rho) L_n(f)$, we have 
			\begin{equation}\label{eq_ln_dom_min_aux_22}
				\rho_2 \rho L_n(f) \in {\rm{Dom}}_{\min}^{vN}(\nabla^F).
			\end{equation}						
			\par Now, by (\ref{eq_supp_rho}), the support of $(1 - \rho_2) \rho L_n(f)$ is contained in a union of rectangles.
			But since over rectangles, by (\ref{eq_lin_rect_defn}), the function $L_n(f)$ varies only in one direction - parallel to the boundary, the inclusion $(1 - \rho_2) \rho L_n(f) \in {\rm{Dom}}_{\min}^{vN}(\nabla^F)$ follows from the fact that on an interval, a piecewise linear function $l$ can be approximated by smooth functions $g_n$ such that $g_n \to l$ and $d(g_n) \to h$ in $L^2(\Psi, F)$ for some $h \in L^2$ (which follows again from the density results for Sobolev space $W^1_{\infty}$, (cf.  \cite[Theorem 1.4.2.1]{Grisvard})).
			By this and (\ref{eq_ln_dom_min_aux_22}), we deduce
			\begin{equation}\label{eq_ln_dom_min_aux_3}
				\rho L_n(f) \in {\rm{Dom}}_{\min}^{vN}(\nabla^F),
			\end{equation}
			We conclude by (\ref{eq_ln_dom_min_aux_1}), (\ref{eq_ln_dom_min_aux_2}) and (\ref{eq_ln_dom_min_aux_3}). 
		\end{proof}
		\begin{proof}[Proof of Proposition \ref{prop_dmin_equal}.]
			The proof is a rather boring verification.
			We decompose 
			\begin{equation}\label{eq_dmin_eq_aux_1}
			\begin{aligned}
				\scal{ \nabla^F \big(L_n(f) \big) }{ \nabla^F \big( L_n(g) \big) }_{L^2(\Psi, F)}
				=
				&
				\\
				\sum_{P \in {\rm{Con}}(\Psi) \cup {\rm{Ang}}(\Psi)} & \int_{U_n(P)} \scal{ \nabla^F \big(L_n(f) \big) }{ \nabla^F \big( L_n(g) \big) }(z) dv_{\Psi}(z)
				\\
				+
				\sum_{R} & \int_{R} \scal{ \nabla^F \big(L_n(f) \big) }{ \nabla^F \big( L_n(g) \big) }(z) dv_{\Psi}(z)
				\\
				+
				\sum_{T} & \int_{T} \scal{ \nabla^F \big(L_n(f) \big) }{ \nabla^F \big( L_n(g) \big) }(z) dv_{\Psi}(z),
			\end{aligned}
			\end{equation}
			where $\sum_{R}$ means the sum over the rectangles  near $\partial \Psi$ from the definition of $L_n$ and $\sum_{T}$ means a sum over the triangles of the triangulation used in the definition of $L_n$, see Figure \ref{fig_ln_func}.
			\par First, by the construction, $L_n(f)$ is constant near $U_n(P)$ for $P \in {\rm{Con}}(\Psi) \cup {\rm{Ang}}(\Psi)$.	
			Thus 
			\begin{equation}\label{eq_dmin_eq_aux_2}
				\int_{U_n(P)} \scal{ \nabla^F \big(L_n(f) \big) }{ \nabla^F \big( L_n(g) \big) }(z) dv_{\Psi}(z) = 0.
			\end{equation}
			\par 
			Now, let $P, Q \in V(\Psi_n)$ correspond to a rectangle $R$ and $x, y$ be the coordinates from the construction of $L_n$.
			Normalize the coordinate $y$ so that it takes value $0$ on the boundary of rectangle corresponding to $\partial \Psi$ and value $\frac{1}{2}$ on the parallel boundary. 
			Then by (\ref{eq_lin_rect_defn}), the following holds
			\begin{equation}\label{eq_dmin_eq_aux_3}
			\begin{aligned}
				\int_{R} \scal{ \nabla^F \big(L_n(f) \big) }{ \nabla^F \big( L_n(g) \big) }&(z) dv_{\Psi}(z)
				\\
				=
				&
				\int_0^{\frac{1}{2}} \int_0^{1} 
				\scal{(f(P) - f(Q)) dx}{ (g(P) - g(Q)) dx}(z) dv_{\Psi}(z)
				\\
				=
				&
				\frac{1}{2}(f(P) - f(Q)) \cdot (g(P) - g(Q)).
			\end{aligned}
			\end{equation}
			\par Finally, let $P, Q, R \in V(\Psi_n)$ be the vertices of a triangle $T$, appearing in the definition of $L_n$.
			Let $x, y$ be the coordinates as in the definition of $L_n$.
			Assume $P$ and $Q$ share the same $x$-coordinate and $P$ and $R$ share the same $y$-coordinate.
			Then similarly to (\ref{eq_dmin_eq_aux_3}), we have
			\begin{equation}\label{eq_dmin_eq_aux_4}
			\begin{aligned}
				\int_{T} \scal{ \nabla^F & \big(L_n(f) \big) }{ \nabla^F \big( L_n(g) \big) }(z) dv_{\Psi}(z)
				\\						
				&
				=
				\int_T \scal{(f(P) - f(Q) ) dy}{(g(P) - g(Q))dy}(z) dv_{T}(z)
				\\
				&
				\phantom{= \,}
				+
				\int_T \scal{(f(P) - f(R) ) dx}{(g(P) - g(R))dx}(z) dv_{T}(z)
				\\
				&
				=
				\frac{1}{2}(f(P) - f(Q)) \cdot (g(P) - g(Q))
				+
				\frac{1}{2}(f(P) - f(R)) \cdot (g(P) - g(R)).
			\end{aligned}
			\end{equation}
			By (\ref{eq_scal_laplg_dg}) and (\ref{eq_dmin_eq_aux_1}) - (\ref{eq_dmin_eq_aux_4}), we conclude.
		\end{proof}
		\begin{proof}[Proof of Proposition \ref{prop_ln_appr_easy}.]
			Note that by Theorem \ref{thm_bnd_fin_cont_above}, for any $i \in \nat$, there is a constant $C_i > 0$ such that for any $n \in \nat^*$, we have
			\begin{equation}\label{eq_un_bnd_eigenv_up}
				\lambda_{i}^{n} < C_i.
			\end{equation}
			\par 
			By the construction of the linearization functional $L_n$, there is $C$, which depends only on the set $\angle({\rm{Con}}(\Psi) \cup {\rm{Ang}}(\Psi))$ and on $\norm{\phi}_{\ccal^{1}(\Psi, F)}$, such that for $P \in V(\Psi_n)$ and $x \in \Psi$, satisfying $\dist_{\Psi}(P, x) < \frac{1}{n}$, the following holds
			\begin{equation}\label{eq_ln_appr_triv}
				\Big| (\phi \cdot L_n(f_i^{n}))(x) - (\phi \cdot f_i^{n})(P) \Big| \leq C \max_{{\rm{dist}}_{\Psi_n}(P, Q) \leq 2} \big| f_n(P) - f_n(Q) \big|
				+
				\frac{C}{n} |f_i^{n}(P)|.
			\end{equation}
			From (\ref{eq_ln_appr_triv}), and the fact that our normalization is chosen so that the area of the initial tiles is $1$, we see that there is $C > 0$ such that for any $i, j \in \nat$, $n \in \nat^*$, we have
			\begin{equation}\label{eq_ln_eig_appr_aux_1}
			\begin{aligned}
				\Big| \scal{\phi L_n(f_i^{n}) }{L_n(f_j^{n}) }_{L^2(\Psi, F)}& - \frac{1}{n^2} \scal{\phi f_i^{n} }{f_j^{n}}_{L^2(\Psi_n)} \Big|
				\leq
				\frac{C}{n^3}\scal{\phi f_i^{n} }{f_j^{n}}_{L^2(\Psi_n)}
				\\
				&
				+
				\frac{C}{n^2}
				\sum_{P \in V(\Psi_n)} |f_j^{n}(P)| \cdot \max_{{\rm{dist}}_{\Psi_n}(P, Q) \leq 2} |f_i^{n}(P) - f_i^{n}(Q)|
				\\
				&
				+
				\frac{C}{n^2}
				\sum_{P \in V(\Psi_n)} |f_i^{n}(P)| \cdot \max_{{\rm{dist}}_{\Psi_n}(P, Q) \leq 2} |f_j^{n}(P) - f_j^{n}(Q)|.
			\end{aligned}
			\end{equation}
			Now, by Cauchy inequality, there is a constant $C$, which depends only on the set $\angle({\rm{Con}}(\Psi) \cup {\rm{Ang}}(\Psi))$ and on $\norm{\phi}_{\ccal^{1}(\Psi, F)}$, such that we have
			\begin{equation}\label{eq_ln_eig_appr_aux_2}
				\sum_{P \in V(\Psi_n)} |f_j^{n}(P)| \cdot  \max_{{\rm{dist}}_{\Psi_n}(P, Q) \leq 2}   \big| f_i^{n}(P) - f_i^{n}(Q) \big| 
				\leq 
				C
				\big\| f_j^{n} \big\|_{L^2(\Psi_n)} \cdot \big\| \nabla_{\Psi_n}^{F_n} f_i^{n} \big\| _{L^2(\Psi_n)}.
			\end{equation}
			However by (\ref{eq_un_bnd_eigenv_up}) and the bound on the norm of $f_i^{n}$, we see that
			\begin{equation}\label{eq_ln_eig_appr_aux_3}
				\big\| \nabla_{\Psi_n}^{F_n} f_i^{n} \big\| _{L^2(\Psi_n)}
				=
				\scal{\laplcomp_{\Psi_n}^{F_n} f_i^{n}}{f_i^{n}}_{L^2(\Psi_n)}^{1/2}
				\leq
				C_i.
			\end{equation}
			Also, from the bound on the norm of $f_i^{n}$, $f_j^{n}$ and Cauchy inequality, we have
			\begin{equation}\label{eq_ln_eig_appr_aux_4}
				\frac{C}{n^3}\scal{\phi f_i^{n} }{f_j^{n}}_{L^2(\Psi_n)} \leq \frac{C}{n} \norm{\phi}_{L^{\infty}(\Psi, F)}.
			\end{equation}
			From (\ref{eq_ln_eig_appr_aux_1}), (\ref{eq_ln_eig_appr_aux_2}), (\ref{eq_ln_eig_appr_aux_3}) and (\ref{eq_ln_eig_appr_aux_4}), we conclude.
		\end{proof}
		
	\subsection{Convergence of the eigenvectors, a proof of Theorem \ref{thm_eigvec_convergence}}\label{sect_pf_eigenvect}
	The main goal of this section is to prove Theorem \ref{thm_eigvec_convergence}.
	We will see that Theorem \ref{thm_eigvec_convergence} follows almost formally from the approximation theory we developed in Sections \ref{sect_eigvec_reg} and \ref{sect_harn_in}.
	\par 
	\begin{sloppypar}
	Assume that the eigenvalue $\lambda_i$, $i \in \nat^*$ of $\laplcomp_{\Psi}^{F}$ has multiplicity $m_i$.
	By Theorem \ref{thm_eigval_convergence}, we see that there is a series of eigenvalues $\lambda_{i, j}^{n}$, $j = 1, \ldots, m_i$ of $n^2 \cdot \laplcomp_{\Psi_n}^{F_n}$ (possibly equal, but in general not), converging to $\lambda_i$, as $n \to \infty$. Moreover, no other eigenvalue of $n^2 \cdot \laplcomp_{\Psi_n}^{F_n}$ come close to $\lambda_i$ asymptotically.
	We denote by $V_{i}^{n}$ the vector space spanned by the eigenvectors $f_{i, j}^{n}$ of $n^2 \cdot \laplcomp_{\Psi_n}^{F_n}$, corresponding to the eigenvalues $\lambda_{i, j}^{n}$, $j = 1, \ldots, m_i$, and let 
	\begin{equation}
		\pi_{V_{i}^{n}} : {\rm{Map}}(V(\Psi_n), F_n) \to V_{\lambda_i}^{n}
	\end{equation}
	be the orthogonal projection onto this space.
	Recall that the “linearization" $L_n : {\rm{Map}}(V(\Psi_n), F_n) \to L^2(\Psi, F)$ functional was defined in the beginning of the Section \ref{sect_harn_in}, and the restriction functional $R_n : \ccal^{0}(\Psi, F) \to {\rm{Map}}(V(\Psi_n), F_n)$ was defined in (\ref{eq_defn_restr_fun}).
	The main result of this section is the following
	\end{sloppypar}
	\begin{thm}\label{thm_eigvec_convergence_precised}
		For fixed $i, j \in \nat^*$, $j \leq m_i$, as $n \to \infty$, in $L^2(\Psi, F)$, we have
		\begin{equation}\label{eq_eigvec_convergence_precised}
			L_n( \pi_{V_{i}^{n}} (R_n f_{i, j})) \to f_{i, j},
		\end{equation}
		where $f_{i, j} \in L^2(\Psi, F)$, $j = 1, \ldots, m_i$ are mutually orthogonal eigenvectors of $\laplcomp_{\Psi}^{F}$, corresponding to the eigenvalue $\lambda_i$.
	\end{thm}
	\begin{rem}
		Clearly Theorem \ref{thm_eigvec_convergence_precised} implies Theorem \ref{thm_eigvec_convergence} by Proposition \ref{prop_ln_appr_easy} and (\ref{eq_lim_cont_linf}).
	\end{rem}
	\begin{proof}
		For simplicity of the presentation, we will suppose that the spectra of $n^2 \cdot \laplcomp_{\Psi_n}^{F_n}$ and $\laplcomp_{\Psi}^{F}$ are simple, i.e. the eigenvalues have multiplicity $1$. The proof of the general case remains verbatim, but the notation becomes way more difficult.
		We denote by $f_i$, $\norm{f_i}_{L^2(\Psi, F)} = 1$, the eigenvectors of $\laplcomp_{\Psi}^{F}$ corresponding to the eigenvectors $\lambda_i$, and by $f_i^{n}$, $\norm{f_i^{n}}_{L^2(\Psi_n, F_n)}^{2} = n^2$, the eigenvectors of $n^2 \cdot \laplcomp_{\Psi_n}^{F_n}$ corresponding to the eigenvectors $\lambda_i^{n}$, 
		\par 
		We decompose 
		\begin{equation}\label{eq_thm_vect_pr_dec}
			\begin{aligned}
			& 
			R_n f_{i} = 
			\alpha_i^{n} f_i^{n}
			+
			\sum_{l \neq i} \beta_{l}^{i, n} f_l^{n},
			\\
			&
			L_n f_{i}^{n} = 
			\gamma_i^{n} f_i
			+
			\sum_{l \neq i} \delta_{l}^{i, n} f_l,
			\end{aligned}
		\end{equation}
		for some $\alpha_i^{n}, \gamma_i^{n}$ and $\beta_l^{i, n}, \delta_l^{i, n}$, for $i, n \in \nat^*$.
		Then Theorem \ref{thm_eigvec_convergence_precised} would follow if, as $n \to \infty$, the following convergence holds
		\begin{equation}\label{eq_conv_eigv_main}
			\begin{aligned}
			& 
			\alpha_i^{n}  \to 1, \qquad
			\sum_{l \neq i} (\beta_{l}^{i, n})^2 \to 0.
			\\
			&
			\gamma_i^{n}  \to 1, \qquad
			\sum_{l \neq i} (\delta_{l}^{i, n})^2 \to 0.
			\end{aligned}
		\end{equation}
		Let's show that (\ref{eq_conv_eigv_main}) holds by induction on $i$.
		It clearly holds for $i = 0$ by Corollary \ref{cor_kernel_const}.
		Suppose we proved it for $i - 1$, let's show that it also holds for $i$.
		\par 
		Indeed, we know by (\ref{eq_lim_cont_linf}) that for any $j \leq i - 1$, as $n \to \infty$, the following convergence holds
		\begin{equation}\label{eq_conv_eigv_main_aux_11}
			\frac{1}{n^2} \scal{R_n f_{i}}{R_n f_{j}}_{L^2(\Psi_n, F_n)} \to \scal{f_{i}}{f_{j}}_{L^2(\Psi, F)} =  0.
		\end{equation}
		However, by (\ref{eq_thm_vect_pr_dec}), we have
		\begin{equation}\label{eq_conv_eigv_main_aux_1132}
			\frac{1}{n^2} \scal{R_n f_{i}}{R_n f_{j}}_{L^2(\Psi_n, F_n)}
			=
			\alpha_i^{n} \beta_i^{j, n} 
			+
			\alpha_j^{n} \beta_j^{i, n} 
			+
			\sum_{l \neq i, j} \beta_l^{i, n} \beta_l^{j, n}.
		\end{equation}
		From the induction hypothesis (\ref{eq_conv_eigv_main}), (\ref{eq_conv_eigv_main_aux_11}) and (\ref{eq_conv_eigv_main_aux_1132}), we conclude that for any $j \leq i -1$, as $n \to \infty$, the following limit holds
		\begin{equation}\label{eq_conv_eigv_main_aux_0}
			\beta_{j}^{i, n}  \to 0.
		\end{equation}
		Now, by Theorem \ref{thm_reg_main_limit_2}, as $n \to \infty$, the following convergence holds
		\begin{equation}\label{eq_conv_eigv_main_aux_1}
			\scal{\laplcomp_{\Psi_n}^{F_n} ( R_n f_{i})}{R_n f_{i}}_{L^2(\Psi_n, F_n)} \to \lambda_i.
		\end{equation}
		However, by (\ref{eq_thm_vect_pr_dec}), we have
		\begin{equation}\label{eq_conv_eigv_main_aux_2}
			\scal{\laplcomp_{\Psi_n}^{F_n} (R_n f_{i})}{R_n f_{i}}_{L^2(\Psi_n, F_n)}
			=
			(\alpha_i^{n})^2 \cdot \lambda_i^{n}
			+
			\sum_{l \neq i} (\beta_{l}^{i, n})^2 \lambda_l^{n}.
		\end{equation}
		Since $\lambda_l^{n}$, $l \geq j+1$ is asymptotically bigger than $\lambda_i^{n}$, as $n \to \infty$, by (\ref{eq_conv_eigv_main_aux_0}),  (\ref{eq_conv_eigv_main_aux_1}),  (\ref{eq_conv_eigv_main_aux_2}), we deduce the first statement of (\ref{eq_conv_eigv_main}). The second statement is shown analogically, where one has to use Theorem \ref{thm_ln_eig_appr} instead of Theorem \ref{thm_reg_main_limit_2} and Proposition \ref{prop_ln_appr_easy} instead of (\ref{eq_lim_cont_linf}).
	\end{proof}

\section{Harnack-type inequality for discrete Laplacian}\label{sect_harn_pf}
	The theory of harmonic functions on lattices has received considerable attention in the past; for example, see Duffin \cite{Duffin}, Kenyon \cite{Kenyon2002Invent}, {B\"ucking} \cite{Bucking}... Much less attention has been drawn to functions which are \textit{almost harmonic} in the sense that their Laplacian is not zero, but evaluated at some point, it is close to zero in comparison with the value of the function itself at this point.
	The main goal of this section is to give a proof Theorem \ref{thm_harn_type}, which is a statement in this direction.
	\par
	This section is organized as follows. In Section \ref{sec_pf_harn_modul}, we give a proof of Theorem \ref{thm_harn_type} modulo Theorem \ref{thm_harn_type_sup_bound}, which gives a uniform upper bound on the values of almost-harmonic function in terms of the square root of the logarithm of the number of vertices in a graph. 
	In Section \ref{sec_unif_bound_pf}, we give a proof of Theorem \ref{thm_harn_type_sup_bound} by assuming Theorem \ref{thm_harn_type_sup_bound_1}, which gives uniform constant bound on the values of almost-harmonic function away from the conical and angular singularities. 
	Finally, in Section \ref{sec_inter_est}, we prove Theorem \ref{thm_harn_type_sup_bound_1} by using discrete potential theory.
	
	\subsection{Asymptotic continuity of discrete eigenvectors, a proof of Theorem \ref{thm_harn_type}}\label{sec_pf_harn_modul}
	This section is devoted to the proof of Theorem \ref{thm_harn_type}. The following result is central here
	\begin{thm}\label{thm_harn_type_sup_bound}
		We fix $\lambda > 0$. Suppose that a sequence $f_n \in {\rm{Map}}(V(\Psi_n), F_n)$ satisfies the same assumptions as in Theorem \ref{thm_harn_type}.
		Then, as $n \to \infty$, the following limit holds
		\begin{equation}\label{eq_harn_type_sup_bound}
			\frac{1}{\sqrt{ \log n}} \max_{P \in V(\Psi_n)} \big| f_n(P) \big| \to 0.
		\end{equation}
	\end{thm}
	We prove Theorem \ref{thm_harn_type_sup_bound} in Section \ref{sec_unif_bound_pf}.
	Let's see how Theorem \ref{thm_harn_type_sup_bound} implies Theorem \ref{thm_harn_type}.
	\begin{proof}[Proof of Theorem \ref{thm_harn_type}.]
		First, by (\ref{eq_harn_bound_est}), we see that 
		\begin{equation}\label{eq_harn_1_aux_0}
			\scal{\laplcomp_{\Psi_n}^{F_n} f_n}{ f_n}_{L^2(\Psi_n, F_n)}
			\leq
			\lambda.
		\end{equation}
		By (\ref{eq_scal_laplg_dg}) and (\ref{eq_harn_1_aux_0}), we deduce that
		\begin{equation}\label{eq_harn_1_aux_00}
			\sum_{e \in E(\Psi_n)} (\nabla_{\Psi_n}^{F_n} f)(e)^2 
			\leq			
			\lambda.
		\end{equation}
		As a consequence, we obtain
		\begin{equation}\label{eq_harn_1_aux_000}
			\max_{e \in E(\Psi_n)} \big| (\nabla_{\Psi_n}^{F_n} f)(e) \big| \leq 
			\sqrt{\lambda}.
		\end{equation}
		\par 
		The main idea of the proof is to show that if there is an edge $e$ so that the difference of $f_n$ along $e$ is asymptotically bounded from below by a nonzero constant, then the neighboring edges of $e$ will contribute a lot to the sum in the left-hand side of (\ref{eq_harn_1_aux_00}) up to the point so that the bound (\ref{eq_harn_1_aux_00}) doesn't hold. This gives a contradiction to the initial assumption that (\ref{eq_harn_type_main_res}) doesn't hold.
		\par 
		More rigorously, suppose that the statement of Theorem \ref{thm_harn_type} is false. 
		Then we can find $c > 0$ and $e_n \in E(\Psi_n)$, $n \in \nat$ such that for any $n \in \nat^*$, up to choosing a subsequence, we have
		\begin{equation}\label{eq_harn_1_aux__1}
			\big| (\nabla_{\Psi_n}^{F_n} f_n)(e_n) \big| \geq c.
		\end{equation}
		\par 
		Denote by $\partial B_{\Psi_n}(r, P) \subset E(\Psi_n)$, $P \in V(\Psi_n)$ the subset defined as follows
		\begin{equation}\label{eq_defn_partial_set}
			\partial B_{\Psi_n}(r, P)
			=
			\Big\{
				e \in E(\Psi_n) : h(e) \in B_{\Psi_n}(r, P), t(e) \notin  B_{\Psi_n}(r, P)
			\Big\}.
		\end{equation}
		Let $P_n = h(e_n)$.		 
		For $e \in E(\Psi_n)$, we denote $e \in B_{\Psi_n}(r, P_n)$ if $h(e), t(e) \in B_{\Psi_n}(r, P_n)$.
		Then by (\ref{eq_scal_laplg_dg}), the following identity holds
		\begin{multline}\label{eq_harn_1_aux_1}
			\scal{\laplcomp_{\Psi_n}^{F_n} f_n}{f_n}_{L^2(B_{\Psi_n}(r, P_n), F_n)}
			=
			\sum_{e \in B_{\Psi_n}(r, P_n)} \big| (\nabla_{\Psi_n}^{F_n} f_n)(e) \big|^2 
			\\
			+ 
			\sum_{e \in \partial B_{\Psi_n}(r, P_n)} (\nabla_{\Psi_n}^{F_n} f_n)(e) \cdot  f_n(h(e)). 
		\end{multline}
		However, by (\ref{eq_harn_bound_est}), we see that 
		\begin{equation}\label{eq_harn_2_aux_2}
			\scal{\laplcomp_{\Psi_n}^{F_n} f_n}{ f_n}_{L^2(B_{\Psi_n}(r, P_n), F_n)}
			\leq
			\frac{\lambda}{n^2}
			\big\| f_n \big\|_{L^2(B_{\Psi_n}(r, P_n), F_n)}^{2}.
		\end{equation}
		Now, let $r \leq \sqrt{n}$, then by Theorem \ref{thm_harn_type_sup_bound}, there is $N \geq 0$ such that for $n \geq N$, we have
		\begin{equation}\label{eq_harn_2_aux_3}
			\big\| f_n \big\| _{L^2(B_{\Psi_n}(r, P_n), F_n)}^{2} \leq n^{3/2}.
		\end{equation}
		From (\ref{eq_harn_1_aux__1}), (\ref{eq_harn_1_aux_1}), (\ref{eq_harn_2_aux_2}) and (\ref{eq_harn_2_aux_3}), we deduce that for any $n \geq N$ and $r \leq \sqrt{n}$ we have
		\begin{equation}\label{eq_harn_2_aux_4}
			\sum_{e \in \partial B_{\Psi_n}(r, P_n)} |(\nabla_{\Psi_n}^{F_n} f_n)(e)| \cdot |f_n(h(e))|
			\geq
			c^2 - \frac{\lambda}{\sqrt{n}}.
		\end{equation}
		However, by Theorem \ref{thm_harn_type_sup_bound}, for any $\epsilon > 0$, there is $N$ such that for all $n \geq N$, we have
		\begin{equation}\label{eq_harn_2_aux_5}
			\sum_{e \in \partial B_{\Psi_n}(r, P_n)} |(\nabla_{\Psi_n}^{F_n} f_n)(e)| \cdot |f_n(h(e))|
			\leq
			\epsilon \sqrt{\log(n)}
			\sum_{e \in \partial B_{\Psi_n}(r, P_n)} | (\nabla_{\Psi_n}^{F_n} f_n)(e)|.
		\end{equation}
		Now, from trivial geometric considerations, there exists $C > 0$, which depends only on the sets $\angle({\rm{Con}}(\Psi))$ and $\angle({\rm{Ang}}(\Psi))$, such that for any $r \leq n$, we have
		\begin{equation}\label{eq_harn_2_aux_6}
			\# \partial B_{\Psi_n}(r, P) \leq C r.
		\end{equation}
		By (\ref{eq_harn_2_aux_6}) and mean inequality, we have
		\begin{equation}\label{eq_harn_2_aux_7}
			\Big( \sum_{e \in \partial B_{\Psi_n}(r, P_n)} \big| (\nabla_{\Psi_n}^{F_n} f_n)(e) \big| \Big)^2
			\leq
			Cr \sum_{e \in \partial B_{\Psi_n}(r, P_n)} \big| (\nabla_{\Psi_n}^{F_n} f_n)(e) \big|^2.
		\end{equation}
	 	From (\ref{eq_harn_2_aux_4}), (\ref{eq_harn_2_aux_5}) and (\ref{eq_harn_2_aux_7}), we deduce that for any $\epsilon > 0$, there is $N$ such that for all $n \geq N$ and $r \leq \sqrt{n}$, we have
	 	\begin{equation}\label{eq_harn_2_aux_8}
	 		\sum_{e \in \partial B_{\Psi_n}(r, P_n)} \big| (\nabla_{\Psi_n}^{F_n} f_n)(e) \big|^2
			\geq
			\frac{c^2}{2 \log(n)} \frac{1}{\epsilon^2 C r}.
	 	\end{equation}
	 	From (\ref{eq_harn_1_aux_00}) and (\ref{eq_harn_2_aux_8}), we see that for any $\epsilon > 0$, there is $N \geq 0$ such that for any $n \geq N$, we have
	 	\begin{equation}\label{eq_harn_2_aux_9}
	 		\sum_{r = 1}^{\sqrt{ n}} \frac{1}{r} \leq \epsilon \log(n),
	 	\end{equation}
	 	which is a nonsense. 
	 	Thus, the initial assumption (\ref{eq_harn_1_aux__1}) is false and Theorem \ref{thm_harn_type} holds.
	\end{proof}
\subsection{Uniform bound on discrete eigenvectors, a proof of Theorem \ref{thm_harn_type_sup_bound}}\label{sec_unif_bound_pf}
	The main goal of this section is to give a uniform bound on values of discrete eigenvectors, i.e. to prove Theorem \ref{thm_harn_type_sup_bound}.
	We note that all the statements of this section are local in nature, so we will suppose without loosing the generality that $(F, h^F, \nabla^F)$ is trivial of rank $1$.
	\par
	For $c > 0$, we define
		\begin{equation}
			V_c(\Psi_n) := \Big\{ 
				P \in V(\Psi_n) 
				:
				\dist_{\Psi}(P, {\rm{Con}}(\Psi) \cup {\rm{Ang}}(\Psi)) > c
			\Big\}.
		\end{equation}
	The following theorem says that in the interior of $\Psi$, the bound (\ref{eq_harn_type_sup_bound}) can be significantly improved.
	The proof of it uses the potential theory on $\integ^2$ and it is given Section \ref{sec_inter_est}. 
	\begin{thm}\label{thm_harn_type_sup_bound_1}
		We fix $\lambda > 0$.  
		For any $c > 0$, there is a constant $C(\lambda) > 0$, such that for any sequence $f_n \in {\rm{Map}}(V(\Psi_n), F_n)$  as in Theorem \ref{thm_harn_type}, we have
		\begin{equation}
			\max_{P \in V_c(\Psi_n)} \big| f_n(P) \big| \leq C(\lambda).
		\end{equation}
	\end{thm}
	Our proof of Theorem \ref{thm_harn_type_sup_bound} relies on Theorem \ref{thm_harn_type_sup_bound_1}, maximum principle and some ad-hoc constructions.
	The first construction is the construction of the discrete flow on the $n \times n$ subgraph inside of $\integ^2$, which has divergence $-1$ at one corner point and divergence $\frac{1}{n+1}$ on $n+1$ points, located at the distance $n$ from the corner point (see (\ref{eq_div_defn}) for the definition of divergence).
	This flow will be used to show that if the values of $|f_n|$, evaluated at $P_n \in V(\Psi_n)$ grow at least as $\sqrt{\log n}$, as $n \to \infty$, then there is another sequence $Q_n \in V(\Psi_n)$, satisfying $\dist_{\Psi_n}(P_n, Q_n) > cn$ for some $c > 0$, such that $|f_n(Q_n)|$  grow at least as $\sqrt{\log n}$, as $n \to \infty$ as well.
	\par 
	To be more precise, we define $E_{\integ^2}^{n}$ as follows: for $a, b \in \nat$, $a + b \leq n$, we put
	\begin{equation}\label{eq_harn_pf_aux_2}
	\begin{aligned}
		& E_{\integ^2}^{n}((a, b), (a+1, b)) = (a + 1) \Big( \frac{1}{a + b + 1} - \frac{1}{a + b + 2} \Big),
		\\
		& E_{\integ^2}^{n}((a, b), (a, b+1)) = \frac{1}{a + b + 2} - a \Big( \frac{1}{a + b + 1} - \frac{1}{a + b + 2} \Big).
	\end{aligned}
	\end{equation}
	And for other edges, $E_{\integ^2}^{n}$ is defined to be zero.  
	Then it is a matter of a tedious calculation, to see that for $Q \in V(\Psi_n)$, the following holds
	\begin{equation}\label{eq_harn_pf_aux_3}
		(d_{\integ^2}^{*} E_{\integ^2}^{n})(Q) = 
		\begin{cases} 
			\hfill -1, & \text{ for } Q = (0,0), \\
			\hfill (n + 1)^{-1}, & \text{ for $Q = (a, b) \in \nat^2$ with $a + b = n$,}  \\
			\hfill 0, & \text{ otherwise.}
 			\end{cases}
	\end{equation}
	Also there is $C > 0$ such that for any $n \in \nat$, we have
	\begin{equation}\label{eq_harn_pf_aux_6}
		\big\|  E_{\integ^2}^{n} \big\|_{L^2(\Psi_n)}
		\leq 
		\Big( 2 \sum_{i = 1}^{n} \frac{1}{i}\Big)^{1/2}
		\leq
		C \sqrt{\log(n)}.
	\end{equation}
	\par 
	The second construction is of a certain convex function on infinite cones and corners.
	\begin{proof}[Proof of Theorem \ref{thm_harn_type_sup_bound}.]
		The proof has two parts. In the first part, by using the flow $E_{\integ^2}^{n}$, considered above, we prove that there is a constant $C_0(\lambda)$ such that for any $n \in \nat^*$, the following bound holds
		\begin{equation}\label{eq_harn_pf_aux_1}
			\frac{1}{\sqrt{ \log n}} \max_{P \in V(\Psi_n)} \big| f_n(P) \big| \leq C_0(\lambda).
		\end{equation}
		In the second step, by using certain convex function, we “bootstrap" the estimate (\ref{eq_harn_pf_aux_1}) to get Theorem \ref{thm_harn_type_sup_bound}.
		\par We start by giving a proof of (\ref{eq_harn_pf_aux_1}).
		The main idea is to estimate the change of $f_n$ along a flow $E_{P}^{n}$, $P \in V(\Psi_n)$, constructed similarly to $E_{\integ^2}^{n}$, to relate the values of $ f_n$ near the conical points and angular singularities with the values of $f_n$ at $V_1(\Psi_n)$. By using the bound on the values of $f_n$ on $V_1(\Psi_n)$ from Theorem \ref{thm_harn_type_sup_bound_1}, we deduce (\ref{eq_harn_pf_aux_1}).
		\par More precisely, let $P \in V(\Psi_n)$. Choose a horizontal and vertical rays from $P$ such that the quadrant with side $1$ between them is isomorphic to a standard quadrant in $\real^2$ and all the points at the distance bigger than $1$ from $P$ in this quadrant lie in $V_1(\Psi_n)$. 
		This is always possible up to a multiplication of the initial metric by $4$.
		Then the part of the graph $\Psi_n$, which gets trapped inside of the quadrant of size $1$ is isomorphic to the part of $\integ^2$.
		Our flow $E_P^{n} \in {\rm{Map}}(E(\Psi_n), \comp)$ has support inside this quadrant and it coincides with $E_{\integ^2}^{n}$, normalized so that $P$ corresponds to $(0, 0)$.
		\par 
		As a consequence of (\ref{eq_harn_pf_aux_3}), we see that
		\begin{equation}\label{eq_harn_pf_aux_5}
			\scal{d_{\Psi_n} f_n}{E_P^{n}}_{L^2(\Psi_n)}
			=
			\scal{f_n}{d_{\Psi_n}^{*} E_P^{n}}_{L^2(\Psi_n)}
			=
			-f_n(P)
			+
			\frac{1}{n + 1} \sum_{Q \in {\rm{supp}}(d_{\Psi_n}^{*} E_P^{n}) \setminus P} f_n(Q).
		\end{equation}
		By Cauchy inequality, (\ref{eq_harn_1_aux_0}) and  (\ref{eq_harn_pf_aux_6}), we deduce
		\begin{equation}\label{eq_harn_pf_aux_66}
			\scal{d_{\Psi_n} f_n}{E_P^{n}}_{L^2(\Psi_n)}
			\leq 
			C \sqrt{\lambda} \sqrt{\log(n)}.
		\end{equation}
		By (\ref{eq_harn_pf_aux_5}), (\ref{eq_harn_pf_aux_66}), we see that
		\begin{equation}\label{eq_harn_pf_aux_7}
			| f_n(P) | \leq \frac{1}{n + 1} \sum_{Q \in {\rm{supp}}(d_{\Psi_n}^{*} E_P^{n}) \setminus P} |f_n(Q)| + C \sqrt{\lambda} \sqrt{\log(n)}.
		\end{equation}
		By (\ref{eq_harn_pf_aux_7}), the fact that ${\rm{supp}}(d_{\Psi_n}^{*} E_P^{n}) \setminus P \in V_1(\Psi, F)$, $\# {\rm{supp}}(d_{\Psi_n}^{*} E_P^{n}) = n + 2$ and Theorem \ref{thm_harn_type_sup_bound_1}, we see that the bound (\ref{eq_harn_pf_aux_1}) holds for $C_0(\lambda) = C \sqrt{\lambda} + C(\lambda)$, where $C(\lambda)$ is from Theorem \ref{thm_harn_type_sup_bound_1}.
		\par Now, let's show that the bound (\ref{eq_harn_pf_aux_1}) can be “bootstraped" to Theorem \ref{thm_harn_type_sup_bound}.
		In fact, by Theorem \ref{thm_harn_type_sup_bound_1}, we only need to prove Theorem \ref{thm_harn_type_sup_bound} in the neighborhood of ${\rm{Con}}(\Psi) \cup {\rm{Ang}}(\Psi)$.
		\par 
		Let's fix $P \in {\rm{Con}}(\Psi) \cup {\rm{Ang}}(\Psi)$.
		Recall that the set $V_n(P)$ was defined in (\ref{eq_defn_vp}).
		Consider a function $h_P$, given by
		\begin{equation}\label{eq_harn_pf_aux_77}
			h_P(Q) = \dist_{\Psi_n}(Q, V_n(P))^2, \qquad Q \in V(\Psi_n).
		\end{equation}
		We verify that inside of $B_{\Psi_n}(n, V_n(P))$, we have
		\begin{equation}\label{eq_harn_pf_aux_8}
			\laplcomp_{\Psi_n}^{F_n} h_P \leq -1.
		\end{equation}
		Let $C_0(\lambda)$ be from (\ref{eq_harn_pf_aux_1}).
		Consider a function
		\begin{equation}\label{eq_harn_pf_aux_9}
			g_n = f_n + 2 C_0(\lambda) \sqrt{\log(n)} \frac{\lambda}{n^2} h_P.
		\end{equation}
		Then by (\ref{eq_harn_bound_est}), (\ref{eq_harn_pf_aux_1}) and (\ref{eq_harn_pf_aux_8}), we see that over $B_{\Psi_n}(n, V_n(P))$, we have
		\begin{equation}
			\laplcomp_{\Psi_n}^{F_n} g_n < 0.
		\end{equation}		 
		Thus, by the maximum principle, the maximum of $g_n$ on $B_{\Psi_n}(r,  V_n(P))$, $r \leq n$ is attained on $S_{\Psi_n}(r,  V_n(P))$.
		Thus, for any $\epsilon$, satisfying $0 < \epsilon < 1$, we have
		\begin{multline}\label{eq_harn_pf_aux_10}
			\max_{Q \in B_{\Psi_n}(\epsilon n,  V_n(P))} \Big( f_n(Q) + 2C_0(\lambda) \sqrt{\log(n)} \frac{\lambda}{n^2} h_P(Q)
			\Big)
			\\
			\leq
			\max_{R \in S_{\Psi_n}(\epsilon n,  V_n(P))} f_n(R)
			+
			 2C_0(\lambda) \epsilon^2 \lambda \sqrt{\log(n)}.
		\end{multline}
		From Theorem \ref{thm_harn_type_sup_bound_1}, we see that there is $C_{0}(\lambda, \epsilon) > 0$ such that
		\begin{equation}\label{eq_harn_pf_aux_11}
			\max_{\substack{R \notin B_{\Psi_n}(\epsilon n, P) \\ P \in {\rm{Con}}(\Psi) \cup {\rm{Ang}}(\Psi)}} | f_n(R) | \leq C_{0}(\lambda, \epsilon).
		\end{equation}
		In particular, we have
		\begin{equation}
			 \max_{R \in S_{\Psi_n}(\epsilon n,  V_n(P))} f_n(R) \leq C_{0}(\lambda, \epsilon).
		\end{equation}
		From (\ref{eq_harn_pf_aux_77}), (\ref{eq_harn_pf_aux_10}) and (\ref{eq_harn_pf_aux_11}), we conclude that for any $0 < \epsilon < 1$ and $n \in \nat^*$, we have
		\begin{equation}\label{eq_harn_pf_aux_121}
			\max_{Q \in B_{\Psi_n}(\epsilon n,  V_n(P))} f_n(Q) 
			\leq
			2 C_0 (\lambda) \epsilon^2 \sqrt{\log(n)} + C_0(\lambda, \epsilon).
		\end{equation}
		Similarly, we can bound $f_n$ from below.
		By this, (\ref{eq_harn_pf_aux_11}) and (\ref{eq_harn_pf_aux_121}), we conclude.
	\end{proof}		
	\par 
	\subsection{Interior bounds on almost harmonic functions, a proof of Theorem \ref{thm_harn_type_sup_bound_1}}\label{sec_inter_est}
	The main goal of this section is to get interior bounds on almost harmonic functions and to prove Theorem \ref{thm_harn_type_sup_bound_1}.
	The proof of this theorem relies heavily on the discrete potential theory introduced by Duffin, \cite{Duffin}, in dimension 3 and developed by Kenyon, \cite{Kenyon2002Invent}, in dimension 2.
	\par 
	Since Theorem \ref{thm_harn_type_sup_bound_1} is a statement about the value of $f_n$ at points far away from conical and angle singularities, it would follow from the following two theorems.
	The first one corresponds to the bound on the interior points of $\Psi$.
	\begin{thm}\label{thm_harn_type_sup_bound_2}
		We fix $\lambda > 0$. 
		There is a constant $C(\lambda) > 0$, such that for any sequence $f_n \in {\rm{Map}}(V(\integ^2), \comp)$ satisfying $\norm{f_n}_{L^2(B_{\integ^2}(n, 0))}^{2} \leq n^2$, and 
		\begin{equation}\label{eq_f_n_z_sq_ass}
			\big| \laplcomp_{\integ^2} f_n \big| \leq \frac{\lambda}{n^2} |f_n|, \qquad \text{over} \qquad \qquad B_{\integ^2}(n, 0),
		\end{equation}
		we have the following bound
		\begin{equation}
			\big| f_n(0) \big| \leq C(\lambda).
		\end{equation}
	\end{thm}
	And the second one corresponds to the bound on the points near the boundary $\partial \Psi$.
	\begin{thm}\label{thm_harn_type_sup_bound_3}
		We fix $\lambda > 0$. 
		There is a constant $C(\lambda) > 0$, such that for any sequences $P_n \in \nat \times \integ$ and $f_n \in {\rm{Map}}(V(\nat \times \integ), \comp)$ satisfying $\norm{f_n}_{L^2(B_{\nat \times \integ}(n, P_n))}^{2} \leq n^2$ and
		\begin{equation}\label{eq_f_n_nat_z_ass}
			\big| \laplcomp_{\nat \times \integ} f_n \big| \leq \frac{\lambda}{n^2} |f_n|, \qquad \text{over} \qquad B_{\nat \times \integ}(n, P_n),
		\end{equation}
		we have the following bound
		\begin{equation}
			\big| f_n(P_n) \big| \leq C(\lambda).
		\end{equation}
	\end{thm}
	The proofs of those statements are similar and use Green functions on balls in $\integ^2$ and $\nat \times \integ$.
	\par 
	Let's start with Theorem \ref{thm_harn_type_sup_bound_2}. Consider $\integ^2$ with its nearest-neighbor graph structure and a standard injection $\integ^2 \hookrightarrow \comp$.
	For simplicity, for $n \in \nat$, we denote
	\begin{equation}
		B_{\comp}^{\integ^2}(n, 0) := \{ z \in \integ^2  : |z| \leq n\},
	\end{equation}
	where $|z|$ is the modulus of $z$ (and not the distance between $0$ and $z$ in $\integ^2$).
	\par 
	None of the Propositions \ref{prop_gr_fun}, \ref{prop_ken_green}, \ref{prop_bound_gr_integ2} are new.
	We, nevertheless, recall their proofs for completeness and to use them later in the corresponding statements on $\nat \times \integ$.
	\begin{prop}[{Duffin \cite[p. 242]{Duffin}}]\label{prop_gr_fun}
		There exists a function $G_n^{\integ^2} : \integ^2 \to \real$ satisfying
		\begin{equation}\label{eq_prop_gr_fun}
		\begin{aligned}
			&
			{\rm{supp}}(G_n^{\integ^2}) \subset B_{\comp}^{\integ^2}(n, 0),
			\\
			&
			(\laplcomp_{\integ^2} G_n^{\integ^2})(z) 
			=
			\begin{cases} 
				\hfill 1, & \text{ for } z = 0, \\
				\hfill 0, & \text{ for  } z \in B_{\comp}^{\integ^2}(n, 0) \setminus \{0\}.
 			\end{cases}
		\end{aligned}
		\end{equation}
	\end{prop}
	\begin{proof}
		We prove Proposition \ref{prop_gr_fun} iteratively.
		Let $G_n^{1} : \integ^2 \to \real$ be a function, equal to $1$ at $0$ and $0$ otherwise.
		Construct $G_n^{2} : \integ^2 \to \real$ by 
		\begin{equation}
			G_n^{i}(z) 
			=
			\begin{cases} 
				\hfill 1, & \text{ for } z = 0, \\
				\hfill \frac{1}{4} \big( \sum_{(z_1, z) \in E(\integ^2)} G_n^{i - 1}(z_1)  \big), & \text{ for $z \in B_{\comp}^{\integ^2}(n, 0) \setminus \{0\}$,}\\
				\hfill 0, & \text{ otherwise.}
 			\end{cases}
		\end{equation}
		Continuing this process gives a sequence of functions $G_n^{i}$, $i \in \nat$.
		By induction, it's easy to see that this sequence is increasing and bounded above by $1$.
		Thus, the sequence $G_n^{i}$, $i \in \nat^*$, converges to a function $G_n^{\infty}$, satisfying $0 \leq G_n^{\infty} \leq 1$.
		Clearly, $G_n^{\infty}$ is harmonic on $B_{\integ^2}(n, 0) \setminus \{0\}$.
		If $G_n^{\infty}$ were harmonic at $0$, then by the maximum principle it follows that $G_n^{\infty}$ takes maximal value at $S_{\integ^2}(n + 1, 0)$.
		But it is clearly false, since it takes zero values there.
		Thus, we see that $(\laplcomp_{\integ^2} G_n^{\infty})(0) > 0$, and so the function $G_n^{\integ^2}$, defined by
		\begin{equation}
			G_n^{\integ^2} := \frac{1}{(\laplcomp_{\integ^2} G_n^{\infty})(0)} G_n^{\infty},
		\end{equation}
		satisfies the hypothesizes (\ref{eq_prop_gr_fun}).
	\end{proof}
	\begin{prop}[{Duffin \cite[Theorem 2, \S 3]{Duffin} in dimension $3$, Kenyon \cite[Theorem 7.3]{Kenyon2002Invent} in dimension $2$}]\label{prop_ken_green}
		There exists a unique function $G_{\integ^2} : \integ^2 \to \real$ such that for any $z, w \in \integ^2$, we have
		\begin{equation}
		\begin{aligned}
			&
			(\laplcomp_{\integ^2, z}) G_{\integ^2}(z, w) = 
			\begin{cases} 
				\hfill 1, & \text{ for } z = w, \\
				\hfill 0, & \text{ otherwise,}
 			\end{cases}
			\\
			&
			G_{\integ^2}(z, z) = 0,
			\\
			&
			G_{\integ^2}(z, w) = O(\log |z - w|), \text{ as } |z - w| \to \infty,
		\end{aligned}
		\end{equation}	
		where by $\laplcomp_{\integ^2, z}$ we mean a Laplacian evaluated with respect to the $z$-coordinate.
		Moreover, as $|z - w| \to \infty$, such a function satisfies the following asymptotic expansion
		\begin{equation}\label{eq_ken_bound_g_integ}
			G_{\integ^2}(z, w) = -\frac{1}{2 \pi} \log |z - w| - \frac{\gamma_{EM}}{2 \pi} + O\Big( \frac{1}{|z - w|} \Big),
		\end{equation}			 
		where $\gamma_{EM}$ is the Euler-Mascheroni constant.
	\end{prop}
	\begin{rem}
		The precise value of $\gamma_{EM}$ is not important in what follows. It is, however, of vital importance that the expansion (\ref{eq_ken_bound_g_integ}) is up to the term $O(1/ |z - w|)$.
	\end{rem}
	For $n \in \nat$, we denote 
	\begin{equation}\label{eq_defn_sphere}
		S_{\comp}^{\integ^2}(n, 0) := B_{\comp}^{\integ^2}(n, 0) \setminus B_{\comp}^{\integ^2}(n-1, 0).
	\end{equation}
	Using Proposition \ref{prop_ken_green}, we can now get
	\begin{prop}[{Duffin \cite[Lemma 2]{Duffin} in dimension $3$, B\"ucking \cite[Proposition A.3]{Bucking} in dimension $2$}]\label{prop_bound_gr_integ2}
		Let $G_n$ be a function from Proposition \ref{prop_gr_fun}.
		There is a constant $C > 0$ such that for any $n \in \nat^*$
		\begin{equation}
			\max_{z \in S_{\comp}^{\integ^2}(n, 0)} G_n(z) \leq \frac{C}{n}.
		\end{equation}
	\end{prop}
	\begin{proof}
		Consider a function $f_n : \integ^2 \to \real$, given by
		\begin{equation}\label{eq_bound_g_n_fun_aux_1}
			f_n(z) 
			=
			G_n^{\integ^2}(z)
			-
			G_{\integ^2}(0, z)
			-
			\frac{1}{2 \pi} \log(n) 
			-
			\frac{\gamma_{EM}}{2 \pi}.
		\end{equation}
		From (\ref{eq_prop_gr_fun}) and (\ref{eq_ken_bound_g_integ}), we see that there is $C > 0$ such that for any $n \in \nat^*$, we have
		\begin{equation}\label{eq_bound_g_n_fun_aux_2}
			\max_{z \in S_{\comp}^{\integ^2}(n + 1, 0)} |f_n(z)|
			\leq 
			\frac{C}{n}.
		\end{equation}
		However, by Propositions \ref{prop_gr_fun}, \ref{prop_ken_green}, the function $f_n$ satisfies
		\begin{equation}\label{eq_bound_g_n_fun_aux_3}
			(\laplcomp_{\integ^2} f_n)(z) = 0, \qquad \text{for } z \in B_{\integ^2}(n, 0).
		\end{equation}
		Thus, by the maximal principle and (\ref{eq_bound_g_n_fun_aux_2}), we conclude that 
		\begin{equation}\label{eq_bound_g_n_fun_aux_4}
			\max_{z \in S_{\comp}^{\integ^2}(n, 0)} |f_n(z)|
			\leq 
			\frac{C}{n}.
		\end{equation}
		By (\ref{eq_ken_bound_g_integ}) and (\ref{eq_bound_g_n_fun_aux_4}), we conclude. 
	\end{proof}
	As a consequence of Proposition \ref{prop_bound_gr_integ2}, we get the following
	\begin{prop}\label{prop_l2_bound_green_fun}
		There is a constant $C > 0$ such that for any $n \in \nat^*$, $z \in B_{\comp}^{\integ^2}(n, 0)$, we have
		\begin{equation}\label{eq_l2_bound_green_fun}
		\begin{aligned}
			&
			G_n^{\integ^2}(z) \leq C \log\Big( \frac{n + 2}{|z| + 1} \Big),
			\\
			&
			\| G_n^{\integ^2} \|_{L^2(\integ^2)} \leq C n.
		\end{aligned}
		\end{equation}
	\end{prop}
	\begin{proof}
		We begin with a proof of the first bound  of (\ref{eq_l2_bound_green_fun}).
		Consider $f_n : \integ^2 \to \real$, defined by
		\begin{equation}
			f_n = G_{n+1}^{\integ^2} - G_{n}^{\integ^2}.
		\end{equation}
		By (\ref{eq_prop_gr_fun}), the function $f_n$ is harmonic in $B_{\comp}^{\integ^2}(n, 0)$.
		Thus, by the maximal principle, it attains its maximum on $S_{\comp}^{\integ^2}(n + 1, 0)$.
		By this and Proposition \ref{prop_bound_gr_integ2}, we see that there is $C > 0$ such that for any $n \in \nat^*$, we have
		\begin{equation}
			|f_n| \leq \frac{C}{n}.
		\end{equation}
		From this, we conclude by induction that
		\begin{equation}\label{eq_g_n_part_harm}
			G_n^{\integ^2}(z) \leq \sum_{i = |z|}^{n} \frac{C}{i}.
		\end{equation}
		By standard bounds on the partial harmonic sums, we obtain from (\ref{eq_g_n_part_harm}) the first estimate in (\ref{eq_l2_bound_green_fun}).
		Let's prove that the second estimate in  (\ref{eq_l2_bound_green_fun}) follows from the first.
		\par 
		In fact, by the first bound in (\ref{eq_l2_bound_green_fun}), we see that there is a constant $C > 0$, which depends only on the sets $\angle({\rm{Con}}(\Psi))$ and $\angle({\rm{Ang}}(\Psi))$, such that for any $n \in \nat^*$, the following holds
		\begin{equation}
			\| G_n^{\integ^2} \|_{L^2(\integ^2)}^{2} \leq C \sum_{i = 1}^{n} \log\Big( \frac{n+1}{i} \Big)^2 i.
		\end{equation}
		It is only now left to prove that there is a constant $C > 0$ such that for any $n \in \nat$, we have
		\begin{equation}\label{eq_l2_bound_green_fun_aux_01}
			A_n := \sum_{i = 1}^{n} \log\Big( \frac{n+1}{i} \Big)^2 i \leq C n^2.
		\end{equation}
		To show (\ref{eq_l2_bound_green_fun_aux_01}), it is enough to establish that there is a constant $C > 0$ such that for any $n \in \nat^*$:
		\begin{equation}\label{eq_l2_bound_green_fun_aux_001}
			B_n := A_n - A_{n - 1} \leq C n.
		\end{equation}
		We trivially have
		\begin{equation}\label{eq_l2_bound_green_fun_aux_1}
			B_n = \log\Big( \frac{n+1}{n} \Big)^2 n
			+
			\sum_{i = 1}^{n-1} \bigg( \log\Big( \frac{n+1}{i} \Big)^2 - \log\Big( \frac{n}{i} \Big)^2 \bigg)   i.
		\end{equation}
		Now, clearly, there is a constant $C > 0$ such that for any $n \in \nat^*$, we have
		\begin{equation}\label{eq_l2_bound_green_fun_aux_2}
		\begin{aligned}
			&
			\log\Big( \frac{n+1}{i} \Big)^2 - \log\Big( \frac{n}{i} \Big)^2 \leq \frac{C}{n} \log\Big( \frac{n}{i} \Big),
			\\
			&
			\log\Big( \frac{n+1}{n} \Big)^2 n \leq \frac{C}{n}.
		\end{aligned}
		\end{equation}
		By (\ref{eq_l2_bound_green_fun_aux_1}) and (\ref{eq_l2_bound_green_fun_aux_2}), to prove (\ref{eq_l2_bound_green_fun_aux_001}), it is enough to prove that there is a constant $C > 0$ such that for any $n \in \nat^*$, we have
		\begin{equation}\label{eq_l2_bound_green_fun_aux_22}
			C_n := \sum_{i = 1}^{n} \log\Big( \frac{n+1}{i} \Big)   i \leq Cn^2.
		\end{equation}
		To show (\ref{eq_l2_bound_green_fun_aux_22}), it is obviously enough to establish that there is a constant $C > 0$ such that for any $n \in \nat^*$, we have
		\begin{equation}\label{eq_l2_bound_green_fun_aux_3}
			D_n := C_n - C_{n - 1} \leq C n.
		\end{equation}
		But we have the following identity
		\begin{equation}\label{eq_l2_bound_green_fun_aux_4}
			D_n = 
			 \log\Big( \frac{n+1}{n} \Big)   n
			 +
			\sum_{i = 1}^{n-1} \log\Big( \frac{n+1}{n} \Big)   i,
		\end{equation}
		from which, along with the second bound from (\ref{eq_l2_bound_green_fun_aux_2}), we deduce (\ref{eq_l2_bound_green_fun_aux_3}).
	\end{proof}
	Finally we are ready to give
	\begin{proof}[Proof of Theorem \ref{thm_harn_type_sup_bound_2}.]
		Similarly to (\ref{eq_defn_partial_set}), we define $\partial B_{\comp}^{\integ^2}(r, 0) \subset E(\integ^2)$. 
		By (\ref{eq_scal_laplg_dg}), we see that for any $r \in \nat^*$, $r \leq n$, the following discrete analogue of Green identity holds
		\begin{multline}\label{eq_green_sup_bound_2_aux_1}
			\scal{\laplcomp_{\integ^2} f_n}{G_r^{\integ^2}}_{L^2(B_{\comp}^{\integ^2}(r, 0))} 
			- 
			\scal{f_n}{\laplcomp_{\integ^2}G_r^{\integ^2}}_{L^2(B_{\comp}^{\integ^2}(r, 0))}
			\\
			=
			\sum_{e \in \partial B_{\comp}^{\integ^2}(r, 0)}
			\Big(
			f_n(h(e)) G_r^{\integ^2}(t(e)) - f_n(h(e)) G_r^{\integ^2}(t(e))
			\Big).
		\end{multline}
		However, by Proposition \ref{prop_l2_bound_green_fun} and (\ref{eq_f_n_z_sq_ass}), we deduce that there is a constant $C > 0$ such that for any $n \in \nat^*$, $r \in \nat^*$, $r \leq n$, we have
		\begin{equation}\label{eq_green_sup_bound_2_aux_2}
			\scal{\laplcomp_{\integ^2} f_n}{G_r^{\integ^2}}_{L^2(B_{\comp}^{\integ^2}(r, 0))} 
			\leq
			C \lambda.
		\end{equation}
		By (\ref{eq_prop_gr_fun}), we have
		\begin{equation}\label{eq_green_sup_bound_2_aux_3}
			\scal{f_n}{\laplcomp_{\integ^2}G_r^{\integ^2}}_{L^2(B_{\comp}^{\integ^2}(r, 0))}
			=
			f_n(0).
		\end{equation}
		From Proposition \ref{prop_bound_gr_integ2} and (\ref{eq_prop_gr_fun}), we see that there is a constant $C > 0$ such that for any $n \in \nat^*$, $r \in \nat^*$, $n/2 \leq r \leq n$, we have
		\begin{equation}\label{eq_green_sup_bound_2_aux_4}
			\Big|
				\sum_{e \in \partial B_{\comp}^{\integ^2}(r, 0)}
				\Big(
				f_n(h(e)) G_r^{\integ^2}(t(e)) - f_n(h(e)) G_r^{\integ^2}(t(e))
				\Big)
			\Big|
			\leq
			\frac{C}{n}
			\sum_{e \in \partial B_{\comp}^{\integ^2}(r, 0)}
			|f_n(h(e))|
		\end{equation}
		By combining (\ref{eq_green_sup_bound_2_aux_1})-(\ref{eq_green_sup_bound_2_aux_4}), we see that for any $n \in \nat^*$, $r \in \nat^*$, $n/2 \leq r \leq n$, we have
		\begin{equation}\label{eq_green_sup_bound_2_aux_5}
			|f_n(0)| 
			\leq 		
			C \lambda
			+	
			\frac{C}{n}
			\sum_{e \in \partial B_{\comp}^{\integ^2}(r, 0)}
			|f_n(h(e))|.
		\end{equation}
		By taking average of (\ref{eq_green_sup_bound_2_aux_5}) for $r \in \nat^*$, $n/2 \leq r \leq n$, we have
		\begin{equation}\label{eq_green_sup_bound_2_aux_6}
			|f_n(0)| 
			\leq 		
			C \lambda
			+	
			\frac{2 C}{n^2}
			\sum_{P \in B_{\comp}^{\integ^2}(n, 0)}
			|f_n(P)|.
		\end{equation}
		By mean inequality and the assumptions on the norm of $f_n$, for some constant $C > 0$, we have
		\begin{equation}\label{eq_green_sup_bound_2_aux_7}
			\bigg(			
			\frac{1}{n^2}
			\sum_{P \in B_{\comp}^{\integ^2}(n, 0)}
			|f_n(P)|
			\bigg)^2	
			\leq
			\frac{C}{n^2}
			\sum_{P \in B_{\comp}^{\integ^2}(n, 0)}
			|f_n(P)|^2
			\leq
			C.
		\end{equation}
		By (\ref{eq_green_sup_bound_2_aux_6}) and (\ref{eq_green_sup_bound_2_aux_7}), we conclude.
	\end{proof}
	Now, let's prove Theorem \ref{thm_harn_type_sup_bound_3}.
	The proof is parallel to Theorem \ref{thm_harn_type_sup_bound_2} modulo minor modifications.
	\par We denote by $\hh := \{ z \in \comp : \Im(z) \geq 0 \}$ the upper half-plane.
	We will implicitly use the following (non-standard) injection
	\begin{equation}\label{eq_nat_integ_inj}
		\nat \times \integ \hookrightarrow \hh, \qquad (a, b) \mapsto b + \imun \Big(a + \frac{1}{2} \Big).
	\end{equation}	 
	The shifting by $\frac{\imun}{2}$ will be used Proposition \ref{prop_ken_green_halfpla}.
	For $P \in \nat \times \integ$, $r \in \nat$, we denote 
	\begin{equation}
		QB(r, P)
		=
		\Big\{
			z \in \hh : 
			\big|
				(z - P) (z - \overline{P})
			\big|
			\leq 
			r^2
		\Big\}.
	\end{equation}
	We stress our that $z - P, z - \overline{P} \in \comp$ should be interpreted through (\ref{eq_nat_integ_inj}).
	The sets $\partial QB(r, P_n) \in E(\nat \times \integ)$, $QS(r, P_n) \in V(\nat \times \integ)$ are defined similarly to (\ref{eq_defn_partial_set}) and (\ref{eq_defn_sphere}).
	We have the following analogue of Proposition \ref{prop_gr_fun}.
	\begin{prop}\label{prop_gr_fun_half_plane}
		For any $P \in \nat \times \integ$, $r \in \nat^*$, there is a function $G_{P, r}^{\nat \times \integ} : \nat \times \integ \to \real$, satisfying
		\begin{equation}\label{eq_prop_gr_fun_half_plane}
		\begin{aligned}
			&
			{\rm{supp}}(G_{P, r}^{\nat \times \integ}) \subset QB(r, P),
			\\
			&
			(\laplcomp_{\nat \times \integ} G_{P, r}^{\nat \times \integ})(z) 
			=
			\begin{cases} 
				\hfill 1, & \text{ for } z = P, \\
				\hfill 0, & \text{ for  } z \in QB(r, P) \setminus \{P\}.
 			\end{cases}
		\end{aligned}
		\end{equation}
	\end{prop}
	\begin{proof}
		The proof is a formal repetition of the proof of Proposition \ref{prop_gr_fun}.
	\end{proof}
	We have the following weak analogue of Proposition \ref{prop_ken_green}.
	\begin{prop}\label{prop_ken_green_halfpla}
		There exists a function $G_{\nat \times \integ} : \nat \times \integ \to \real$ such that for any $z, w \in \nat \times \integ$:
		\begin{equation}\label{eq_ken_bound_g_integ_halfpla0}
			\laplcomp_{\nat \times \integ, z} G_{\nat \times \integ}(z, w) = 
			\begin{cases} 
				\hfill 1, & \text{ for } z = w, \\
				\hfill 0, & \text{ otherwise,}
 			\end{cases}
		\end{equation}	
		and the following asymptotic expansion, as $|z - w| \to \infty$, holds
		\begin{equation}\label{eq_ken_bound_g_integ_halfpla}
			G_{\nat \times \integ}(z, w) = -\frac{1}{2 \pi} \Big( \log |z - w| + \log |z - \overline{w}| \Big)  - \frac{\gamma_{EM}}{\pi} + O\Big( \frac{1}{|z - w|} \Big),
		\end{equation}			 
		where $\gamma_{EM}$ is the Euler-Mascheroni constant.
	\end{prop}
	\begin{proof}
		It suffices to take 
		\begin{equation}
			G_{\nat \times \integ}(z, w) 
			:=
			G_{\integ^2}(z - w, 0) 
			+
			G_{\integ^2}(z - \overline{w}, 0),
		\end{equation}
		where the points $z - w, z - \overline{w} \in \integ^2$ should be evaluated through the injection (\ref{eq_nat_integ_inj}).  The shifting by $\frac{\imun}{2}$ assures that the resulting function $G_{\nat \times \integ}(z, w)$ satisfies (\ref{eq_ken_bound_g_integ_halfpla0}).
	\end{proof}
	Now, by repeating line in line the proofs of Propositions \ref{prop_bound_gr_integ2}, \ref{prop_l2_bound_green_fun}, only by replacing each use of Proposition \ref{prop_ken_green} by Proposition \ref{prop_ken_green_halfpla}, we get the following two propositions 
	\begin{prop}\label{prop_bound_gr_integ2_hp}
		Let $P \in \nat \times \integ$, $r \geq 4 \dist_{\nat \times \integ}(P, \{ 0 \} \times \integ)$, and $G_{r, P}$ be as in Proposition \ref{prop_gr_fun_half_plane}.
		There is a constant $C > 0$ such that for any $n \in \nat^*$, we have
		\begin{equation}
			\max_{z \in QS(r, P)} G_{r, P}(z) \leq \frac{C}{r}.
		\end{equation}
	\end{prop}
	\begin{prop}\label{prop_l2_bound_green_fun_hp}
		Let $P \in \nat \times \integ$, $r \geq 4 \dist_{\nat \times \integ}(P, \{ 0 \} \times \integ)$, and $G_{r, P}$ be as in Proposition \ref{prop_gr_fun_half_plane}.
		There is a constant $C > 0$ such that for any $z \in QB(r, P)$, the following bounds hold
		\begin{equation}\label{eq_l2_bound_green_fun_hp}
		\begin{aligned}
			&
			G_{r, P}(z) \leq C \log\Big( \frac{r + 2}{|P-z| + 1} \Big),
			\\
			&
			\| G_{r, P} \|_{L^2(\nat \times \integ)} \leq C r.
		\end{aligned}
		\end{equation}
	\end{prop}
	By repeating formally the proof of Theorem \ref{thm_harn_type_sup_bound_2}, only replacing the use of Propositions \ref{prop_bound_gr_integ2}, \ref{prop_l2_bound_green_fun} by Propositions \ref{prop_bound_gr_integ2_hp}, \ref{prop_l2_bound_green_fun_hp}, we get Theorem \ref{thm_harn_type_sup_bound_3}.

		\bibliographystyle{abbrv}

\Addresses

\end{document}

%% file: packages.tex
\usepackage{times}

\usepackage{mathrsfs}
\usepackage{subfig}

\usepackage[shortlabels]{enumitem}
\usepackage[utf8]{inputenc}
\usepackage{commath}
\usepackage{amsthm}
\usepackage{wrapfig}
\usepackage{amsmath}
\usepackage{amssymb}
\usepackage{cite}
\usepackage{tikz}
\usepackage{setspace}

\usepackage{bigints}
\usepackage{comment}
\usepackage{mathtools}
\usepackage{mathrsfs}
\usepackage{fancyhdr}

\allowdisplaybreaks[4]

\numberwithin{equation}{section}
\usepackage{etoolbox}
\patchcmd{\thebibliography}
  {\settowidth}
  {\setlength{\itemsep}{0pt plus -10pt}\settowidth}
  {}{}
\apptocmd{\thebibliography}
  {
  }
  {}{}

%% file: theorem.tex
\makeatletter
\newtheorem*{rep@theorem}{\rep@title}
\newcommand{\newreptheorem}[2]{%
\newenvironment{rep#1}[1]{%
 \def\rep@title{#2 \ref{##1}}%
 \begin{rep@theorem}}%
 {\end{rep@theorem}}}
\makeatother

\theoremstyle{theorem}

\newreptheorem{theorem}{Theorem}
\newtheorem{thm}{Theorem}[section]
\newtheorem*{thm*}{Theorem}
\theoremstyle{definition}
\newtheorem{prop}[thm]{Proposition}
\newtheorem*{prop*}{Proposition}

\newtheorem{cor}[thm]{Corollary}
\newtheorem*{cor*}{Corollary}
\theoremstyle{remark}
\newtheorem{rem}[thm]{Remark}

%% file: title.tex
\title{Finite difference method on flat surfaces \\
 with a flat unitary vector bundle.} 

\author
{Siarhei Finski
}

\date{}

%% file: paddings.tex
\usepackage[%
    left=1in,%
    right=1in,%
    top=1.1in,%
    bottom=0.8in,%
    paperheight=11in,%
    paperwidth=8.5in%
]{geometry}

%% file: commands.tex
\newcommand{\imun} {\sqrt{-1}}

\newcommand{\comp}{\mathbb{C}}
\newcommand{\real}{\mathbb{R}}

\newcommand{\nat}{\mathbb{N}}
\newcommand{\integ}{\mathbb{Z}}

\newcommand{\spec}{{\rm{Spec}}}
\newcommand{\dist}{{\rm{dist}}}

\newcommand{\ccal}{\mathscr{C}}

\newcommand{\laplcomp}{\Delta}

\newcommand{\rk}[1]{{\rm{rk}} ( #1 )}

\renewcommand{\Im}{\operatorname{Im}}
\newcommand{\scal}[2]{\big< #1, #2 \big>}

\newcommand{\hh}{\mathbb{H}}